\documentclass[aos, preprint]{imsart}

\usepackage{fancyhdr}
\usepackage{amsmath}
\usepackage{amsthm}
\usepackage{amssymb,amsfonts}
\usepackage{threeparttable}
\usepackage[mathscr]{eucal}
\usepackage[dvips]{graphicx}
\usepackage{graphics}
\usepackage{geometry}
\usepackage{times}
\usepackage[dvipsnames]{xcolor}
\usepackage{natbib}
\usepackage{subfigure}
\usepackage{multirow}
\usepackage[normalem]{ulem}
\setattribute{journal}{name}{}

\RequirePackage[OT1]{fontenc}
\RequirePackage[colorlinks,citecolor=blue,urlcolor=blue]{hyperref}

\newtheorem{theorem}{Theorem}
\newtheorem{lemma}{Lemma}
\newtheorem{proposition}{Proposition}

\newtheorem{remark}{Remark}

\def\argmax{\mathop{\rm arg\, max}}
\def\argmin{\mathop{\rm arg\, min}}

\newcommand{\bel}{\begin{eqnarray}\label}
\newcommand{\eel}{\end{eqnarray}}
\newcommand{\bes}{\begin{eqnarray*}}
\newcommand{\ees}{\end{eqnarray*}}
\newcommand{\bei}{\begin{itemize}}
\newcommand{\eei}{\end{itemize}}
\newcommand{\beiftnt}{\begin{itemize}\footnotesize}

\def\benu{\begin{enumerate}}
\def\eenu{\end{enumerate}}

\def\argmax{\mathop{\rm arg\, max}}
\def\argmin{\mathop{\rm arg\, min}}
\def\real{{\mathbb{R}}}
\def\R{{\real}}

\def\E{{\mathbb{E}}}

\def\N{{\mathbb{N}}}
\def\P{{\mathbb{P}}}

\def\complex{\mathop{{\rm I}\kern-.58em\hbox{\rm C}}\nolimits}
\def\pa{\partial}

\def\Var{\hbox{\rm Var}}

\def\mathbold{\boldsymbol} 

\def\ba{\mathbold{a}}

\def\bb{\mathbold{b}}

\def\bc{\mathbold{c}}

\def\bfe{\mathbold{e}}

\def\bff{\mathbold{f}}\def\fhat{\widehat{f}}
\def\hbf{{\widehat{\bff}}}
\def\scrF{{\mathscr F}}

\def\bg{\mathbold{g}}
\def\hbg{{\widehat{\bg}}}
\def\scrG{{\mathscr G}}

\def\Htil{{\widetilde H}}

\def\bi{\mathbold{i}}

\def\bj{\mathbold{j}}

\def\bk{\mathbold{k}}

\def\scrK{{\mathscr K}}

\def\scrL{{\mathscr L}}

\def\bm{\mathbold{m}}

\def\scrM{{\mathscr M}}

\def\bn{\mathbold{n}}\def\ntil{\widetilde{n}}

\def\bs{\mathbold{s}}

\def\bt{\mathbold{t}}

\def\scrT{{\mathscr T}}

\def\bu{\mathbold{u}}

\def\scrU{{\mathscr U}}

\def\bv{\mathbold{v}}

\def\bw{\mathbold{w}}

\def\bx{\mathbold{x}}

\def\by{\mathbold{y}}\def\ybar{{\overline y}}

\def\Ybar{{\overline Y}}

\def\veps{\varepsilon}

\def\btheta{\mathbold{\theta}}\def\htheta{\widehat{\theta}}

\def\hbtheta{{\widehat{\btheta}}}

\def\hmu{\widehat{\mu}}

\def\fbar{{\overline f}}

\def\obu{{\overline{\bu}}}
\def\ubu{{\underline{\bu}}}
\def\ou{{\overline{u}}}
\def\uu{{\underline{u}}}
\def\argsup{\mathop{\rm arg\, sup}}

\begin{document}

\begin{frontmatter}
\title{Isotonic Regression in Multi-Dimensional Spaces and Graphs}

\begin{aug}
\author{\fnms{Hang} \snm{Deng}\ead[label=e1]{hdeng@stat.rutgers.edu}},
\and
\author{\fnms{Cun-Hui} \snm{Zhang}\thanksref{t1}\ead[label=e2]{czhang@stat.rutgers.edu}}

\thankstext{t1}{Partially supported by NSF grants DMS-1513378, IIS-1407939, DMS-1721495, IIS-1741390 and CCF-1934924.}

\affiliation{Department of Statistics, Rutger University}

\end{aug}

\begin{abstract}
In this paper we study minimax and adaptation rates in general isotonic regression. For uniform deterministic and random designs in $[0,1]^d$ with $d\ge 2$ and $N(0,1)$ noise, the minimax rate for the $\ell_2$ risk is known to be bounded from below by $n^{-1/d}$ when the unknown mean function $f$ is non-decreasing and its range is bounded by a constant, while the least squares estimator (LSE) is known to nearly achieve the minimax rate up to a factor $(\log n)^\gamma$ where $n$ is the sample size, $\gamma = 4$ in the lattice design and $\gamma = \max\{9/2, (d^2+d+1)/2 \}$ in the random design. Moreover, the LSE is known to achieve the adaptation rate $(K/n)^{-2/d}\{1\vee \log(n/K)\}^{2\gamma}$ when $f$ is piecewise constant on $K$ hyper-rectangles in a partition of $[0,1]^d$.

Due to the minimax theorem, the LSE is identical on every design point to both the max-min and min-max estimators over all upper and lower sets containing the design point. This motivates our consideration of estimators which lie in-between the max-min and min-max estimators over possibly smaller classes of upper and lower sets, including a subclass of block estimators. Under a $q$-th moment condition on the noise, we develop $\ell_q$ risk bounds for such general estimators for isotonic regression on graphs. For uniform deterministic and random designs in $[0,1]^d$ with $d\ge 3$, our $\ell_2$ risk bound for the block estimator matches the minimax rate $n^{-1/d}$ when the range of $f$ is bounded and achieves the near parametric adaptation rate $(K/n)\{1\vee\log(n/K)\}^{d}$ when $f$ is $K$-piecewise constant. Furthermore, the block estimator possesses the following oracle property in variable selection: When $f$ depends on only a subset $S$ of variables, the $\ell_2$ risk of the block estimator automatically achieves up to a poly-logarithmic factor the minimax rate based on the oracular knowledge of $S$.
\end{abstract}

\end{frontmatter}


\section{Introduction}

Let $G=(V,E)$ be a directed graph with vertex set $V$ and edge set $E$.
For $\ba$ and $\bb$ in $V$, we say that $\ba$ is a descendant of $\bb$ if $E$ contains a chain of edges from $\bv_j$ to $\bv_{j+1}$ such that $\bb=\bv_0$ and $\ba=\bv_m$ for some finite $m\ge 0$.
We write $\ba \preceq \bb$ if $\ba = \bb$ or $\ba$ is a descendant of $\bb$.
A function $f: V\to\R$ is non-decreasing on the graph $G$ if $f(\ba)\le f(\bb)$ whenever $\ba \preceq \bb$.
Let $\scrF$ be the class of all non-decreasing functions on $G$.
In isotonic regression, we observe $\bx_i\in V$ and $y_i\in \R$ satisfying
\bel{iso-reg}
y_i = f(\bx_i) + \veps_i,\ i=1,\ldots, n,\ \hbox{for some } f\in \scrF,
\eel
where $\veps_1,\ldots,\veps_n$ are independent noise variables with $\E\,\veps_i = 0$ and $\Var(\veps_i) \le \sigma^2$
given the (deterministic or random) design points $\{\bx_i\}$.
Note that we allow $|V|>n$.

An interesting special case of (\ref{iso-reg}) is the multiple isotonic regression
where $V \subset \R^d$ is a subset of a certain Euclidean space of dimension $d$,
and for $\ba = (a_1,\ldots,a_d)^T\in \R^d$ and
$\bb = (b_1,\ldots, b_d)^T\in \R^d$, $\ba \preceq \bb$ iff
$a_j\le b_j$ for all $1\le j\le d$.
In this case, $\scrF$ is the class of all non-decreasing functions on $V$.

Let $\bff_n = (f(\bx_1),\ldots,f(\bx_n))^T$ and
$\hbf_n = (\fhat_n(\bx_1),\ldots,\fhat_n(\bx_n))^T$ for any estimator $\fhat_n$ of $f$.
We are interested in the estimation of $f$ under 
the (normalized) $\ell_q$ risk
\bel{empirical-risk}
R_q(\hbf_n, \bff_n) = \frac{1}{n}\E \Big\|\hbf_n - \bff_n\Big\|_q^q = \frac{1}{n} \sum_{i=1}^n \E \Big| \fhat_n(\bx_i) - f(\bx_i)\Big|^q.
\eel
In this case, a specification of $\hbf_n$ is sufficient for the definition of $\fhat_n$.
For multiple isotonic regression with random design in $V\subseteq \R^d$,
we are also interested in the $L_q$ risk
\bel{integrated-risk}
R^*_q(\fhat_n, f) = \E \|\fhat_n - f\|_{L_q(V)}^q
= \E \int_V \Big|\fhat_n(\bx) - f(\bx)\Big|^q d\bx.
\eel

The literature of univariate isotonic regression ($d=1$)
encompasses at least the past six decades; See for example
\cite{brunk1955}, \cite{ayer1955empirical}, \cite{grenander1956theory}, \cite{rao1969estimation}, \cite{groeneboom1984estimating}, \cite{vandegeer1990, vandegeer1993}, \cite{donoho1990gelfand}, \cite{birge1993rates}, \cite{woodroofe1993penalized}, \cite{wang1996l2risk}, \cite{durot2007}, \cite{durot2008monotone}, and \cite{yang2017contraction} among many others for some key developments.
The least squares estimator (LSE), say $\fhat_n^{(lse)}$, has been the focus of this literature.
We describe in some detail here existing results on minimax and adaptation rates
as they are directly related to our study.
For any $a<b$, the $\ell_q$ risk of the LSE in the interval $[a,b]$ is bounded by
\bel{risk-bd-univ}
\E \sum_{a\le x_i\le b}\Big|\fhat_n^{(lse)}(x_i)-f(x_i)\Big|^q
\le C_q \sigma^q \bigg\{n_{a,b}\bigg(\frac{\Delta_{a,b}(\bff_n/\sigma)}{n_{a,b}}\wedge 1\bigg)^{q/3}
+\sum_{j=1}^{n_{a,b}} j^{-q/2} \bigg\},
\eel
where
$\Delta_{a,b}(\bff_n/\sigma) = \max_{a\le x_i<x_j\le b}\{f(x_j)-f(x_i)\}/\sigma$
is the range-to-noise ratio for the mean
vector $\bff_n$ in $[a,b]$, $n_{a,b}=\#\{j: a\le x_j\le b\}$ is the number of design points
in the interval, and $C_q$ is a constant depending on $q$ only.
This result can be found in \cite{meyer2000}
for $n_{a,b}=n$, $q=2$ and $\veps_i\sim N(0,\sigma^2)$,
and in \cite{zhang2002} for general $a<b$ and $1\le q < 3$
under a $(q\vee 2)$-th moment condition on $\veps_i$.
For $\Delta_{-\infty,\infty}(\bff_n/\sigma) \le \Delta^*_n \asymp 1$,
(\ref{risk-bd-univ}) yields the cube-root rate $\sigma^q(\Delta^*_n/n)^{q/3}$
for the LSE in terms of the $\ell_q$ risk in (\ref{empirical-risk}).
By summing over the risk bound (\ref{risk-bd-univ}) over $K$ intervals $[a_k,b_k]$ with
$\Delta_{a_k,b_k}(\bff_n/\sigma)=0$,
the LSE can be seen to achieve the near parametric adaptation rate $(K/n)\{1\vee \log(n/K)\}$
in the mean squared risk when the unknown $f$ is piecewise constant on the $K$ intervals
and $x_i\in\cup_{k=1}^K[a_k,b_k]$ for all $i\le n$.
This adaptation rate was explicitly given in \cite{chatterjee2015risk}.
However, \cite{gao2017minimax} proved that the sharp adaptation rate
in the mean squared risk, achieved by a penalized LSE, is $(K/n)\log\log(16n/K)$ in the
piecewise constant case. Moreover, by summing over the risk bound (\ref{risk-bd-univ})
over a growing number of disjoint intervals, the LSE has been shown to converge
faster than the cube root rate when the measure $f(dx)$ is singular to the Lebesgue measure \citep{zhang2002}.

Compared with the rich literature on univariate isotonic regression,
our understanding of the multiple isotonic regression,
i.e. $V\subset \R^d$ with $d>1$, is quite limited.
A major difficulty is that the design points are typically only partially ordered.
Univariate risk bounds can be directly applied to linearly ordered paths in $V$,
but this typically does not yield a nearly minimax rate.
However, significant advances have been made recently
on the minimax and adaptation rates for the LSE.
For $n_1\times\cdots\times n_d$ lattice designs with $n=\prod_{j=1}^d n_j$,
the LSE provides
\bel{lse-risk-bd}
R_2(\hbf_n^{(lse)},\bff_n)
\le C_d\sigma^2\bigg\{ \Delta(\bff_n/\sigma) n^{-1/d}(\log n)^\gamma
+ n^{-2/d}(\log n)^{2\gamma}\bigg\}
\eel
in certain settings, where $\Delta(\bff_n/\sigma) = \max_{1\le i<j\le n}|f(\bx_i)-f(\bx_j)|/\sigma$ is the range-to-noise ratio
of the mean over the design points.
For Gaussian $\veps_i$ and $n_1 = \cdots = n_d$,
the minimax rate is bounded from below by
\bel{lower-bd-univ}
\inf_{\hbf_n}\sup_{\Delta(\bff_n/\sigma) \le \Delta^*_n} R_2(\hbf_n,\bff_n)
\ge \sigma^2 \min\Big\{1, C_0 n^{-1/d}\Delta^*_n\Big\}.
\eel
Moreover, when $f$ is piecewise constant on
$K$ hyper-rectangles in a partition of the lattice,
\bel{lse-adapt}
R_2(\hbf_n^{(lse)},\bff_n) \le C_d\sigma^2(K/n)^{2/d}\{1\vee \log(n/K)\}^{2\gamma}.
\eel
For $d=2$ and Gaussian noise,
\cite{chatterjee2018} proved the above mean squared risk bounds with $\gamma=4$.
Thus, up to a logarithmic factor, the LSE is nearly rate minimax
for a wide range of $\Delta^*_n$ and also nearly adaptive
to the parametric rate $\sigma^2 K/n$ when
$f$ is piecewise constant on $K$ rectangles.
\cite{han2017isotonic} extended the results of \cite{chatterjee2018} from $d=2$ to $d>2$ under the conditions $n_1=\cdots =n_d$ and $\Delta(\bff_n/\sigma) \le \Delta^*_n=1$ in (\ref{lse-risk-bd}) and (\ref{lower-bd-univ}), and also proved parallel results for random designs with a larger
$\gamma = \max\{9/2, (d^2+d+1)/2\}$.
However, there is still a gap of a poly-logarithmic factor between such
upper and lower minimax bounds for $d\ge 2$,
and it is still unclear from (\ref{lse-adapt}) the feasibility
of near adaptation to the parametric rate $\sigma^2 K/n$ for $d\ge 3$
when $f$ is piecewise constant on $K$ hyper-rectangles.

We have also seen some progresses in adaptive estimation to variable selection
in isotonic regression on lattices with $\max_{j\le d} n_j\le C_d n^{1/d}$.
When the unknown mean function depends on only a \emph{known} subset of $s$ variables, say $f(\bx) = f_S(\bx_S)$ where $\bx_S=(x_j, j\in S)^T$ with $|S|=s$,
one may use the LSE, say $\fhat_{n,S}^{(lse)}$, based on the average of $y_i$ given $\bx_S$
to attain
\bel{oracle-selection}
R_2(\hbf_{n,S}^{(lse)},\bff_n)
\le \begin{cases}
C_d\sigma_S^2\Big[\Delta(\bff_n/\sigma_S) n^{-1/d}(\log n)^\gamma
+ n^{-2/d}(\log n)^{2\gamma}\Big], & s\ge 2,
\cr C_d\sigma_S^2\Big[\{(\Delta(\bff_n/\sigma_S)n^{-1/d})\wedge 1\}^{2/3}
+ n^{-1/d}\log n\Big], & s=1,
\end{cases}
\eel
with $\sigma_S^2 = \sigma^2/\prod_{j\not\in S}n_j\le C_d \sigma^2/n^{1-s/d}$,
which would match the minimax rate for Gaussian $\veps_i$ for a proper
range of $\Delta(\bff_n/\sigma_S)$ as we discussed in the previous paragraph.
For unknown $S$ with $d\ge 2$ and
$\Delta(\bff_n/\sigma) \le 1 = \sigma$,
\cite{han2017isotonic} proved that the LSE $\fhat_n^{(lse)}$
for the general $f$ automatically achieves the rate
$n^{-4/(3d)}(\log n)^{16/3}$ for $s=d-1$ and $n^{-2/d}(\log n)^8$ for $s\le d-2$.
As $\Delta(\bff_n/\sigma_S)\asymp n^{(d-s)/(2d)}$ in their setting,
(\ref{oracle-selection}) would yield the rates
$n^{-(d-s)/(2d)-1/d}$ for $s\ge 2$ and
$n^{-(d-1)/d-(3-d)_+/(3d)}$ for $s = 1$ up to a logarithmic factor. These oracle minimax rates nearly match the adaptation rates in \cite{han2017isotonic} for
$d-s=2$ or $(d,s)=(2,1)$, but not for other configurations of $(d,s)$.

We consider isotonic regression on directed graphs, i.e.
with general domain $V$ in (\ref{iso-reg}), including $V\subset \R^d$ as a special case.
In this general setting, \cite{robertson1988order} proved the following minimax formula
for the LSE on the design points:
\bel{minimax-formula}
\fhat_n^{(lse)}(\bx) = \max_{U\ni \bx}\min_{L\ni \bx}\ybar_{U\cap L}
= \min_{L\ni \bx}\max_{U\ni \bx}\ybar_{U\cap L}
\eel
for $\bx = \bx_i$, $i=1,\ldots,n$, where the maximum is taken over all upper sets $U$
containing $\bx$, the minimum over all lower sets $L$ containing $\bx$,
and $\ybar_A$ is the
average of the observed $y_i$ over $\bx_i\in A$ for any $A\subseteq V$.
As the high complexity of the upper and lower sets for $d\ge 2$
could be the culprit behind the
possible suboptimal performance of the LSE in convergence and adaptation rates,
we consider a class of block estimators involving rectangular upper and lower sets.
As the minimax theorem no longer holds in this setting in general, the block estimator, say $\fhat_n^{(block)}(\bx)$, is defined as any estimator in-between the following max-min and min-max estimators,
\bel{max-min-min-max}
\fhat_n^{(max-min)}(\bx) &=& \max_{\bu\preceq \bx, n_{\bu,*}>0}\ \min_{\bx\preceq\bv, n_{\bu,\bv}>0}
\ybar_{[\bu,\bv]},\quad\forall\ \bx\in V, 
\cr \fhat_n^{(min-max)}(\bx) &=&\min_{\bx\preceq\bv, n_{*,\bv}>0}\ \max_{\bu\preceq \bx, n_{\bu,\bv}>0}
\ybar_{[\bu,\bv]},\quad\forall\ \bx\in V, 
\eel
where $[\bu,\bv]=\{\bx: \bu\preceq\bx\preceq\bv\}$, $n_{\bu,\bv}=\#\{i\le n: \bx_i\in [\bu,\bv]\}$,
$n_{\bu,*}=\#\{i\le n: \bu\preceq \bx_i\}$ and $n_{*,\bv}=\#\{i\le n: \bx_i\preceq \bv\}$.
The idea of replacing the general level sets $U\cap L$ by rectangular blocks
$[\bu,\bv]$ is not new as 
a preliminary version of the block estimator in the case of $V=[0,1]^d$ was considered in \cite{fokianos2017integrated}.
Some more delicate details of different versions of the block estimator are discussed in Section 2.

We derive in Section 3 a general
$\ell_q$ risk bound for the above block estimator on graphs.
For $n_1\times\cdots\times n_d$ lattice designs with $d\ge 2$,
our general risk bound yields
\bel{risk-bd-block}
R_2\Big(\hbf_n^{(block)},\bff_n\Big)
\le C_d \sigma^2\min\Big\{1,\Delta(\bff_n/\sigma) n^{-1/d}(\log n)^{I\{d=2\}}+ n^{-1}(\log n)^d\Big\}
\eel
when $\max_{j\le d} n_j\le C_d n^{1/d}$, compared with \eqref{lse-risk-bd} and \eqref{lower-bd-univ}, and the adaptation rate
\bel{adaptation-block}
R_2\Big(\hbf_n^{(block)},\bff_n\Big) \le C_d\sigma^2(K/n)\{1\vee \log(n/K)\}^d
\eel
when the true $f$ is non-decreasing and piecewise constant on $K$ hyper-rectangles, compared with \eqref{lse-adapt}. 

We also explore the phase transition of the risk bounds, both the minimax lower bound and the upper risk bound for the block estimator, by presenting them using its effective dimension $s$ in the sense that the risk bound only depends on the largest $s$ $n_j$'s.
For example, when $n_1 \ge n_2 \ge \cdots \ge n_d$ and $n_2^{3/2}/n_1^{1/2} \le \Delta(\bff_n/\sigma)$, we show that the risk bound for the block estimator in $d$-dimensional isotonic regression with $n$ design points is almost no different from that in univariate isotonic regression with $n_1$ design points.
This phase transition, captured by effective dimension, proved for $d=2$ in \cite{chatterjee2018}, is new for $d>2$. 

Moreover, perhaps more interestingly,
we prove that when the unknown $f$ depends on an unknown set of $s$ variables,
the block estimator achieves near adaptation to the oracle selection
in the sense that for $\Delta(\bff_n/\sigma)\le \Delta^*_n$, 
\bel{adapt-selection}
&& R_2(\hbf_{n}^{(block)},\bff_n)
\\ \nonumber &\le& \begin{cases}
C_d \sigma_S^2\min\Big[(\log n)^{d-s}, \Delta^*_n n^{(d-s-2)/(2d)}(\log n)^{I\{s=2\}}
+ n^{-s/d}(\log n)^d\Big], & s\ge 2,
\cr C_d\sigma_S^2
\min\Big[(\log n)^{d-1}, \big(\Delta^*_n n^{(d-s-2)/(2d)}\big)^{2/3} + n^{-1/d}(\log n)^d\Big], & s=1,
\end{cases}
\eel
with $\sigma_S^2 = \sigma^2/\prod_{j\not\in S}n_j\le C_d \sigma^2/n^{1-s/d}$,
while the oracle minimax rate with the knowledge of $S$ is bounded from below by
\bel{minimax-selection}
&& \inf_{\hbf_n}\sup_{\bff_n}\Big\{R_2(\hbf_n,\bff_n):
\bff_n \in\scrF_n,\, f(\bx)=f_S(\bx_S),\Delta(\bff_n/\sigma)\le\Delta^*_n\Big\}
\\ \nonumber &\ge& \begin{cases}
C_d \sigma^2n^{-1+s/d}\min\Big[1,\Delta^*_n n^{(d-s-2)/(2d)}\Big], & s\ge 2,
\cr C_d\sigma^2n^{-1+1/d}
\min\Big[1, \big(\Delta^*_n n^{(d-3)/(2d)}\big)^{2/3}\Big], & s=1,
\end{cases}
\eel
where $\scrF_n = \{\bff_n: f\in\scrF\}$. 

Let ${\overline\bff}^*_n$ be the noiseless version of the block estimator. 
When the isotonic regression model is misspecified in the sense of having a non-monotone 
regression function, we prove that the error bounds discussed above still hold if 
${\overline\bff}^*_n$ is treated as the estimation target; 
\eqref{risk-bd-block}, \eqref{adaptation-block} and \eqref{adapt-selection} are valid 
with $\bff_n$ replaced by ${\overline\bff}^*_n$ when $f\not\in\scrF$ in \eqref{iso-reg}. 
However, such results are of a less ideal form compared with the 
existing oracle inequalities for the LSE under misspecified monotonicity assumption   
\citep{chatterjee2015risk, bellec2018, chatterjee2018, han2017isotonic}. 

We summarize our main results as follows.
In terms of the mean squared risk, the block estimator is rate minimax for
$\Delta(\bff_n/\sigma)\le \Delta^*_n$ with
a wide range of $\Delta^*_n$ (with no extra logarithmic factor for $d\neq 2$),
achieves near parametric adaptation in the piecewise constant case,
and also achieves near adaptation to the oracle minimax rate in variable selection.
Furthermore, we prove parallel results for the integrated risk
for i.i.d. random designs in $[0,1]^d$ when the joint density of the design point
is uniformly bounded away from zero and infinity.
In addition to Sections 2 and 3, we present in Section 4 some simulation results
to demonstrate the advantage of the block estimator over the LSE in multiple isotonic regression. The full proofs of all theorems, propositions and lemmas in this paper are relegated to the supplement \citep{supplement}.

Here and in the sequel, the following notation is used.
For $\{\ba,\bb\}\subset V$, we say $\bb$ is larger than $\ba$ when $\ba\preceq \bb$,
and we set $[\ba,\bb] =\{\bx\in V: \ba\preceq\bx\preceq\bb\}$ as a block in $G=(V,E)$.
We denote by $n_A$ the number of sampled points in $A$,
i.e. $n_A=\#\{i\le n: \bx_i \in A\}$,
and set $n_{\ba,\bb}=n_{[\ba, \bb]}$, 
$n_{\ba,*} = \# \{i\le n: \ba \preceq \bx_i\}$,
and 
$n_{*, \bb} =\# \{i \le n: \bx_i \preceq \bb\}$.
 For $\ba = (a_1,\ldots,a_d)^T\in \R^d$ and
$\bb = (b_1,\ldots, b_d)^T\in \R^d$, $\ba \preceq \bb$ iff $a_j\le b_j$ for all $1\le j\le d$, and this is also expressed as $\ba\le \bb$.
We denote by $C$ a positive numerical constant, and $C_{\rm index}$ a positive constant depending
on the ``index'' only. For example, $C_{q,d}$ is a positive constant depending on $(q,d)$ only.
For the sake of convenience, the value of such a constant with the same subscript may change from one appearance to the next. We may write $x \lesssim_{\rm index} y $ when $x \le C_{\rm index}\,y$.
Finally, we set $\log_+(x) = 1 \vee \log x$.

\section{The least squares and block estimators}

Given design points $\bx_i\in V$ and responses $y_i\in \R$, the isotonic LSE is formally defined as
\bes
\fhat_n^{(lse)} = \argmin_{f\in\scrF} \sum_{i=1}^n\big\{y_i - f(\bx_i)\big\}^2,
\ees
where $\scrF=\{f: f(\bu)\le f(\bv)\ \forall\ \bu\preceq \bv\}$ is the set of all non-decreasing functions on the directed graph $G = (V,E)$. As the squared loss only involves the value of $f$ at the design points, this LSE is any non-decreasing extension of the LSE of the mean vector $\bff_n = (f(\bx_1),\ldots,f(\bx_n))^T$
in (\ref{iso-reg}),
\bel{lse}
\hbf_n^{(lse)} = \argmin_{\bff_n\in\scrF_n} \|\by-\bff_n\|_2^2,
\eel
where $\by = (y_1,\ldots,y_n)^T$ and $\scrF_n = \{\bff_n: f\in\scrF\}\subset \R^n$. 
As $\scrF_n $ is defined with no more than ${n\choose 2}$ linear constraints, 
$\hbf_n^{(lse)}$ can be computed with quadratic programming. 
Potentially more efficient algorithms for the LSE have been developed in \cite{dykstra1983algorithm}, \cite{kyng2015fast} and \cite{Stout2015}, among others.

As mentioned in the introduction, the LSE $\hbf_n^{(lse)}$ has an explicit representation in
the minimax formula (\ref{minimax-formula}) for isotonic regression on graphs in general
\citep{robertson1988order}, although this fact is better known in the univariate case.
As the high complexity of the general upper and lower sets in the minimax formula seems to be the
cause of the analytical or possibly real gap between the risk of the LSE
and the optimal minimax and adaptation rates, we consider in this paper
block estimators $\fhat_n^{(block)}$ of the form
\bel{block}
&& \min\Big\{\fhat_n^{(max-min)}(\bx), \fhat_n^{(min-max)}(\bx)\Big\}
\cr
&\le& \fhat_n^{(block)}(\bx)
\\ \nonumber
&\le& \max\Big\{\fhat_n^{(max-min)}(\bx), \fhat_n^{(min-max)}(\bx)\Big\},
\quad \forall\ \bx \in V,
\eel
where $\fhat_n^{(max-min)}$ and $\fhat_n^{(min-max)}$ are
the block max-min and min-max estimators given in (\ref{max-min-min-max}).
It is clear from (\ref{max-min-min-max}) that both the max-min and min-max estimators
are non-decreasing on the graph $G=(V,E)$ as the maximum is taken over increasing
classes indexed by $\bx\in V$ and the minimum over decreasing classes.
However, the monotonicity of the block estimator,
$\fhat_n^{(block)} \in \scrF$ or even $\hbf_n^{(block)} \in \scrF_n$, is optional in our analysis.
A practical monotone solution is
\bel{block-mid}
\fhat_n^{(block)}(\bx)  =
\frac{1}{2}\Big\{\fhat_n^{(max-min)}(\bx) + \fhat_n^{(min-max)}(\bx)\Big\},\quad \forall \bx\in V.
\eel

\begin{figure}[h!]
\centering
\includegraphics[width=0.5\textwidth]{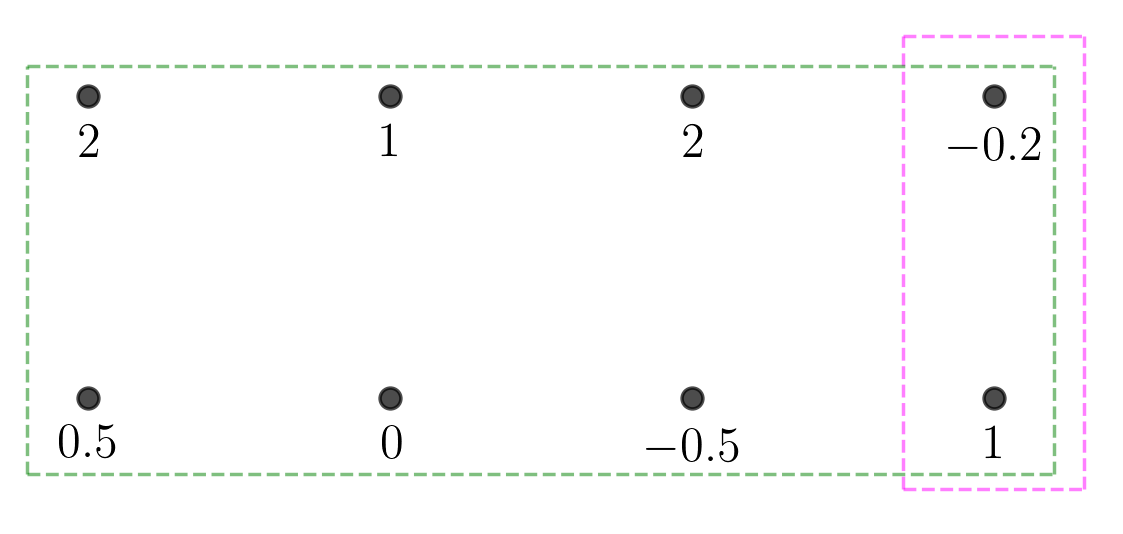}
\caption{\it Responses $y_i$ on a $4\times2$ lattice design:
At design point $\bx=(4,1)$, $\fhat_n^{(max-min)} (\bx) = 0.4$ is attained by the mean
inside the magenta box and $\fhat_n^{(min-max)} (\bx) =  0.725$
attained by the mean inside the green box.
}
\label{fig:example}
\end{figure}

We note that the estimator (\ref{block}) is defined on the entire $V$.
This is needed as we shall consider
the $L_q$ risk (\ref{integrated-risk}) as well as the $\ell_q$ risk (\ref{empirical-risk}).
It would be tempting to define the block estimator by
\bes
\max_{\bu\preceq\bx}\min_{\bx\preceq\bv} \ybar_{[\bu,\bv]}
\le \fhat^{(block)}(\bx)\le \min_{\bx\preceq \bv}\max_{\bu\preceq \bx}\ybar_{[\bu,\bv]}
\ees
\citep{fokianos2017integrated}. However, unfortunately, when $\bx$ is not a design point,
$\ybar_{[\bu,\bv]}$ is undefined when $[\bu,\bv]$ contains no data point,
and $\fhat_n^{(max-min)}(\bx) \le \fhat_n^{(min-max)}(\bx)$ is not guaranteed to hold
even for properly defined max-min and min-max estimators in (\ref{max-min-min-max}),
even in the univariate case.
For example, for $V = [0,1]$ with two data points $(x_1,y_1)=(0,1)$ and $(x_2,y_2)=(1,2)$,
(\ref{max-min-min-max}) gives
$\fhat_n^{(max-min)}(0.5) = 2 > 1 = \fhat_n^{(min-max)}(0.5)$.
We do have
\bel{max-min-inequality}
\fhat_n^{(max-min)}(\bx_i) \le \fhat_n^{(min-max)}(\bx_i),\quad i=1, \ldots, n,
\eel
but the minimax formula $\fhat_n^{(max-min)} = \fhat_n^{(min-max)}$ may fail
even on the design points as the example in Figure \ref{fig:example} demonstrates.

In the rest of this section, we prove that the max-min and min-max
estimators defined with upper and lower sets in a graph $G$, including the LSE,
can always be expressed as the block estimators defined
as in (\ref{block}) but over a larger graph than $G$, so that our analysis of
general block estimators is also relevant to the LSE.
We present our argument in a more general setting as follows.

Formally, a subset of vertices $U\subseteq V$ is called an upper set if
$U=\{\bx: f(\bx) > t\}$ for some $f\in\scrF$ and real $t$, or equivalently the indicator
function $1_U$ is non-decreasing on $G$, i.e. $1_U\in \scrF$;
a subset $L\subseteq V$ is called a lower set if $L=\{\bx: f(\bx) \le t\}$ for some
$f\in\scrF$ and $t\in\R$, i.e. the complement of an upper set.
Let $\scrU$ be the collection of all upper sets, $\scrL$ the collection of all lower sets, and
\bes
\scrU_{\bx} \subseteq \{U\in\scrU: \bx\in U\}\ \hbox{ and }\
\scrL_{\bx} \subseteq \{L\in\scrL: \bx\in L\}
\ees
be certain subsets of the collections of upper and lower sets containing $\bx$.
The max-min and min-max estimator can be defined in general as
\bel{max-min-min-max*}
\fhat_n^{(max-min)}(\bx) &=& \max_{U \in \scrU_{\bx}, n_U >0 }
\min_{L \in \scrL_{\bx}, n_{U \cap L}>0}
\ybar_{U \cap L}, \quad \bx \in V,
\cr \fhat_n^{(min-max)}(\bx) &=&\min_{L \in \scrL_{\bx}, n_{L}>0}
\max_{U \in \scrU_{\bx}, n_{U \cap L} >0 }
\ybar_{U \cap L}, \quad \bx \in V,
\eel
where $n_A=\{i\le n: \bx_i\in A\}$.
These max-min and min-max estimators are non-decreasing in $\bx$ on the entire graph if
$\scrU_{\bx}$ is non-decreasing in $\bx$ and $\scrL_{\bx}$ non-increasing in $\bx$:
$\scrU_{\bx} \subseteq \scrU_{\bx'}$ and $\scrL_{\bx} \supseteq \scrL_{\bx'}$
for all ordered pairs $\bx \preceq \bx'$.

\begin{figure}
\centering
\includegraphics[width=0.7\textwidth]{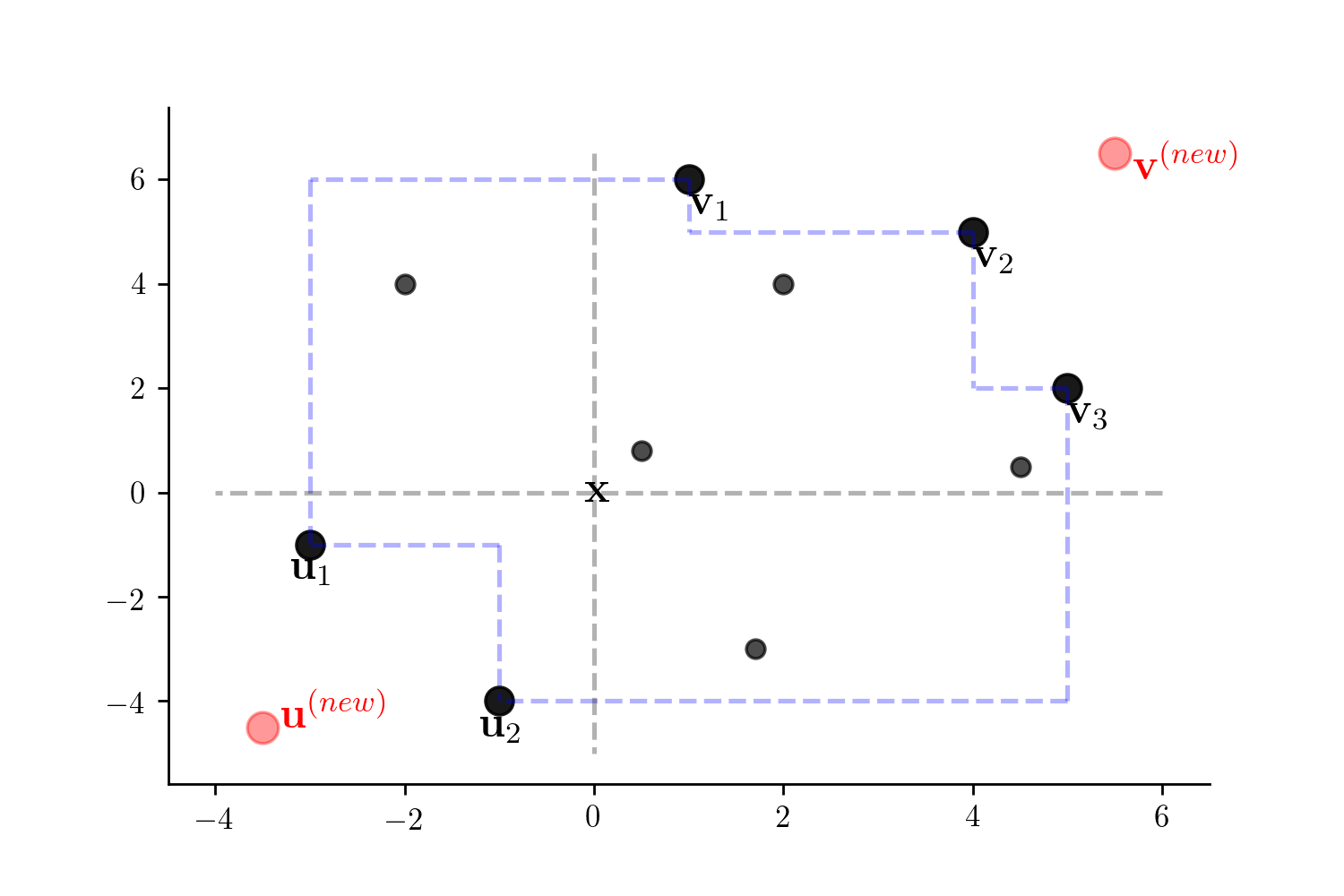}
\caption{\it
Amendment of $G$ to include
$U\cap L= \cup_{j \in \{1,2\}, k \in \{1,2,3\}}[\bu_j,\bv_k]$
where $\bu^{(new)}$ has two inbound edges from $\bu_1$ and $\bu_2$
and $\bv^{(new)}$ has three outbound edges to $\bv_1$, $\bv_2$ and $\bv_3$
}
\label{fig:amend}
\end{figure}

By (\ref{minimax-formula}), the LSE is a special case of (\ref{max-min-min-max*})
when $\scrU_{\bx}$ and $\scrL_{\bx}$ are taken to be the largest possible.
The block max-min and min-max estimators (\ref{max-min-min-max}) are special cases of
(\ref{max-min-min-max*}) with $\scrU_{\bx} = \{[\bu, *]: \bu \preceq \bx \} $
and $\scrL_{\bx} = \{[*, \bv]: \bx \preceq \bv\}$.
Conversely, the LSE, and more generally (\ref{max-min-min-max*}), can be written as
\bel{max-min-min-max-block}
\fhat_n^{(max-min)}(\bx) &=& \max_{\bu\in A_{\bx}, n_{\bu,*} >0 }
\min_{\bv\in B_{\bx}, n_{\bu,\bv}>0}
\ybar_{[\bu,\bv]}, \quad \bx \in V,
\cr \fhat_n^{(min-max)}(\bx) &=&\min_{\bv \in B_{\bx}, n_{*,\bv}>0}
\max_{\bu \in A_{\bx}, n_{\bu,\bv} >0 }
\ybar_{[\bu,\bv]},\quad \bx \in V,
\eel
based on the average response in blocks $[\bu,\bv]$
for suitable $A_{\bx}$ and $B_{\bx}$ in a larger graph $G^*$ in which $G$ is a subgraph.
We define $G^*$ by amending $G$ with new nodes and edges as follows.
For each upper set $U$, we amend $G$ with node
$\bu^{(new)} = \bu^{(new, U)}$
and edges $\{\bu \to \bu^{(new)}: \bu\in U\}$, whereas
for each lower set $L$, we amend $G$ with node
$\bv^{(new)} = \bv^{(new, L)}$
and edges $\{\bv^{(new)}\to \bv: \bv\in L\}$.
Define in the new graph $G^*$ the estimators (\ref{max-min-min-max-block}) with
$A_{\bx} = \{\bu^{(new,U)}: U\in \scrU_{\bx}\}$ and
$B_{\bx} = \{\bv^{(new,L)}: L\in \scrL_{\bx}\}$.
Then, the restriction of (\ref{max-min-min-max-block}) on $G$ is identical to
(\ref{max-min-min-max*}) as $[\bu^{(new,U)},\bv^{(new,L)}]$ contains the same set
of design points as $U\cap L$.
This can be seen as follows. For any pair of upper and lower sets $U$ and $L$,
$[\bu^{(new,U)},\bv^{(new,L)}]\supset U\cap L$ by the definition of $\bu^{(new,U)}$ and
$\bv^{(new,L)}$ and the associated collections of new edges.
On the other hand, for any design point $\bx_i\in [\bu^{(new,U)},\bv^{(new,L)}]$,
$\bu^{(new,U)}\preceq \bx_i$ could happen only if $\bu\preceq \bx_i$
for some $\bu\in U$ as there is no other way to connect to $\bu^{(new,U)}$ in $G^*$,
while $\bx_i \preceq \bv^{(new,L)}$ could happen only if $\bx_i \preceq \bv$ for
some $\bv\in L$. Thus, $\ybar_{U\cap L}=\ybar_{[\bu^{(new,U)},\bv^{(new,L)}]}$.
Figure \ref{fig:amend} demonstrate a $[\bu^{(new)},\bv^{(new)}]$ when
$G$ is a 2-dimensional lattice.

Our theoretical results on general graph in Subsection 3.1 below are applicable to the LSE 
by writing the LSE as a block estimator on a much larger amended graph. 
However, the more specific results in multiple isotonic regression 
in Subsections 3.2-3.7 are not application to the LSE as they are based the calculation of the 
variability bounds in \eqref{r-plus} and \eqref{r-minus} below for the lattice and random designs, 
not on the enlarged graph. 

\section{ Theoretical results }

In this section, we first analyze the block estimator $\fhat_n^{(block)}(\bx)$ in (\ref {block}) for graphs under the most general setting.
Specific risk bounds are then given for multiple isotonic regression with fixed lattice designs and random designs.

\subsection{General isotonic regression on graph}
We shall extend the risk bounds of \cite{zhang2002} from the real line to general graphs.
To this end, we first derive an upper bound for the total risk in subsets $V_0\subset V$,
\bes
T_q(V_0) = \sum_{\bx_i\in V_0} \E \Big| \fhat_n^{(block)}(\bx_i) - f(\bx_i)\Big|^q,
\ees
based on the value of the true $f$ on $V_0$. Such bounds automatically produce adaptive
risk bounds when the true $f$ is ``piecewise constant'' in a partition of $V$.
Given $V_0$, let $r_{q,+}(m)$ be a non-increasing function of $m\in \N^+$ satisfying
\bel{r-plus}
r_{q,+}(m) \ge \max\bigg\{
\E\bigg(\max_{\bu \preceq \bx} \sum_{\bx_i\in [\bu,\bv]} \frac{\veps_i}{n_{\bu,\bv}}\bigg)_+^q:
n_{\bx,\bv}=m, \bx \preceq \bv \hbox{ and } \bv \in V_0\bigg\}.
\eel
This function bounds the error of the block estimator from the positive side
when the positive part of its bias is no greater than the positive part of the maximum average of at least $m$ noise variables.
Similarly, to control the estimation error from the negative side,
let $r_{q,-}(m)$ be a non-increasing function satisfying
\bel{r-minus}
r_{q,-}(m) \ge  \max\bigg\{
\E\bigg(\min_{\bv\succeq \bx}\sum_{\bx_i\in [\bu,\bv]} \frac{\veps_i}{n_{\bu,\bv}}\bigg)_-^q:
n_{\bu,\bx}=m, \bu \preceq \bx \hbox{ and } \bu \in V_0 \bigg\}.
\eel
With the above functions $r_{q,\pm}(m)$, we define for $\bx\in V_0$
\bel{k-m}
m_{\bx,-}&=&  \max \Big\{n_{\bu,\bx}: f(\bu) \ge f(\bx) - r_{q,-}^{1/q}(n_{\bu,\bx}),
\bu \preceq \bx \hbox{ and } \bu \in V_0  \Big\},
\cr \bu_{\bx} &=& \argmax_{\bu\in V_0 :\, \bu \preceq \bx} \Big\{n_{\bu,\bx}:
f(\bu) \ge f(\bx) - r_{q,-}^{1/q}(n_{\bu,\bx})\Big\},
\\ \nonumber m_{\bx} = m_{\bx,+} &=& \max \Big\{n_{\bx,\bv}: f(\bv) \le f(\bx)+r_{q,+}^{1/q}(n_{\bx,\bv}),
\bx \preceq \bv \hbox{ and } \bv \in V_0 \Big\},
\\ \nonumber \bv_{\bx} &=& \argmax_{\bv\in V_0:\, \bx \preceq \bv} \Big\{n_{\bx,\bv}:
f(\bv) \le f(\bx)+r_{q,+}^{1/q}(n_{\bx,\bv})\Big\}.
\eel
Roughly speaking, the above quantities provide configurations in which the bias of $\fhat_n(\bx_i)$ is of no greater order than its variability from the negative and positive sides,
so that the error of the block estimator is of no greater order than
an average of $m_{\bx_i,-}$ noise variables on the negative side and
the average of $m_{\bx} = m_{\bx_i,+}$ noise variables on the positive side.
Thus, it makes sense to count the frequencies of $m_{\bx_i,-}$ and $m_{\bx_i} $
as follows,
\bel{ell-pm}
\ell_-(m) = \#\big\{i: \bx_i \in V_0, m_{\bx_i,-}\le m\big\},\ \
\ell_+(m) = \#\big\{i: \bx_i \in V_0, m_{\bx_i}\le m\big\}.
\eel
We note that the functions $r_{q,\pm}$ in (\ref{r-plus}) and (\ref{r-minus}) do not depend on $f$, and
all the quantities in (\ref{k-m}) and (\ref{ell-pm})
depend on the true $f$ only through $\{f(\bx): \bx\in V_0\}$.

\begin{theorem}\label{th-1} 
Assume $f$ is non-decreasing on a graph $G=(V,E)$.  
Let $r_{q,\pm}(m)$ be given by (\ref{r-plus}) and (\ref{r-minus}),
and $\ell_\pm(m)$ by (\ref{ell-pm}). Then it holds for any block estimator $\fhat_n^{(block)}(\bx)$ in (\ref {block}) 
that
\bel{th-1-1}
&& \E \Big\{\fhat_n^{(block)}(\bx_i) - f(\bx_i)\Big\}_+^q\le 2^qr_{q,+}(m_{\bx_i}),\ \forall \bx_i\in V_0, 
\\
\nonumber &&
\E \Big\{\fhat_n^{(block)}(\bx_i) - f(\bx_i)\Big\}_-^q\le 2^qr_{q,-}(m_{\bx_i,-}),\ \forall \bx_i\in V_0.
\eel
Consequently, for any upper bounds $\ell^*_{\pm}(m) \ge \ell_\pm(m)$ with $\ell^*_\pm(0)=0$,
\bel{th-1-2}
T_q(V_0)
&\le& \sum_{m=1}^{\infty}2^q r_{q,+}(m)\Big\{\ell^*_+(m)-\ell^*_+(m-1)\Big\}
\\
\nonumber && \qquad \qquad
 + \sum_{m=1}^{\infty}2^q r_{q,-}(m)\Big\{\ell^*_-(m)-\ell^*_-(m-1)\Big\}.
\eel
\end{theorem}

Theorem \ref{th-1} provides risk bound for the block estimator (\ref{block}) over a
subset $V_0$ of design points in terms of upper bound functions $r_{q,\pm}(m)$ and $\ell^*_\pm(m)$.
Ideally, we would like to have
\bel{r-bd}
r_{q,\pm}(m) = C_{q,d}\sigma^q m^{-q/2}
\eel
in (\ref{r-plus}) and (\ref{r-minus}).
When the design points in $V_0$ are linear and the $(q \vee 2)$-th moment of the noise variable
is uniformly bounded, (\ref{r-plus}) and (\ref{r-minus}) hold for the above choice of $r_{q,\pm}(m)$.
This choice of $r_{q,\pm}(m)$ is also valid when $V$ is a lattice in $\R^d$
and $\veps_i$ are independent variables with uniformly bounded $(q \vee 2)$-th moment, as we will prove in Subsection 3.3. 

\subsection{Minimax lower bound in multiple isotonic regression with lattice designs}

We study in the rest of this section multiple isotonic regression in $V \subseteq \R^d$ where $\ba \preceq \bb$ iff $\ba\le\bb$,
i.e. $a_j\le b_j\ \forall\ 1\le j\le d$, for all $\ba=(a_1,\ldots,a_d)^T$ and $\bb= (b_1,\ldots, b_d)^T$,
and $\scrF$ is the class of all non-decreasing functions $f(t_1,\ldots,t_d)\uparrow t_j,\, \forall\ j=1,\ldots,d$.

\medskip
The lattice design we are considering 
is given by
\bel{lattice}
V = \big\{\bx_i: 1\le i\le n\big\} = [{\bf 1},\bn]=\prod_{j=1}^d \{1,\ldots,n_j\},  
\eel
where $\bn=(n_1,\ldots,n_d)^T$ with positive integers $n_j$ and $n = \prod_{j=1}^d n_j$. 
Here $[{\bf 1},\bn]$ is treated as a set of integer-valued vectors in $\N^d$, forming a lattice. 
Occasionally, we may also use $[\bu,\bv]$ to denote a hyper-rectangle of real numbers in continuum. 
This slight abuse of notation typically would not lead to confusion, for example in $\bx_i\in [\bu,\bv]$, but we would be specific if necessary. 
Without loss of generality, we assume in this subsection $n_1 \ge n_2 \ge \cdots \ge n_d$.
In the above lattice design, we provide a minimax lower bound in multiple isotonic regression as follows.

\begin{proposition}\label{prop-lower-bd}
Suppose $\veps_i\sim N(0,\sigma^2)$.
Let $\Delta(\bff_n/\sigma)=\{f(\bn)-f(\bf 1)\}/\sigma$,
$n_{d+1}=1$, $n_s^* = \prod_{j=1}^s n_j$, $t_s = n^*_s/n_s^s$,
$t_{d+2}=\infty$ and $s_q = \lceil 2/(q-1)\rceil \wedge (d+1)$.
Let $h_0(t)=\Delta^*_n\sqrt{t}$ and define piecewise
$H(t) = \min\big\{1, h_0(t)/(n^*_s/t)^{1/(s\wedge d)}\big\}, t\in [t_s,t_{s+1}], s=1,\ldots,d+1$. Then,

\bel{prop-lower-bd-1}
&& \inf_{\hbf}\, \sup\Big\{ R_q(\hbf, \bff_n):
\bff_n \in\scrF_n, \Delta(\bff_n/\sigma)\le\Delta_n^*\Big\}
\\ \nonumber &\gtrsim_{q,d}& \sigma^q\, \max
\Big\{(t\wedge n)^{-q/2}H(t):\ t \wedge h_0(t) \ge 1\Big\}
\\ \nonumber &=&
\sigma^q \times
\left\{ \begin{aligned}
& 1, && n_1 \le \Delta^*_n, && (s=0)
\cr
& \displaystyle \big(\Delta^*_n/(n_s^*)^{1/s}\big)^{qs/(2+s)}
&& \displaystyle n_{s+1}/t_{s+1}^{1/2}
\le \Delta^*_n \le n_{s}/t_{s}^{1/2},
&& (1\le s < s_q)
\cr
& \displaystyle \Delta_n^*/\big(n_s t_{s}^{(q-1)/2}\big),
&& \displaystyle t_{s}^{-1/2}
\le \Delta^*_n \le n_{s}/t_{s}^{1/2},
&& (s=s_q\le d)
\cr
& \displaystyle (\Delta_n^*)^{q-2/s}/(n^*_s)^{1/s},
&& \displaystyle t_{s+1}^{-1/2}
\le \Delta^*_n \le t_{s}^{-1/2},
&& (s_q\le s\le d)
\cr
& n^{-q/2}, &&  0\le \Delta^*_n\le n^{-1/2}.
&& (s = d+1)
\end{aligned}
\right.
\eel
In particular, when $n_1=\cdots=n_d=n^{1/d}$ and $\Delta_n^* \ge n^{-1/2}$,
the right-hand side of (\ref{prop-lower-bd-1}) is
\bel{prop-lower-bd-2}
\sigma^q \times \begin{cases}
\min\big\{1, \big(\Delta_n^*/n^{1/d}\big)^{qd/(d+2)} \big\}, &
q \le 1+2/d,
\cr \min\big\{1,\Delta_n^*/n^{1/d}, (\Delta_n^*)^{q-2/d}\big/n^{1/d}\big\},
& q \ge 1 + 2/d.
\end{cases}
\eel
\end{proposition}

On the right-hand side of (\ref {prop-lower-bd-1}), the breaking points on $[0, \infty)$ for $\Delta_n^*$ are
\bes
0, n^{-1/2}=t_{d+1}^{-1/2}, t_d^{-1/2}, \ldots, t_{s_q}^{-1/2}, n_{s_q}/t_{s_q}^{1/2}, \ldots, n_1/t_1^{1/2} = n_1.
\ees
Note that $1$ lies in between $t_{s_q}^{-1/2} $ and $ n_{s_q}/t_{s_q}^{1/2}$.
The above minimax lower bound also depends on the loss function through $q$ and the dimension of the lattice.
For $q\ge 3$, we have $s_q=1$, so that
\bes
\inf_{\hbf}\, \sup\Big\{ R_q(\hbf, \bff_n):
\bff_n \in \scrF_n, \Delta(\bff_n/\sigma)\le\Delta_n^*\Big\}
\gtrsim_{q,d} \sigma^q \min \Big(1, \Delta_n^*/n_1\Big)
\ees
for $\Delta^*_n\ge 1$.
However for $q=2$, we have $s_q=2$, so that (\ref{prop-lower-bd-1}) yields
\bel{minimax-q2}
&& \inf_{\hbf}\, \sup\Big\{ R_2(\hbf, \bff_n):
\bff_n \in \scrF_n, \Delta(\bff_n/\sigma)\le\Delta_n^*\Big\}
\\ \nonumber &\gtrsim_{d}& \sigma^2 \times
\left\{ \begin{aligned}
& 1, && n_1 \le \Delta_n^* , && (s=0)
\cr
& \displaystyle \big(\Delta_n^*/n_1\big)^{2/3}, && \displaystyle n_2^{3/2} / n_1^{1/2} \le \Delta_n^* \le n_1, && (s=1)
\cr
& \Delta_n^*/(n_1n_2)^{1/2},
&& \sqrt{n_2/n_1}\le \Delta_n^* \le n_2^{3/2} / n_1^{1/2}. && (s=2)
\end{aligned}
\right.
\eel
For $\Delta_n^*\asymp 1$,
this matches the lower bound for the $\ell_2$ minimax rate in
\cite{chatterjee2018} for $d=2$ and \cite{han2017isotonic} for $d \ge 3$.
For $5/3\le q<2 \le d$, we have $s_q=3$.

If (\ref{prop-lower-bd-1}) is achievable, the integer parameter $s$
can be viewed as the effective dimension of the isotonic regression problem
as the rate depends on $\bn$ only through $n_1,\ldots, n_s$ when
$n_{s+1}$ is sufficiently small; the rate would also be
achievable by separate $s$-dimensional isotonic regression in the
$\prod_{j=s+1}^d n_j = n/n^*_s$ individual $s$-dimensional sheets
with fixed $x_{s+1},\ldots,x_d$.
For example, in (\ref{minimax-q2}), the minimax rate can be achieved by
$\hbf_n=\by$ for $s=0$, by the row-by-row univariate isotonic regression for $s=1$,
and by individual bivariate isotonic least squares
up to a factor of $(\log n)^4$ for $s=2$ \citep{chatterjee2018}.
We will prove in the next subsection that the block estimator (\ref{block}) achieves
the rate in (\ref{prop-lower-bd-1}) for a wide range of $\Delta^*_n$, so that
Proposition \ref{prop-lower-bd} indeed provides the minimax rate.

In the proof of Proposition \ref{prop-lower-bd}, we divide $[{\bf 1},\bn']\subset V =[{\bf 1}, \bn]$ into a
$K_1\times\cdots\times K_d$ lattice of hyper-rectangles of size 
$m_1\times\cdots\times m_d$, indexed by $\bk=(k_1,\ldots,k_d)^T$, 
$k_j = 1,\ldots,K_j$, $j=1,\ldots,d$, 
and consider the class of piecewise constant functions $f(\bx) = g(\bk)$
satisfying
\bes
g(\bk) = \sigma \min\Big\{\Delta_n^*, (m^*)^{-1/2}
\big[\theta(\bk) + (k_1+\cdots+k_d - k^*)_+\big] \Big\},\ \theta(\bk) \in \{0,1\},
\ees
and $f(\bx)=\sigma\Delta^*_n$ for $\bx \in [{\bf 1},\bn]\setminus [{\bf 1},\bn']$, where $m^*=\prod_{j=1}^dm_j$ is the size of the hyper-rectangle. 
As $g(\bk)$ is non-decreasing in $k_j$ for each $j$ for all $\theta(\bk)\in \{0,1\}$,
this construction provides a lower bound for the $\ell_q$ risk proportional to the product of
$\sigma^q(m^*)^{-q/2}$ and the number of free $\theta(\bk)$.
This is summarized in the following lemma.

\begin{lemma}\label{lm-lower-bd}
Under the conditions of Proposition \ref{prop-lower-bd},
\bel{lm-lower-bd-1}
&& \inf_{\hbf}\, \sup\Big\{\E \big\|\hbf - \bff_n \big\|_q^q:
\bff_n\in\scrF_n, \Delta(\bff_n/\sigma)\le\Delta_n^*\Big\}
\\
\nonumber &\ge& c_qc_d\sigma^q n
\max_{\bm \in \scrM} \bigg\{ \frac{1}{(m^*)^{q/2}}
\min\bigg( \frac{\sqrt{m^*}\Delta_n^*}{\max_j \lfloor n_j/m_j\rfloor}, 1\bigg) \bigg\},
\eel
where $c_q = \inf_{\delta}\E_{\mu\sim \hbox{\rm\footnotesize Bernoulli}(1/2)} \big| \delta(N(\mu,1)) - \mu\big|_q^q$
is the Bayes risk for estimating $\mu$ with the Bernoulli$(1/2)$ prior based on a single $N(\mu,1)$ observation,
$c_d$ is a constant depending on $d$ only,
\bes
\scrM = \Big\{\bm = (m_1,\cdots, m_d): m_j \in \N_+,  m_j\le n_j\, \forall j \le d, \sqrt{m^*}\Delta_n^* \ge 1\Big\},
\ees
and $m^* = \prod_{j \le d}m_j$.
Moreover, the optimal configuration of $\bm$ in (\ref{lm-lower-bd-1}) must satisfy either
$m_j=1$ or $\lfloor n_j/m_j\rfloor  = \max_{1\le j\le d}\lfloor n_j/m_j\rfloor$ for each $j$.
\end{lemma}

\subsection{The block estimator in multiple isotonic regression with lattice designs}

We further divide this subsection into three separate sub-subsections to study the performance of 
the block estimator at a single design point $\bx_i$, 
in an arbitrary sub-block $[\ba,\bb]\subset [{\bf 1},\bn]$, 
and on the entire lattice $[{\bf 1},\bn]$.  
It is of great interest to show that the block estimator in (\ref {block}) 
matches the minimax lower bound given in Proposition \ref{prop-lower-bd},
which will be done in the third sub-subsection for general $q$ and $d$.

\subsubsection{Risk of the block estimator at a single design point} 
For any given point in the design lattice, 
the following proposition asserts that the block estimator matches 
certain one-sided oracle estimators in the rate of one-sided $L_q$ risks.

\begin{proposition}\label{prop-single-point} 
Let $\fhat_n^{(block)}(\bx)$ be the block estimator in (\ref {block})
with the lattice design $V=[{\bf 1},\bn]$ in (\ref{lattice}).  
Let $q\ge 1$ and $r_{q,\pm}(m)$ be as in (\ref{r-plus}) and (\ref{r-minus}). 
Assume $\veps_i$ are independent $N(0,\sigma^2)$ random variables. 
Then, for any design point $\bx_i\in [{\bf 1},\bn]$, 
\bel{new-prop-2-1}
\E \Big(\fhat_n^{(block)}(\bx_i) - f(\bx_i)\Big)_+^q 
\le 2^qr_{q,+}(m_{\bx_i}) 
\le C_{q,d} \min_{\bx_i \le \bv\le \bn}\E\Big(\ybar_{[\bx_i,\bv]}-f(\bx_i)\Big)_+^q\,,
\eel
where $\ybar_{[\bu,\bv]}=\sum_{\bu\le \bx_i\le \bv}y_i/n_{\bu,\bv}$, and 
\bel{new-prop-2-2}
\E \Big(\fhat_n^{(block)}(\bx_i) - f(\bx_i)\Big)_-^q 
\le 2^qr_{q,-}(m_{\bx_i})
\le C_{q,d} \min_{{\bf 1}\le \bu \le \bx_i}\E\Big(\ybar_{[\bu,\bx_i]}-f(\bx_i)\Big)_-^q\,.
\eel
Consequently, with $\E_g$ being the expectation under which $y_i = g(\bx_i)+\veps_i$, 
\bel{new-prop-2-3}
&& \E \Big|\fhat_n^{(block)}(\bx_i) - f(\bx_i)\Big|^q 
\\ \nonumber &\le& C_{q,d} \min_{\bu \le \bx_i\le\bv}
\bigg\{\E_g\Big|\ybar_{[\bu,\bv]}-g(\bx_i)\Big|^q: g\in\scrF, g(\bv)=f(\bv)\ \forall \bv\ge \bx_i \bigg\}
\\ \nonumber && + C_{q,d} \min_{\bu \le \bx_i\le\bv}
\bigg\{\E_g\Big|\ybar_{[\bu,\bv]}-g(\bx_i)\Big|^q: g\in\scrF, g(\bu)=f(\bu)\ \forall \bu\le \bx_i \bigg\}. 
\eel
\end{proposition} 

Suppose we are confined to consider only block mean estimators $\ybar_{[\bu,\bv]}$ 
with no negative bias in the estimation of $f(\bx_i)$ 
but we also want to control the positive side of the error. 
As $f$ is non-decreasing but otherwise unknown, 
we are thus forced to choose $\bu\ge \bx_i$. As $\ybar_{[\bu,\bv]}$ with $\bx_i\le\bu\le\bv$ 
would have larger bias and variance than $\ybar_{[\bx_i,\bv]}$, the optimal $[\bu,\bv]$ is given by 
\bes
\min_{ \bu = \bx_i  \le \bv\le \bn}\E\Big(\ybar_{[\bx_i,\bv]}-f(\bx_i)\Big)_+^q. 
\ees
The above minimum can be viewed as an oracle benchmark under the no-negative-bias constraint 
as the solution of the optimal $\bv$ still depends on $f$. 
Although the block estimator (\ref {block}) is unlikely to be unbiased, 
\eqref{new-prop-2-1} and \eqref{new-prop-2-2} 
assert that its one-sided risks match the rates of 
such oracle benchmarks from both the positive and negative sides. 
Another interpretation of the performance of the block estimator is 
\eqref{new-prop-2-3} in which the oracle expert has to guard against 
the worst case scenarios in the uncertainty of $f$ on either sides, but not simultaneously on both. 

We prove Proposition \ref{prop-single-point} with an application of Theorem \ref{th-1}. 
This requires more explicit variability bounds $r_{q,\pm}(m)$ in (\ref{r-plus}) and (\ref{r-minus}) 
as in (\ref {r-bd}). This validity of (\ref {r-bd}) is a consequence of the following lemma, 
which extends Doob's inequality to certain multiple indexed sub-martingales. It plays a key role in removing the normality assumption on the noise
$\veps_1,\ldots, \veps_n$ in our analysis. 

\begin{lemma}\label{lm-doob} Let $\scrT=\scrT_1\times\cdots \times \scrT_d \subseteq \R^d$
be an index set with $\scrT_j\subseteq \R$.
Let $\{f_{\bt}, \bt\in\scrT\}$ be a collection of random variables.
Suppose for each $j$ and each $(s_1,\ldots,s_{j-1}, t_{j+1},\ldots,t_d)$,
$\{f_{s_1,\ldots,s_{j-1}, t, t_{j+1},\ldots,t_d}, t\in \scrT_j\}$ is a sub-martingale
with respect to certain filtration $\{\scrF_t^{(j)}, t\in \scrT_j\}$.
Then, for all $q>1$ and $\bt\in \scrT$,
\bes
\E \max_{\bs\in\scrT, \bs\le\bt}\big|f_{\bs}\big|^q \le (q/(q-1))^{qd}\,\E\big|f_{\bt}\big|^q.
\ees
In particular when $\veps_i$'s are independent random variables with $\E \veps_i=0$,
\bes
\E \max_{\bs\le\bt}\bigg|\sum_{\bx_i\le \bs} \veps_i\bigg|^q
\le \begin{cases} (q/(q-1))^{qd}\,\E\big|\sum_{\bx_i\le \bt} \veps_i\big|^q,& q\ge 2,
\cr
\Big(4^d\E\big|\sum_{\bx_i\le \bt} \veps_i\big|^2\Big)^{q/2}, & 1\le  q < 2.
\end{cases}
\ees
\end{lemma}

\subsubsection{Risk of the block estimator in a sub-block}
To automatically deal with adaptation which gives better risk bound when $f(\cdot)$ is piecewise constant, we first consider the risk in one of such ``piece'', 
a hyper-rectangle $[\ba, \bb] \subseteq V = [{\bf 1}, \bn]$.

\begin{theorem}\label{th-lattice-general}
Let $\fhat_n^{(block)}(\bx)$ be the block estimator in (\ref {block})
with the lattice design $V=[{\bf 1},\bn]$ in (\ref{lattice}).
Assume $\veps_i$ are independent random variables with $\E\, \veps_i =0$ and
$\E |\veps_i|^{q\vee 2} \le \sigma^{q\vee2}$.
Let $\ba\le \bb$ be integer vectors in $V=[{\bf 1},\bn]$ and $\ntil_j=b_j-a_j+1$.
Suppose $\ntil_1 \ge  \cdots \ge \ntil_d$.
Define $\ntil = n_{\ba,\bb}$, $\ntil_{d+1}=1$, $\ntil_s^*=\prod_{j=1}^s \ntil_j$ and
$t_s = \ntil_s^*/\ntil_s^s$ (with $1=t_1\le \cdots \le t_d\le t_{d+1}=\ntil$). 
Then, for $q\ge 1$ and any $f\in\scrF$ with
$\Delta_{\ba, \bb}(\bff_n/\sigma) = \{f(\bb)-f(\ba)\}/\sigma \le\Delta_n^*$,
\bel{th2-main}\qquad
T_{q}([\ba,\bb])
&=&\sum_{\bx_i\in [\ba,\bb]}\E \big|\fhat_n^{(block)}(\bx_i) - f(\bx_i)\big|^q
\\ \nonumber
&\le& C^*_{q,d}\,
n_{\ba, \bb}\,\sigma^q \bigg( \Htil(1) + \int_1^{n_{\ba,\bb}}\frac{\Htil(dt)}{t^{q/2}} + \frac{1}{n_{\ba,\bb}} 
\prod_{j=1}^d \int_0^{\ntil_j}\frac{dt}{(t\vee 1)^{q/2}}\bigg),
\eel
where $\Htil(t)$ is a non-decreasing and continuous function of $t$, defined piecewise by
$\Htil(t) = \min\big\{1, \Delta^*_n t^{1/2}(t/\ntil^*_s)^{1/s}\big\}$ for
$t_s \le t\le t_{s+1}$, $s=1,\ldots, d$,
and $C^*_{q,d}$ is continuous in $q\in [1,\infty)$ and non-decreasing in $d$.
Moreover,
  \bel{th2-1}
&& \Htil(1)+ \int_1^{n_{\ba,\bb}}t^{-q/2}\Htil(dt)
\\
\nonumber &\lesssim_{q,d}& 
\left\{
  \begin{aligned}
& 1, && \ntil_1 \le \Delta_n^*, && (s = 0)
\cr \displaystyle
& \big(\Delta^*_n/(\ntil_s^*)^{1/s}\big)^{qs/(2+s)},
&& \displaystyle \ntil_{s+1}/t_{s+1}^{1/2}
\le \Delta^*_n \le \ntil_{s}/t_{s}^{1/2}, &&(1\le s <s_q)
\cr \displaystyle
& \big(\Delta_n^*/\big(\ntil_s t_{s}^{(q-1)/2}\big)\big)\Lambda_s,
&& \displaystyle
\Delta_n^* \le \ntil_{s}/t_{s}^{1/2}, && (s=s_q\le d)
\end{aligned}
\right.
\eel
where $s_q = \lceil 2/(q-1) \rceil \wedge (d+1)$ is as in Proposition \ref{prop-lower-bd} and
\bel{lambda}
\Lambda_s = \bigg[\log_+ \Big( \min \Big\{\frac{\ntil_s}{\ntil_{s+1}},
\frac{\ntil_s/(\ntil_{s}^*)^{1/(s+2)}}{(\Delta_n^*)^{2/(s+2)}} \Big\} \Big)\bigg]^{I\{2/(q-1)=s\}}.
\eel
\end{theorem}

\begin{remark}
The last component on the right-hand side of (\ref {th2-main}) is bounded by
\bel{th2-2}\qquad 
\sigma^q \prod_{j=1}^d \int_0^{b_j-a_j+1}\frac{dt}{(t\vee 1)^{q/2}} 
\,\lesssim_{q,d}\, \sigma^q \bigg[n_{\ba, \bb}^{1-q/2} + \Big(\prod_{j=1}^d \log_+(b_j -a_j+1)\Big)^{I\{q=2\}}\bigg].
\eel
When $\Delta_{\ba,\bb}(\bff_n/\sigma)=0$, $\Htil(t)=0$ for all $t$, 
so that (\ref {th2-2}) is an upper bound 
for the rate of the total risk $T_q([\ba, \bb])$ in the block $[\ba,\bb]$ by Theorem \ref{th-lattice-general}, 
for any $\ba \le \bb$. This yields the adaptation rate stated in Subsection 3.4. 
\end{remark}

\begin{remark}
The function $\Htil(t)$ is defined in the same way as $H(t)$ is in Proposition~\ref{prop-lower-bd} but for the dimensions $\{\ntil_j = b_j-a_j+1, j\le d\}$ of $[\ba, \bb]$ and range-to-noise ratio within $[\ba, \bb]$. 
When $[\ba, \bb] = [{\bf 1}, \bn]$, we have $\Htil(t) = H(t)$ for all $t \in [1, n]$. 
Thus, as discussed below (\ref{minimax-q2}), the integer parameter $s$ in \eqref{th2-1}, completely determined by $\{\ntil_j\}$, $\Delta^*_n$ and $q$, has the interpretation as the effective dimension for the estimation of $f$ in $[\ba,\bb]$ subject to $\{f(\bb)-f(\ba)\}/\sigma \le\Delta_n^*$. 
We note that as $\Htil(t)$ is a smooth fit of pieces proportional to $t^{1/2+1/s}$ or 1, the upper limit of the integration is actually 
$t_*=\min\{t\ge 1: \Htil(t)=1$ or $t=\ntil_{\ba,\bb}\}$, which depends on $\Delta_n^*$, 
and the effective dimension $s$ is then determined by the comparison between $t_*$ 
and $t_s$ and the critical $s_q$. 
\end{remark}

\begin{figure}
  \centering
  \includegraphics[width=0.8\textwidth]{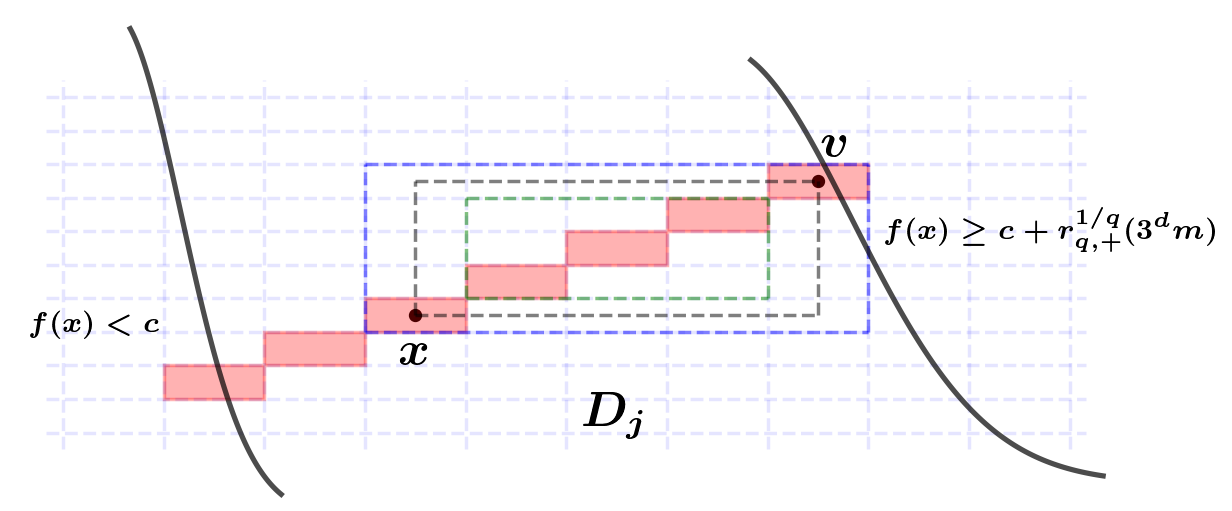}
  \caption{\it 
  Upper bound for the number of design points with $m_{\bx_i}\le m$, an example: 
  $d=2$, $t_d=\ntil/\ntil_d^d$, $m=k^d t_d$ with $k=3$, 
  a line segment of unit blocks in the anti-diagonal direction is colored in red, depicting its intersection with the region $D_j$ between two contours of $f$;
  $\bx$ is a design point $k$ blocks away from the upper boundary of $D_j$, $\bv\in D_j$; 
  $m$, $n_{\bx,\bv}$ and the upper bound $(k+2)^dm$ are respectively the number of points inside the rectangles colored in dashed green, gray and blue;
  as $m_{\bx}\ge n_{\bx,\bv}>m$ in this example, design points inside the intersection of $D_j$ and these red unit blocks with $m_{\bx_i}\le m$ must belong to one of the $k+1=4$ upper-right unit blocks colored in red, and there are at most $(k+1)t_d=4t_d$ such points in this example with $k=3$. 
  For general $k$ and $m=k^dt_d$, $(k+1)t_d\le 2m^{1/d}t_d^{1-1/d}$. }
  \label{fig:proof}
  \end{figure}

In addition to the validity of (\ref {r-bd}) as variability bounds in (\ref{r-plus}) and (\ref{r-minus}), which follows from Lemme \ref{lm-doob}, the proof of Theorem \ref{th-lattice-general} requires   
the complexity bounds for the $\ell_{\pm}(m)$ in (\ref{ell-pm}). We outline here an analysis of the count $\ell_+(m)$ in (\ref{ell-pm}) in the 
case where $\ntil_j/\ntil_d$ are integers and $m\ge t_d = \ntil/\ntil_d^d$. 
We note that $t_d=1$ when $\ntil_j = \ntil_d$ for all $j$. 
Upper bounds for both $\ell_{\pm}(m)$ in the general setting 
are given in the proof of Theorem \ref{th-lattice-general} in subsection A3.3 of the supplement. 

To find upper bounds for $\ell_{+}(m)$, we partition $V_0 = [\ba,\bb]$ into an $\ntil_d\times\cdots \times \ntil_d$ lattice of small ``unit blocks" of size $(\ntil_1/\ntil_d)\times\cdots\times (\ntil_d/\ntil_d)$, each composed of $t_d = \ntil/\ntil_d^d$ design points. 
Consider a line of such unit blocks $L_{\bk}$ in the ``anti-diagonal'' direction and a region $D_j$ between two contours of the unknown $f(\bx)$ at the levels $c$ and $c+r_{q,+}^{1/q}(3^dm)$. 
In Figure \ref{fig:proof}, we color in red the unit blocks in $L_{\bk}$ with nonempty intersection with $D_j$. 
Due to the monotonicity of the $\ell_+(m)$, it suffices to consider $m = k^dt_d$ for some integer $k\ge 1$. 
If $\bx\le \bv$ in $L_{\bk}\cap D_j$ are separated by $k$ unit blocks as depicted in Figure \ref{fig:proof}, then $m=k^dt_d < n_{\bx,\bv} \le (k+2)^dt_d \le 3^dm$ and $f(\bv) - f(\bx) \le r_{q,+}^{1/q}(3^dm)\le r_{q,+}^{1/q}(n_{\bx,\bv})$, so that $m_{\bx}\ge n_{\bx,\bv}>m$. 
Thus, the intersection contains no more than $(k+1)t_d\le 2m^{1/d} t_d^{1-1/d}$ design points $\bx_i$ with $m_{\bx_i}\le m$, all within $k$ unit blocks from the upper contour. 
Let $J = \lceil\{f(\bb)-f(\ba)\}/r_{q,+}^{1/q}(3^dm)\rceil$.  
We divide $[\ba,\bb]$ into $J$ such regions $D_j$ between consecutive contours with $\ba\in D_1$ and $\bb\in D_J$. The last region $D_J$ is special. 
For $\bx\in D_J$ with $n_{\bx,\bb}>m$, there must exist $\bv\in [\bx,\bb]$ such that $m < n_{\bx,\bv}\le 2m$, so that $m_{\bx} \ge n_{\bx,\bv}>m$ due to $f(\bv)\le f(\bx)+r_{q,+}^{1/q}(3^dm)\le f(\bx)+r_{q,+}^{1/q}(n_{\bx,\bv})$. 
Thus, as there are no more than $d \ntil_d^{d-1}$ such $L_{\bk}$ and $J-1\le \{f(\bb)-f(\ba)\}/r_{q,+}^{1/q}(3^dm)\le \Delta^*_n\sigma/r_{q,+}^{1/q}(3^dm)$ 
regions $D_j$ not containing $\bb$, for $m=k^dt_d$ with integer $k\ge 1$ 
\bes
&& \ell_+(m) 
\cr &\le& \min\Big\{\ntil, d\ntil_d^{d-1}\Big(\Delta^*_n\sigma/r_{q,+}^{1/q}(3^dm)\Big)
\Big(2m^{1/d} t_d^{1-1/d}\Big)\Big\}
+ \#\big\{\bx_i\in [\ba,\bb]: n_{\bx_i,\bb}\le m\big\}
\cr &=& \ntil \min\Big\{1, m^{1/d+1/2}\big(\Delta^*_n/\ntil^{1/d}\big)\big(2d 3^{d/2}\big/C_{q,d}^{1/q}\big)\Big\}
+ \#\big\{\bx_i\in [\ba,\bb]: n_{\bx_i,\bb}\le m\big\}
\ees
with the variability bound $r_{q,+}(m)=C_{q,d}\sigma^q m^{-q/2}$ in \eqref{r-bd}. 
It follows that 
\bel{ell-star}
\ell_\pm(m) \le \ell^*_\pm(m) 
= \ntil \Htil(m) + \#\big\{\bx_i\in [\ba,\bb]: n_{\bx_i,\bb}\le m\big\}\ \ \forall\ m\ge t_d
\eel
when $C_{q,d}^{1/q}\ge (2^{1/d+1/2})^d 2d 3^{d/2}$. 
In subsection A3.3 of the supplement, we extend the above inequality to all $m\ge 1$ 
and prove (\ref{th2-main}) by applying (\ref{th-1-2}) of Theorem \ref{th-1} with 
the above $\ell^*_{\pm}(m)$ and the $r_{q,\pm}(m)$ in \eqref{r-bd}. 

Theorem \ref{th-lattice-general} is a comprehensive statement which gives rise to many conclusions. 
In the next sub-subsection, 
we prove that the block estimator
is rate minimax in the $\ell_q$ risk for the entire lattice $[{\bf 1},\bn]$
in a wide range of configurations of $\bn$, $q$ and $\Delta^*_n$. In the next two subsections, we study the 
adaptation rate when $f(\cdot)$ is a piecewise constant function, and the variable selection rate when $f(\cdot)$ only depends on a subset of variables.

\subsubsection{Risk of the block estimator on the entire lattice and rate minimaxity} 
We assume without loss of generality in this sub-subsection $n_1\ge \cdots\ge n_d$. A direct comparison between Proposition~\ref{prop-lower-bd} and Theorem \ref{th-lattice-general} yields the following Theorem \ref{th-fix-worst}.

\begin{theorem}\label{th-fix-worst}
Let $\fhat_n^{(block)}(\bx)$ be the block estimator in (\ref {block})
with the lattice design $V=[{\bf 1},\bn]$ as in (\ref{lattice}).
Assume $\veps_i$ are independent random variables with $\E\, \veps_i =0$ and
$\E |\veps_i|^{q\vee 2} \le \sigma^{q\vee2}$.
Let $s_q = \lceil 2/(q-1) \rceil \wedge (d+1)$,
$n^*_s = \prod_{j=1}^s n_j$ for $s\le d+1$ with $n_{d+1}=1$, and
$\Delta(\bff_n/\sigma)=\{f(\bn)-f({\bf 1})\}/\sigma$.
Then, for $q\ge 1$,
\bel{th-fix-worst-0}
&& \sup\Big\{ R_q(\hbf_n^{(block)}, \bff_n):
\bff_n \in\scrF_n, \Delta(\bff_n/\sigma)\le\Delta_n^*\Big\}
\\ \nonumber &\lesssim_{q,d}&
\Lambda^{(match)} \inf_{\hbf}\, \sup_{\bff_n \in\scrF_n}
\Big\{ R_q(\hbf, \bff_n):
\Delta(\bff_n/\sigma)\le\Delta_n^*\Big\}
+ \frac{\sigma^q}{n}\Big(\prod_{j=1}^d \log_+(n_j)\Big)^{I\{q=2\}}
\eel
holds when $ \Delta_n^* \gtrsim_{q,d} t_{s_q}^{-1/2} = \big(n^*_{s_q}/n_{s_q}^{s_q}\big)^{-1/2}$,
where $\Lambda^{(match)} \le \log n$ is defined by 
\bel{lambda-match}&&
\Lambda^{(match)} = \bigg[\log_+ \Big( \min \Big\{\frac{n_{s_q}}{n_{s_q+1}},
\frac{n_{s_q}/(n_{s_q}^*)^{1/(s_q+2)}}{(\Delta_n^*)^{2/(s_q+2)}} \Big\} \Big)\bigg]^{I\{\frac{2}{q-1}=s_q \le d, \, \Delta_n^* \le n_{s_q}/t_{s_q}^{1/2}\}}.
\eel
Moreover, when $\max_{j \le d} n_j \lesssim_d n^{1/d}$ and $\Delta(\bff_n/\sigma)\le\Delta_n^*$,
\bel{th-fix-worst-1}
&& R_q (\hbf_n^{(block)} , \bff_n)
\\ \nonumber
&\lesssim_{q, d} & \sigma^q \,\min\Bigg\{1, \Big(\frac{\Delta_n^*}{n^{1/d}}\Big)^{\min\{1,\frac{qd}{d+2}\}} \bigg[\log_+ \Big(n \wedge \big(\frac{n^{1/d} }{ \Delta_n^*}\big)^{2d/(d+2)} \Big) \bigg]^{\delta_1}
+ \frac{(\log n)^{d\delta_2}}{n^{(q/2) \wedge 1}}\Bigg\},
\eel
holds for all $\Delta_n^* \ge 0$, where $\delta_1=I\{\frac{qd}{d+2}=1\}$ and $\delta_2=I\{q=2\}$.
\end{theorem}

\begin{remark}
It can be seen in our analysis that the logarithmic term presents for $q=2$, as the last component on 
the right-hand side of \eqref{th2-main}, \eqref{th-fix-worst-0} and \eqref {th-fix-worst-1}, 
due to the lack of data near the extreme points $\{\ba,\bb\}$ or $\{{\bf 1},{\bn}\}$ of the domain. 
\end{remark}

Compared with Proposition \ref{prop-lower-bd}, Theorem \ref{th-fix-worst} shows that
the risk of the block estimator matches the minimax rate when
$\Delta_n^* \ge t_{s_q}^{-1/2}= \big(\prod_{j=1}^{s_q}(n_j/n_{s_q})\big)^{-1/2}$ 
($\Delta_n^* \ge n^{-1/2}$ if $s_q = d+1$) 
possibly up to a logarithmic factor $\Lambda^{(match)} \le \log(n)$,
provided that the minimax rate is no faster than
${\sigma^q}{n^{-1}}\big(\prod_{j=1}^d \log_+(n_j)\big)^{\delta_2}$
due to the edge effect.
The match is always exact when $2/(q-1) \neq s_{q} \le d$, i.e., $2/(q-1)$ is not an integer or an integer greater than $d$.
When $2/(q-1) = s_q \le d-1$ and $n_{s_q}\asymp n_{s_q+1}$,
$\Lambda^{(match)}=O(1)$ and
the match is also exact. 
However, in the interesting setting where $q=d=2$ and $n_1\asymp n_2$, we have $s_q=2$
so that $\Lambda^{(match)} \asymp \log(n)$ when $\Delta_n^*\ll n_2$.

The one-dimensional risk bound for all $q \ge 1$ can be obtained from (\ref {th-fix-worst-1}) as
\bes
&& T_{q}([1,n])
\cr
& \lesssim_{q} 
& \sigma^q n \min \bigg\{1, \Big(\frac{\Delta_n^*}{n}\Big)^{\min\{q/3, 1\}} 
\Big[\log_+(n \wedge \big(\frac{n}{\Delta_n^*})^{2/3}\big)\Big]^{I\{q=3\}}   
+ \frac{(\log_+(n))^{I\{q=2\}}}{n^{ (q/2) \wedge 1}}\bigg\},
\ees
which reproduces (\ref {risk-bd-univ}) for $1\le q < 3$. We note that 
if we view one-dimensional isotonic regression as multi-dimensional on an $n_1\times 1\times \cdots\times 1$ lattice, the general bound yields this one-dimensional $n_1^{-1/3}$-rate. Interestingly, for general $\bn$, we still have the one-dimensional rate 
as long as the effective dimension $s$ is 0 or 1, i.e. $\Delta_n^*\ge n_2/t_2^{1/2} = n_2^{3/2}/n_1^{1/2}$. 
For $q=2$ and $d \ge 2$, it follows from Theorem \ref{th-lattice-general} that when $\Delta_n^* \ge n_2/t_2^{1/2} = n_2^{3/2}/n_1^{1/2} $, we have $s< s_q=2$ and only the first two cases of (\ref {th2-1}) are effective. This implies
\bes
T_2([{\bf 1}, \bn]) \lesssim_d \sigma^2 n \min\Big\{1, (\Delta_n^*/n_1)^{2/3} 
+ {\prod_{j=1}^d\big(\log_+(n_j)/{n_j}\big)}\Big\},
\ees
exactly the same as the bound of $T_2([1,n_1])$ in univariate case when $(\Delta_n^*/n_1)^{2/3}$ is dominant in both rates. 
In this case, our theory does not guarantee an advantage of the multiple isotonic regression on the entire lattice in terms of the $\ell_2$ risk, compared with the row-by-row univariate isotonic regression of length $n_1$. 
This observation agrees with \cite{chatterjee2018} where the $\ell_2$ minimax rate of two-dimensional isotonic regression, $\sigma^2 \Delta_n^*n^{-1/2}$, requires $n_2^{3/2}/n_1^{1/2} \ge \Delta_n^*$.

To conclude this subsection, we compare the $\ell_2$ risk bound for the block estimator in 
Theorem \ref {th-fix-worst} with those for the LSE in the existing literature.
For $d=2$, \cite{chatterjee2018}
gives an upper bound for the LSE as
\bes
R_2 (\hbf_n^{(lse)} , \bff_n) \lesssim  
\sigma^2 \Big(\frac{\Delta_n^*}{\sqrt{n}}(\log n)^4 + \frac{1}{n}(\log n)^8\Big),
\ees
for any $n_1 \times n_2$ lattice and $f$ satisfying $\Delta(\bff_n /\sigma) \le n_2^{3/2}/n_1^{1/2} $,
in contrast to
\bes
R_2 (\hbf_n^{(block)} , \bff_n) \lesssim 
\sigma^2 \Big(\frac{\Delta_n^*}{\sqrt{n}}\log(n) + \frac{1}{n}(\log n)^2\Big)
\ees
in (\ref {th-fix-worst-1}) of Theorem \ref{th-fix-worst} or in the third case of (\ref{th2-1}) of Theorem \ref{th-lattice-general} with $[\ba, \bb] = [{\bf 1}, \bn]$. 
However, for $n_1 = \cdots = n_d = n^{1/d}$ and $\Delta_n^* = 1$
as in \cite{han2017isotonic} for $d \ge 3$, 
(\ref {th-fix-worst-1}) is reduced to
\bes
R_2 (\hbf_n^{(block)} , \bff_n) \lesssim_{d} n^{-1/d},
\ees
which should be compared with the the rate
\bes
R_2 (\hbf_n^{(lse)} , \bff_n) \lesssim_{d} n^{-1/d}\log^4(n)
\ees
for the LSE \citep{han2017isotonic}.

\subsection{Adaptation rate of the block estimator with lattice designs in the piecewise constant case}

We consider here the adaptation behavior of the block estimator in the setting where $f(\cdot)$ is piecewise constant on a union of rectangles, as
a direct consequence of Theorem \ref{th-lattice-general}.

\begin{theorem}\label{coro-fix-piecewise}
Let $\fhat_n^{(block)}(\bx)$ be the block estimator in (\ref{block}). Assume $\veps_i$ are independent variables with $\E\, \veps_i =0$ and
$\E |\veps_i|^{q\vee 2} \le \sigma^{q\vee2}$ and $f$ is non-decreasing and piecewise constant on $V$ in the sense of
$V = \cup_{k=1}^K [\ba_k,\bb_k]$ with $K\le n$ and 
$f(\ba_k) = f(\bb_k)$ for all $k \le K$. Then, 
\bes
R_q (\hbf_n^{(block)}, \bff_n) \lesssim_{q, d} \sigma^q \min\bigg\{1, 
n^{-1}\sum_{k=1}^K n_{\ba_k, \bb_k}^{(1-q/2)_+} 
\Big(\log_+^{s_k}\big(n_{\ba_k,\bb_k}\big)\Big)^{I\{q=2\}}\bigg\}
\ees
with $s_k=\# \{j: b_{k,j} > a_{k,j} \}$. Moreover, if in addition $\{ [\ba_k,\bb_k], k=1, \ldots, K \big\}$ are disjoint, then
\bel{coro-2-2}
R_q (\hbf_n^{(block)}, \bff_n) \lesssim_{q, d} \sigma^q \min\Big\{1, \Big(\frac{K}{n}\Big)^{\min\{1, q/2\}}
\Big(\log^{d_K}_+(n/K)\Big)^{I\{q=2\}}
\Big\},
\eel
where $d_K = \max_{1 \le k \le K} s_k$ is the largest dimension of $[\ba_k,\bb_k]$ in the partition. 
\end{theorem}

The rate in (\ref {coro-2-2}) is consistent with existing results for $d=1$ under which the block estimator is the LSE and the mean squared risk bound is
\bes
R_2 (\hbf_n^{(block)}, \bff_n) \lesssim \sigma^2 \frac{K}{n}\log_+(n/K).
\ees
In general, the risk bound in (\ref {coro-2-2}) under $q=2$ is reduced to at most
\bes
\sigma^2 \frac{K}{n} \log_+^d(n/K),
\ees
which should be compared with
\bes
\sigma^2 \Big(\frac{K}{n}\Big)^{2/d} \log_+^8(n/K)
\ees
for the LSE as in
\cite{chatterjee2018} for $d=2$ and in \cite{han2017isotonic} for $d \ge 3$.

\begin{remark}
\cite{han2017isotonic} proved 
that even when $f(\cdot)$ is a constant function, i.e., $K=1$,
\bes
R_2(\hbf_n^{(lse)}, \bff_n) \gtrsim_d \sigma^2 n^{-2/d}
\ees
so the adaptation rate of the LSE, $(K/n)^{2/d}$, cannot be further improved, which means the LSE is unable to adapt to parametric rate for $d\ge 3$.
\end{remark}

The adaptation rate in (\ref {coro-2-2}) also implies that when $[\ba_k,\bb_k]$ are two-dimensional sheets (i.e. $|\{j: b_{k,j} \neq a_{k,j}\}|\le 2 $), the upper bound turns out to be
\bes
\frac{K}{n}\log_+^2(n/K),
\ees
which again should be compared with
\bes
\frac{K}{n}\log_+^8(n/K)
\ees
in \cite{han2017isotonic}.

\subsection{Adaptive estimation to variable selection with lattice designs}

In this subsection, we consider the case where
the true function of interest, $f(\cdot)$, depends only on a subset $S$ of $s$ variables, i.e., $f(\bx) = f_S(\bx_S)$.
We study the adaptive estimation when $\max_{j\le d} n_j \lesssim_d n^{1/d}$, i.e., $n_j \asymp n^{1/d}$ for all $1\le j \le d$.

\begin{theorem}\label{th-fix-selection}
Assume $f(\cdot)$ is non-decreasing and dependent only on an unknown set $S$ of $s < d$ variables. 
Let $\fhat_n^{(block)}(\bx)$ be the block estimator in (\ref {block})
on the lattice design $V=[{\bf 1},\bn]$. Assume $\max_{1\le j \le d} n_j \lesssim_d n^{1/d}$
and $\veps_i$'s are independent and satisfies $\E\, \veps_i =0$ and
$\E |\veps_i|^{q\vee 2} \le \sigma^{q\vee2}$.
Let $\Delta(\bff_n/\sigma)=\{f(\bn)-f({\bf 1})\}/\sigma$.
Then,
\bel{th-fix-selection-1}
&& \sup\Big\{ R_q(\hbf^{(block)}_n, \bff_n):
\bff_n \in\scrF_n, f(\bx) = f_S(\bx_S), \Delta(\bff_n/\sigma_S)\le \Delta_{n,S}^*\Big\}
\\
\nonumber &\lesssim_d&
\sigma_S^q
\min\bigg\{\Lambda_{s,1}^{(select)}, \Lambda_{s,2}^{(select)} \big(\Delta_{n,S}^*/{n^{1/d}}\big)^{\min\{1, \frac{qs}{s+2}\}}
\\
\nonumber  && \qquad\qquad\qquad + \, \Lambda_{s,1}^{(select)} \big(n^{s/d}\big)^{-\min\{1,  q/2\}} (\log n)^{sI\{q=2\}} \bigg\},
\eel
for all $1\le s\le d$, where $\sigma_S = \sigma / \big(\prod_{j\not\in S}n_j \big)^{1/2} \le C_d \sigma \big/n^{(1-s/d)/2}$ and
\bes
\Lambda_{s,1}^{(select)} &=& \Big(\sum_{j=1}^{n^{1/d}} j^{-q/2} \Big/ \big(n^{1/d}\big)^{1-q/2}\Big)^{d-s},
\cr
\Lambda_{s,2}^{(select)} &=& \Big(\sum_{j=1}^{n^{1/d}} j^{\min\{\frac{1-q}{2}, -\frac{q}{s+2} \}} \Big/ \big(n^{1/d}\big)^{\min\{\frac{1-q}{2}, -\frac{q}{s+2} \}+1}\Big)^{d-s}(\log n)^{I\{\frac{qs}{s+2}=1\}}.
\ees
In particular,
\bel{th-fix-selection-2}
&& R_2(\hbf_{n}^{(block)},\bff_n)
\\ \nonumber &\lesssim_d& \begin{cases}
\sigma^2n^{s/d-1}
\min\Big\{(\log n)^{d-s}, \Delta_{n,S}^* n^{-1/d}(\log n)^{I\{s=2\}}
+ n^{-s/d}(\log n)^d\Big\}, & s\ge 2,
\cr \sigma^2n^{s/d-1}
\min\Big\{(\log n)^{d-1}, (\Delta_{n,S}^*/n^{1/d})^{2/3} + n^{-1/d}(\log n)^d\Big\}, & s=1.
\end{cases}
\eel
\end{theorem}

In the proof of Theorem \ref{th-fix-selection}, the key observation is that in the sheet of $\bx$ with fixed $\bx_{S^c}$, the risk bound is identical
to that of model $S$ with $\sigma^q$ reduced by a factor of $n_{\bx_{S^c},\bn_{S^c}}^{-q/2}$.
The above rate would then become clear after the summation of risk bounds over $\bx_{S^c}$.

Let $n_j=n^{1/d}$ for all $j$. Consider an oracle expert with the extra knowledge of the subset $S$. 
Suppose the oracle expert first computes the average of the $n^{1-s/d}$ values of $y_i$ holding $\bx_S$ fixed and then solves the s-dimensional isotonic regression problem at the noise level 
$\sigma_S=\sigma n^{(s/d-1)/2}$. 
For this oracle expert, the sample size becomes $n^{s/d}$ 
and the condition on the range-to-noise ratio becomes 
$(f(\bn)-f({\bf 1}))/\sigma_S\le \Delta_{n,S}^*$, 
equivalent to $(f(\bn)-f({\bf 1}))/\sigma\le \Delta_n^*$ with $\Delta_{n,S}^* = \Delta_n^*n^{(1-s/d)/2}$. 
It follows from (\ref {prop-lower-bd-2}) in Proposition \ref{prop-lower-bd} that for $\veps_i\sim N(0,\sigma^2)$ and 
$ \Delta_{n,S}^* \ge (n^{-(s/d)/2})\vee(I\{q>1+2/s\})$, the $\ell_q$ minimax lower bound for the oracle expert is 
\bes
&& \inf_{\hbf}\, \sup\Big\{
R_q(\hbf, \bff_n):
\bff_n \in\scrF_n, f(\bx) = f_S(\bx_S), \Delta(\bff_n/\sigma_S) \le  \Delta_{n,S}^* \Big\}
\cr
&& \qquad \qquad\gtrsim \sigma_S^q
\min\Big\{1, \big( \Delta_{n,S}^*/n^{1/d}\big)^{\min\{1, qs/(s+2)\}} \Big\}.
\ees
Hence the variable-selection adaptation rate in (\ref{th-fix-selection-1}) matches 
the oracle minimax lower bound up to some constant or logarithmic factors $\Lambda_{s,1}^{(select)}$, $\Lambda_{s,2}^{(select)}$
and $\Lambda_{s,1}^{(select)}(\log n )^{sI\{q=2\}}$, provided that 
\bes
 \Delta_{n,S}^* \ge \max\Big(n^{-s/(2d)}, I\{ q > 1 + 2/s\}\Big),
\ees
or equivalently $\Delta(\bff_n/\sigma)\le\Delta_n^*$ with 
$\Delta_n^*\ge \max\big(n^{-1/2}, n^{- (1-s/d)/2} I\{q>1+2/s\}\big)$.
The match to the oracle minimax rate is always exact for $q=1$ and any $s$ as both $\Lambda_{s,1}^{(select)}$ and $\Lambda_{s,2}^{(select)}$ are bounded by a constant.
When $q=2$, the match is also exact but up to some logarithmic factors as $\Lambda_{s,1}^{(select)}  \lesssim_d (\log n)^{d-s}$ and $\Lambda_{s,2}^{(select)}  \lesssim_d (\log n)^{I\{s=2\}}$.

\subsection{Multiple isotonic regression with random designs}

In this subsection we consider $V = [{\bf 0}, {\bf 1}]$ in continuum and, same as before, $\ba \preceq \bb$ iff $\ba\le\bb$. 
Different from fixed designs, here $\bx_1,\ldots, \bx_n$ are i.i.d. random vectors from a distribution $\P$ supported on $[{\bf 0}, {\bf 1}]$.
For simplicity we assume the distribution of the design points 
has a Lebesgue density bounded both from above and below; for $\mu_{\bu,\bv}=\P\{\bu\le \bx_i \le\bv\}$ 
and the Lebesgue $\mu^L_{\bu,\bv}=\mu^L ([\bu, \bv]) = \int_{[\bu,\bv]}d\bx$,
\bel{P-condition} 
\rho_1 \mu^L_{\bu,\bv} \le \mu_{\bu,\bv} \le \rho_2 \mu^L_{\bu,\bv}.
\eel
with certain fixed constants $0< \rho_1 \le \rho_2 <\infty$. 
We consider the integrated $L_q$ risk in (\ref {integrated-risk}), i.e.,
\bes
R_q^* (\fhat_n^{(block)}, f) = \int_{\bx \in [{\bf 0}, {\bf 1}]} \E \big|\fhat_n^{(block)}(\bx) - f(\bx)\big|^q d\bx,
\ees
and partial integrated $L_q$ risk on block $[\ba, \bb]$ as
\bes
R^*_{q}([\ba,\bb])
&=&\int_{[\ba,\bb]} \E \big|\fhat_n^{(block)}(\bx) - f(\bx)\big|^q d\bx.
\ees

\begin{theorem}\label{th-random}
Let $\fhat_n^{(block)}(\bx)$ be the block estimator in (\ref {block})
with $V=[{\bf 0},{\bf 1}]$.
Assume $\bx_1,\ldots, \bx_n \in [{\bf 0}, {\bf 1}]$ are i.i.d. random vectors drawn from a distribution satisfying (\ref {P-condition}).
Assume $f$ is non-decreasing and $\veps_i$ are independent random variables with $\E\, \veps_i =0$ and
$\E |\veps_i|^{q\vee 2} \le \sigma^{q\vee2}$.
Let $\{\ba,\bb\}\subset V$ with $\ba\le\bb$.
Then, for $q\ge 1$,
\bel{th-random-main}
R^*_{q}([\ba,\bb])
&=&\int_{[\ba,\bb]} \E \big|\fhat_n^{(block)}(\bx) - f(\bx)\big|^q d\bx
\\
\nonumber &\le& C_{q,d,\rho_1,\rho_2}^* \sigma^q \bigg[ \int_0^{n\mu_{\ba, \bb}} \Big( (t \vee 1)^{-q/2} + \Delta_{\ba, \bb}^q e^{-t}\Big) H^*(dt)
\cr
&& \quad + \, \int_{\bx \in [\ba, \bb]} 
\Big(\big\{(n \mu_{\bx, \bb} )\vee 1\big\}^{-q/2} + \Delta_{{\bf 0}, {\bf 1}}^q  e^{-n\mu_{\bx, \bb}}\Big) d\bx \bigg],
\eel
where $\Delta_{\bu,\bv} = \big(f(\bv) - f(\bu)\big)\big/\sigma$ and 
$\mu_{\bu,\bv}= \P \big\{ \bx_i \in [\bu, \bv] \big\}$ for all $\bu\le\bv$ 
and $H^*(t) = \min\big\{1, \Delta_{\ba, \bb}(n\mu_{\ba, \bb})^{-1/d}t^{1/2+1/d}\big\}$. 
Specifically, (\ref {th-random-main}) is no greater than
\bel{th-random-rate}
&\sigma^q \min \bigg\{(\Delta_{{\bf 0},{\bf 1}}^q+1)\mu_{\ba, \bb}, &
\Big(\frac{\Delta_{\ba,\bb}}{(n\mu_{\ba, \bb})^{1/d}}\Big)^{\min\{1, \frac{qd}{d+2}\}} \Lambda_1^{(random)}
\\
\nonumber && \, + \, \frac{\Delta_{\ba, \bb}^{q+1}}{(n\mu_{\ba, \bb})^{1/d}} 
+ \big(\Delta^q_{{\bf 0}, {\bf 1}} +1\big) \mu_{\ba, \bb} 
\frac{ \Lambda_2^{(random)} 
}{(n\mu_{\ba, \bb})^{(q/2) \wedge 1}}  \bigg\}
\eel
up to a constant depending on $q,d,\rho_1,\rho_2$ only, where 
\bes
\Lambda_1^{(random)} = \Big[\log_+\Big(n\mu_{\ba, \bb} \wedge \Big((n\mu_{\ba, \bb})^{\frac{2}{d+2}}\Big/
\Delta_{\ba, \bb}^{{2d}/{(d+2)}}\Big)\Big)\Big]^{I\{\frac{qd}{d+2}=1\}}
\ees
and $\Lambda_2^{(random)} = \big(\log_+(n \mu_{\ba, \bb} )\big)^{dI\{q = 2\} + (d-1)I\{q>2\}}$. 
\end{theorem}

The $H^*(t)$ here is identical to the $\Htil(t)$ in Theorem \ref{th-lattice-general} in $t\in [t_d,n]$, 
effectively taking $t_d=1$.
This reveals an intrinsic difference between lattice design and random design: the effective dimension of the random design over $[\ba, \bb] \subseteq [{\bf 0},{\bf 1}]$ is always $d$ --- any hyper-rectangle $[\ba, \bb]$ with positive measure behaves similarly to a hyper-cube.
The above rate in (\ref{th-random-rate}) is therefore comparable to the rate in (\ref {th-fix-worst-1}) for the lattice design with $n_j = n^{1/d}$ for all $j$.
In fact, the rate in (\ref {th-random-rate}) can be derived from a scale change of the upper bound for $R^*_q([{\bf 0},{\bf 1}])$.

The study of the integrated $L_q$ risk in isotonic regression is relatively new.
\cite{fokianos2017integrated} gives an asymptotic bound, $O(n^{-1/(d+2)})$, for the $L_1$ risk with 
$[\ba, \bb] = [{\bf 0}, {\bf 1}]$. The $L_1$ error bound in Theorem \ref{th-random} is consistent with their result. 

To fit in with random design, we now define $r_{q,+}(m)$ as a non-increasing function of $m\in [0, n]$ in continuum satisfying
\bel{r-plus-new}&&
r_{q,+}(m) \ge \max\bigg\{
\E\bigg(\max_{\bu \preceq \bx} \sum_{\bx_i\in [\bu,\bv]} \frac{\veps_i}{n_{\bu,\bv} \vee 1}\bigg)_+^q: 
\E[n_{\bx, \bv}] = m , \bx \preceq \bv \hbox{ and } \bv \in V_0\bigg\},
\eel
and modify the definition of $m_{\bx} = m_{\bx, +}$ in (\ref {k-m}) to
\bel{m_x-random}
m_{\bx} 
= n\mu_{\bx, \bv_{\bx}}, \hbox{ where } \bv_{\bx} &=& \argsup_{\bx \le \bv \le \bb}\big\{n\mu_{\bx, \bv}: f(\bv) \le f(\bx) + r_{q,+}^{1/q}(n\mu_{\bx, \bv})\big\}.
\eel
Note $n_{\bx, \bv}$, the number of design points in $[\bx, \bv]$, becomes a
Binomial$(n, \mu_{\bx, \bv})$
random variable.
Here we omit $m_{\bx, -}$ as it can be analyzed by symmetry.
Nevertheless, Theorem \ref{th-random} is still proved in a similar way to Theorem \ref{th-lattice-general}. However, different from (\ref {th-1-1}) in Theorem \ref{th-1}, the point risk bound is given by the following proposition.

\begin{proposition}\label{prop-point-risk}
Assume the conditions of Theorem \ref{th-random}. Then, \eqref{r-plus-new} holds for  
\bel{r-plus-random} 
r_{q, +}(n\mu_{\bx, \bv}) = C_{q,d,\rho_1, \rho_2} \sigma^q(n\mu_{\bx, \bv} \vee 1)^{-q/2} 
\eel
with $C_{q,d,\rho_1, \rho_2}$ continuous in $q\in [1,\infty)$ and for all $\bx \in [\ba, \bb]$
\bel{point-risk}
&& \E \big(\fhat^{(block)}_n(\bx) - f(\bx) \big)_{+}^q
\\ \nonumber &\le& 2^{q} r_{q,+}(m_{\bx})  + 2^{q-1}\sigma^q C_{q, d, \rho_1, \rho_2} \Big( \big(\Delta_{\ba, \bb}^q + 1\big) e^{-m_{\bx}} + \big(\Delta_{{\bf 0}, {\bf 1}}^q + 1\big) e^{-n\mu_{\bx, \bb}}\Big).
\eel
\end{proposition}

As we discussed below (\ref {k-m}), the positive part of the bias of $\fhat_n^{(block)}(\bx)$ is of no greater order than the variability of the noise as measured by $r_{q,+}^{1/q}(n_{\bx, \bv_{\bx}})
\asymp r_{q,+}^{1/q}(m_{\bx})$ provided the presence of at least one design point in $[\bx, \bv_{\bx}]$.
The first term on the right-hand side of (\ref {point-risk}) thus comes from the case of $n_{\bx, \bv_{\bx}} >0$.
However, $[\bx, \bv_{\bx}] $ might be an empty cell with no design points.
We then have to consider points in $[\bx, \bb]$ when $n_{\bx,\bv_{\bx}}=0$ and in $[\bx,{\bf 1}]$ when 
$n_{\bx,\bb}=0$, leading to terms with $\Delta_{\ba, \bb}$ and $\Delta_{{\bf 0}, {\bf 1}}$ respectively. 

\medskip
Corresponding to Theorem \ref{th-fix-worst} and \ref{coro-fix-piecewise}, the following two theorems give the risk bounds for random designs under the general case and the piecewise constant case for the entire $[{\bf 0}, {\bf 1}]$. Due to space limitations, the minimax rate and the adaptation rate to variable selection in random design are not discussed. 

\begin{theorem}\label{coro-random-worst}
Let $\fhat_n^{(block)}(\bx)$, $f$ and $\{\bx_i, \veps_i, i\le n\}$ be as in Theorem \ref{th-random}. 
Suppose $\Delta_{{\bf 0},{\bf 1}} = \big(f({\bf 1}) - f({\bf 0})\big)\big/\sigma$ is bounded by a constant. Then
\bes
R_q^* (\fhat_n^{(block)}, f)
&\lesssim_{q,d,\rho_1, \rho_2}& \sigma^q \Big(\frac{\Delta_{{\bf 0}, {\bf 1}}}{n^{1/d}}\Big)^{\min\{1, \frac{qd}{d+2}\}}\big(\log n\big)^{I\{\frac{qd}{d+2}=1\}}
\cr
&& \qquad \qquad + \,  \frac{ \big(\log n\big)^{dI\{q = 2\} + (d-1)I\{q>2\}} }{n^{(q/2) \wedge 1}}.
\ees
In particular when $q=2$ and $d \ge 2$,
\bel{th-random-worst-2}
R_2^* (\fhat_n^{(block)}, f) \lesssim_{d,\rho_1, \rho_2} \sigma^2 \min \bigg\{1, \frac{\Delta_{{\bf 0}, {\bf 1}}}{n^{1/d}}\big(\log n\big)^{I\{d=2\}} + \frac{ \big(\log n\big)^{d}}{n}  \bigg\}.
\eel
\end{theorem}

\begin{remark}
For simplicity, we here consider the case of bounded $\Delta_{{\bf 0}, {\bf 1}}$. 
Theorem \ref{th-random} also directly yields error bounds for 
general $\Delta_{{\bf 0}, {\bf 1}}$ by setting $[\ba, \bb]=[{\bf 0}, {\bf 1}]$ in (\ref {th-random-main}) and (\ref {th-random-rate}). 
\end{remark}

\begin{theorem}\label{coro-random-piecewise}
Let $\fhat_n^{(block)}(\bx)$, 
$f$ and $\{\bx_i, \veps_i, i\le n\}$ be as in Theorem \ref{th-random}. 
Suppose $V$ has disjoint partition $V=\cup_{k=1}^K [\ba_k,\bb_k]$ with $K \le n$ and $f(\ba_k) = f(\bb_k)$ for all $k \le K$. Then
\bel{coro-3-2}
&& R^*_q (\fhat_n^{(block)}, f)
\\
\nonumber &\lesssim_{q,d,\rho_1,\rho_2}& \sigma^q (\Delta_{{\bf 0}, {\bf 1}}^q +1) \Big(\frac{K}{n}\Big)^{\min\{1, q/2\}}\big(\log_+(n/K)\big)^{dI\{q=2\} + (d-1)I\{q>2\}},
\eel
where $\Delta_{{\bf 0},{\bf 1}} = \big(f({\bf 1}) - f({\bf 0})\big)\big/\sigma$.
In particular, when $q=2$,
\bes
R^*_2 (\fhat_n^{(block)}, f) \lesssim_{d,\rho_1,\rho_2} \sigma^2 (\Delta_{{\bf 0}, {\bf 1}}^2 +1) \frac{K}{n}\log_+^{d}(n/K).
\ees
\end{theorem}

\medskip
We can also derive risk bounds for the empirical $\ell_q$ risk.
As $[\bx_i, \bv_{\bx_i}]$ always has
the
design point $\bx_i$, there is no ``empty cell'' problem as in Proposition \ref{prop-point-risk} when bounding the empirical risk.
It follows that
\bes
\E \big[ \big(\fhat_n^{(block)}(\bx_i) - f(\bx_i) \big)_{+}^q \big| \bx_i = \bx\big] \lesssim_{q,d, \rho_1, \rho_2} r_{q,+}(m_{\bx}),
\ees
so that
\bes
&& R_q(\hbf^{(block)}_n, \bff)
\cr
&\lesssim_{q,d,\rho_1, \rho_2}& \sigma^q \min \bigg\{\mu_{\ba, \bb}, \Big(\frac{\Delta_{\ba,\bb}}{(n\mu_{\ba, \bb})^{1/d}}\Big)^{\min\{1, \frac{qd}{d+2}\}} \Lambda_1^{(random)}
+ \mu_{\ba, \bb} \frac{ \Lambda_1^{(random)}
}{(n\mu_{\ba, \bb})^{(q/2) \wedge 1}}  \bigg\}
\ees
by an almost identical proof. 
It follows that under the conditions of Theorem \ref{th-random} and $\Delta_{{\bf 0}, {\bf 1}}=1$, the worst case upper bound of the mean squared risk is
\bes
R_2 (\hbf^{(block)}_n, \bff) \lesssim_{d,\rho_1, \rho_2} \sigma^2 n^{-1/d} (\log n)^{I\{d=2\}},
\ees
and under the conditions of Theorem \ref{coro-random-piecewise}, the mean squared risk bound in piecewise constant case is
\bes
R_2 (\hbf^{(block)}_n, \bff) \lesssim_{d,\rho_1,\rho_2} \sigma^2 \frac{K}{n}\log^{d}(n/K).
\ees

We shall compare the above two rates with the results for the LSE in
\cite{han2017isotonic} respectively, i.e.,
\bes
\sigma^2n^{-1/d}\log^{\gamma_d}(n)
\ees
and
\bes
\sigma^2\Big(\frac{K}{n}\Big)^{2/d} \log^{2\gamma_d}(en/K),
\ees
where $\gamma_2 = 9/2$ and $\gamma_d = (d^2+d+1)/2$ when $d\ge 3$. It is worth mentioning that \cite{han2017isotonic}
also proved the piecewise constant rate for the LSE, $(K/n)^{2/d}$, is not improvable as when $K=1$,
\bes
R_{2}(\hbf_n^{(lse)}(\bx), \bff) \gtrsim_{d,\rho_1, \rho_2} \sigma^2 n^{-2/d}.
\ees

\subsection{Model misspecification}
We consider in this subsection properties of the block estimator in the nonparametric regression model 
\bel{np-reg}
y_i = f(\bx_i) + \veps_i,\ i=1,\ldots, n,
\eel
for general $f$. When the true regression function $f$ fails to be non-decreasing, 
the isotonic regression model \eqref{iso-reg} is misspecified, 
so that the block estimators actually estimate their noiseless versions, say $\fbar_n^*(\bx)$, 
instead of the true $f$. 
For the block max-min and min-max estimator in \eqref{max-min-min-max},    
\bel{max-min-min-max-mis}
\fbar_n^*(\bx) 
= \fbar_n^{(max-min)}(\bx) &=& \max_{\bu\preceq \bx, n_{\bu,*}>0}\ \min_{\bx\preceq\bv, n_{\bu,\bv}>0}
\fbar_{[\bu,\bv]},\quad\forall\ \bx\in V,
\\ \nonumber 
\fbar_n^*(\bx) = 
\fbar_n^{(min-max)}(\bx) &=&\min_{\bx\preceq\bv, n_{*,\bv}>0}\ \max_{\bu\preceq \bx, n_{\bu,\bv}>0}
\fbar_{[\bu,\bv]},\quad\forall\ \bx\in V,
\eel
are their noiseless versions, where $\fbar_A$ denotes the average of $\{f(\bx_i): 1\le i\le n, \bx_i\in A\}$.  
For the average \eqref{block-mid} of the two estimators, the noiseless version is 
\bel{block-mis}
\fbar_n^*(\bx) = \frac{1}{2}
\Big\{\fbar_n^{(max-min)}(\bx) + \fbar_n^{(min-max)}(\bx)\Big\},
\quad \forall\ \bx \in V. 
\eel
The functions in \eqref{max-min-min-max-mis} and \eqref{block-mis} can be viewed as estimation targets. 

Our results can be summarized as follows.  
If we treat $\fhat^{(block)}_n(\bx) - \fbar_n^*(\bx)$ as the estimation error 
and use $\fbar_n^*/\sigma$ to measure the range-to-noise ratio, 
all the theoretical results we have presented so far 
hold in the nonparametric regression model \eqref{np-reg} for general $f$ with 
the following adjustments of the error bounds $r_{q,\pm}(m)$ 
in \eqref{r-plus} and \eqref{r-minus}, 
\bel{r-plus-mis}\label{r-minus-mis}\qquad
r_{q,+}(m) &\ge& \max\bigg\{
\E\bigg[\max_{\bv'\succeq\bv}
\bigg(\max_{\bu \preceq \bx} \sum_{\bx_i\in [\bu,\bv']} \frac{\veps_i}{n_{\bu,\bv'}}\bigg)_+^q\bigg]:
n_{\bx,\bv}=m, \bx \preceq \bv \hbox{ and } \bv \in V_0\bigg\}, 
\\ \nonumber 
r_{q,-}(m) &\ge&  \max\bigg\{
\E\bigg[\max_{\bu'\preceq\bu}
\bigg(\min_{\bv\succeq \bx}\sum_{\bx_i\in [\bu',\bv]} \frac{\veps_i}{n_{\bu',\bv}}\bigg)_-^q\bigg]:
n_{\bu,\bx}=m, \bu \preceq \bx \hbox{ and } \bu \in V_0 \bigg\}, 
\eel
without changing the notation. 
Both $r_{q,\pm}(m)$ are still required to be 
non-increasing functions of $m\in \N^+$.  
Accordingly, this leads to the following adjustment of the functions in \eqref{k-m}, 
\bel{k-m-mis}
m_{\bx,-}&=&  \max \Big\{n_{\bu,\bx}: \fbar_n^*(\bu) \ge \fbar_n^*(\bx) - r_{q,-}^{1/q}(n_{\bu,\bx}),
\bu \preceq \bx \hbox{ and } \bu \in V_0  \Big\},
\cr \bu_{\bx} &=& \argmax_{\bu\in V_0 :\, \bu \preceq \bx} \Big\{n_{\bu,\bx}:
\fbar_n^*(\bu) \ge \fbar_n^*(\bx) - r_{q,-}^{1/q}(n_{\bu,\bx})\Big\},
\\ \nonumber m_{\bx} = m_{\bx,+} &=& \max \Big\{n_{\bx,\bv}: \fbar_n^*(\bv) \le \fbar_n^*(\bx)+r_{q,+}^{1/q}(n_{\bx,\bv}),
\bx \preceq \bv \hbox{ and } \bv \in V_0 \Big\},
\\ \nonumber \bv_{\bx} &=& \argmax_{\bv\in V_0:\, \bx \preceq \bv} \Big\{n_{\bx,\bv}:
\fbar_n^*(\bv) \le \fbar_n^*(\bx)+r_{q,+}^{1/q}(n_{\bx,\bv})\Big\}, 
\eel
with the error bounds $r_{q,\pm}(m)$ in \eqref{r-plus-mis} and 
the estimation target $\fbar_n^*(\bx)$ in \eqref{max-min-min-max-mis} or \eqref{block-mis}.   

\begin{theorem}\label{th-model-mis} 
Let $\fhat^{(block)}_n$ be as in \eqref{block-mid}, $\fbar_n^*$ as in \eqref{block-mis}, 
$r_{q,\pm}(m)$ as in \eqref{r-plus-mis}, and $\ell_{\pm}(m)$ as in \eqref{ell-pm} with the 
$m_{\bx_i,\pm}$ in \eqref{k-m-mis}. 
Then, the error bounds \eqref{th-1-1} and \eqref{th-1-2} of Theorem \ref{th-1} hold 
with $f$ replaced by $\fbar_n^*$. 
Consequently, for the lattice design and under the $q\vee 2$ moment assumption on the noise $\{\veps_i\}$, 
the error bounds in Theorems \ref{th-lattice-general}, \ref{th-fix-worst}, \ref{coro-fix-piecewise} 
and \ref{th-fix-selection} hold with the same substitution. 
In particular, with $f$ replaced by $\fbar_n^*$ and $\bff_n$ by 
${\overline \bff}_n^* = (\fbar_n^*(\bx_1),\ldots,\fbar_n^*(\bx_n))^T$, 
\eqref{th2-main} holds with the same function $\Htil(t)$ 
when $\{\fbar_n^*(\bb)-\fbar_n^*(\ba)\}/\sigma \le\Delta_n^*$, 
\eqref{th-fix-worst-0} and \eqref{th-fix-worst-1} hold when $\Delta(\fbar_n^*/n)\le \Delta_n^*$, 
\eqref{coro-2-2} holds when $\fbar_n^*(\ba_k) = \fbar_n^*(\bb_k)$ with 
$V = \cup_{k=1}^K [\ba_k,\bb_k]$, 
and \eqref{th-fix-selection-1} and \eqref{th-fix-selection-2} hold when 
$\fbar_n^*(\bx)$ depends on only $s$ of the $d$ variables and $n_j\asymp n^{1/d}$ for all $j$. 
The above results also hold when 
$\{\fhat^{(block)}_n, \fbar_n^*\}=\{\fhat^{(max-min)}_n, \fbar_n^{(max-min)}\}$ 
or $\{\fhat^{(block)}_n, \fbar_n^*\}=\{\fhat^{(min-max)}_n, \fbar_n^{(min-max)}\}$. 
\end{theorem} 

Theorem \ref{th-model-mis} asserts that $\hbf_n$ is close to ${\overline \bff}_n^*$ in many ways 
when the isotonic condition on the unknown $f$ is misspecified. 
However, the interpretation of this result is not as clear as the existing oracle inequality for 
the LSE as ${\overline \bff}_n^*$ is not based on an optimality criterion. 

\section{Simulation results}

In this section, we report the results of several experiments in $d=2$ and $d=3$ to demonstrate the feasibility of the block estimators and to 
compare its estimation performance with the LSE. 
Among potentially many choices of the block estimator, we simply use the block max-min estimator as in (\ref {max-min-min-max}). 
In six simulation settings, the block max-min estimator yields smaller average $\ell_2$ losses than the LSE, with very small $p$-values in piecewise constant and variable selection settings. 
In a seventh setting, the LSE slightly outperforms the block max-min estimator 
but the difference is insignificant.

To compare the LSE and the block estimator, we carry out our experiments as follows. In each experiment, we generate one unknown $\bff$, $5000$ replications of $\by$ with standard Gaussian noise, find the LSE and the block max-min estimator for each $\by$, and compute the mean squared losses $\|\hbf_n-\bff_n\|_2^2/n$ for both estimators. 
We therefore obtain $5000$ simulated losses for each estimator 
and take the averages to approximate their mean squared risks.

We use quadratic programming to compute the LSE in our experiments. We'd like to mention that fast algorithms for the LSE have been developed in the literature: \cite{dykstra1983algorithm}, \cite{kyng2015fast}, \cite{Stout2015}, to name a few. We stick to quadratic programming as it provides somewhat more accurate results, although the difference seems small. The purpose of our experiment is to compare the risk of estimators, not the computational complexity of different algorithms. For the block max-min estimator, we use brute force which exhaustively calculates means over all blocks and finds the max-min value for each lattice point $\bx$. We note again that the computation cost via brute force is of order $n^3$.

\begin{figure}
  \centering
  \subfigure[true $f$ (unknown)]{
    \label{simu-1:a}
    \includegraphics[width=0.48\textwidth]{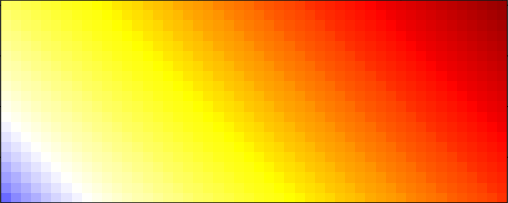}}
  \hspace{0.01\textwidth}
  \subfigure[$\by$ (observed)]{
    \label{simu-1:b}
    \includegraphics[width=0.48\textwidth]{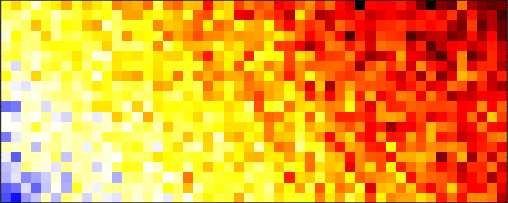}}
    \subfigure[the LSE]{
    \label{simu-1:c}
    \includegraphics[width=0.48\textwidth]{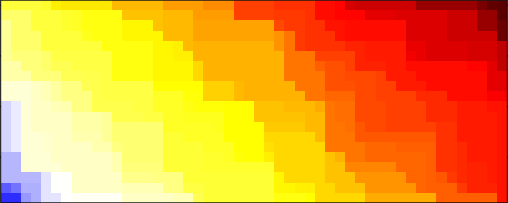}}
     \hspace{0.01\textwidth}
     \subfigure[the block max-min estimator]{
    \label{simu-1:d}
    \includegraphics[width=0.48\textwidth]{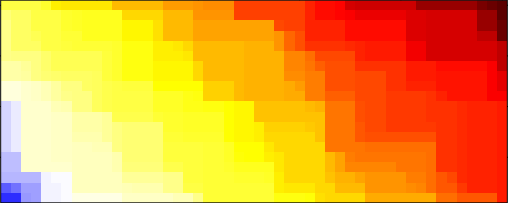}}

  \subfigure{
    \includegraphics[width=0.5\textwidth]{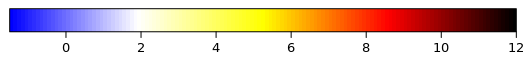}}

\caption{Heatmaps for the true $\bff$, an observed $\by$, and its LSE and max-min estimate in Experiment I.}
\label{simu-1}
\end{figure}

\begin{figure}
  \centering
  \subfigure[true $f$ (unknown)]{
    \label{simu-2:a}
    \includegraphics[width=0.48\textwidth]{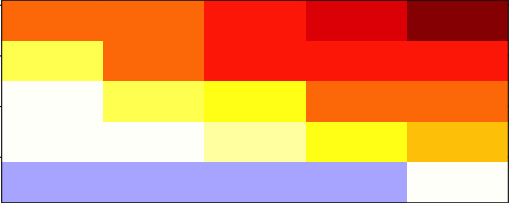}}
  \hspace{0.01\textwidth}
  \subfigure[$\by$ (observed)]{
    \label{simu-2:b}
    \includegraphics[width=0.48\textwidth]{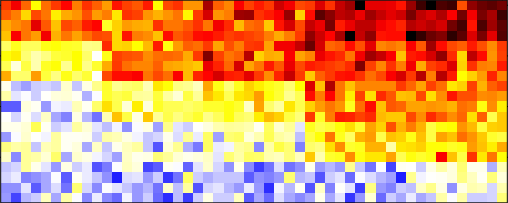}}
    \subfigure[the LSE]{
    \label{simu-2:c}
    \includegraphics[width=0.48\textwidth]{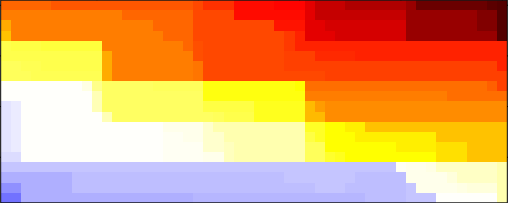}}
     \hspace{0.01\textwidth}
     \subfigure[the block max-min estimator]{
    \label{simu-2:d}
    \includegraphics[width=0.48\textwidth]{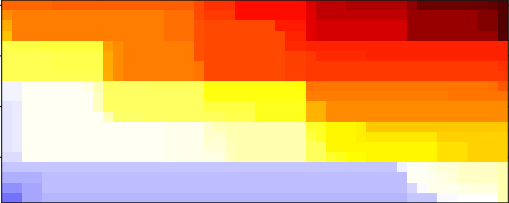}}

  \subfigure{
    \includegraphics[width=0.5\textwidth]{scale.png}}

\caption{Heatmaps for the true piecewise-constant $\bff$, an observed $\by$, and its LSE and max-min estimate in Experiment II.}
\label{simu-2}

\end{figure}

\begin{figure}
  \centering
  \subfigure[true $f$ (unknown)]{
    \label{simu-3:a}
    \includegraphics[width=0.48\textwidth]{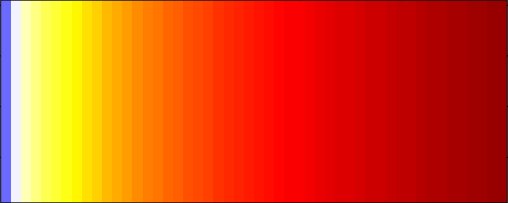}}
  \hspace{0.01\textwidth}
  \subfigure[$\by$ (observed)]{
    \label{simu-3:b}
    \includegraphics[width=0.48\textwidth]{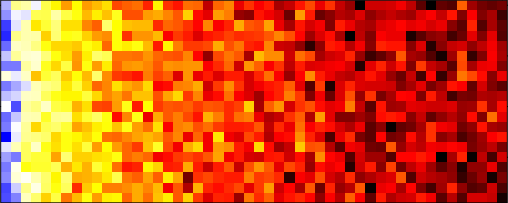}}
    \subfigure[the LSE]{
    \label{simu-3:c}
    \includegraphics[width=0.48\textwidth]{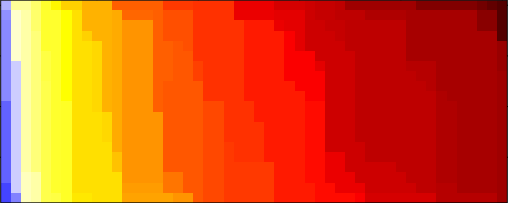}}
     \hspace{0.01\textwidth}
     \subfigure[the block max-min estimator]{
    \label{simu-3:d}
    \includegraphics[width=0.48\textwidth]{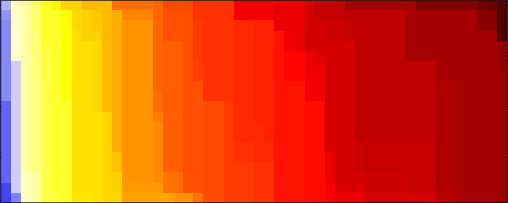}}

  \subfigure{
    \includegraphics[width=0.5\textwidth]{scale.png}}

\caption{Heatmaps for the true $\bff$, an observed $\by$, and its LSE and max-min estimate in Experiment III.}
\label{simu-3}
\end{figure}

\begin{figure}
  \centering
  \subfigure[Experiment I]{
    \label{boxplot-2d:a}
    \includegraphics[width=0.31\textwidth]{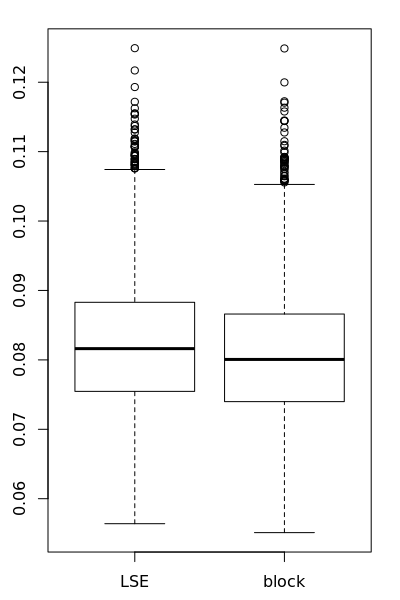}}
  \hspace{0.01\textwidth}
  \subfigure[Experiment II]{
    \label{boxplot-2d:b}
    \includegraphics[width=0.31\textwidth]{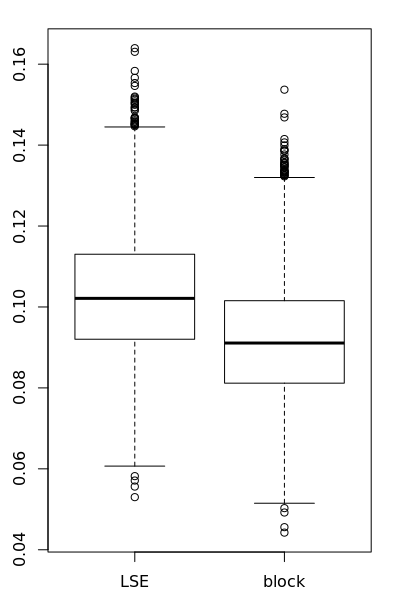}}
      \hspace{0.01\textwidth}
    \subfigure[Experiment III]{
    \label{boxplot-2d:c}
    \includegraphics[width=0.31\textwidth]{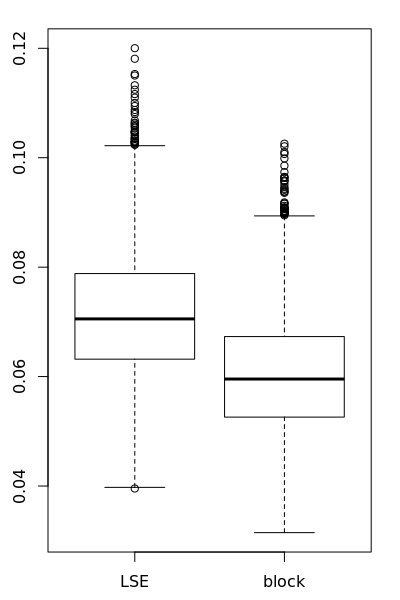}}
\caption{Boxplots for the losses of LSE and block estimator in $d=2$.}
\label{boxplot-2d}
\end{figure}

In $d=2$, we consider isotonic regression with the $n_1 \times n_2$ lattice design $[{\bf 1},\bn]$ with 
$n_1 = 50$ and $n_2 = 20$, so that the number of design points in total is $n = 1000$. 
In Experiment I, we consider the function $f(\bx) = c(x_1 + x_2)^{2/3}$ (here and in the sequel, $c$ is a constant such that $f(\bn)=10$ so that the range of $f$ is about 10 on the lattice). 
As the region between two contours of this $f$ cannot be efficiently represented by rectangular bocks, this example is not expected to favor the block estimator.
In Experiment II, we split the lattice into $5 \times 5$ small blocks of size $10 \times 4$, randomly assign $1, \ldots, 10$ to each small block, conditionally on the realizations satisfying the isotonic constraint. The adaptation of the LSE and the block max-min estimator to piecewise constant $f$ is compared in this experiment. Lastly, we compare the adaptation of the two estimators to variable selection in Experiment III by setting $f(\bx) = f_1(x_1) = c\log(x_1)$. See Figure \ref{simu-1}, \ref{simu-2} and \ref{simu-3} for heat maps in Experiment I, II and III respectively; each figure contains heat maps for the unknown $f$, one example of observed $\by$, the LSE and the block max-min estimator for this $\by$. 
Figure \ref{boxplot-2d} provides boxplots of mean squared losses of both estimators in Experiment I, II and III.

\begin{figure}
  \centering
  \subfigure[Experiment IV]{
    \label{boxplot-3d:a}
    \includegraphics[width=0.31\textwidth]{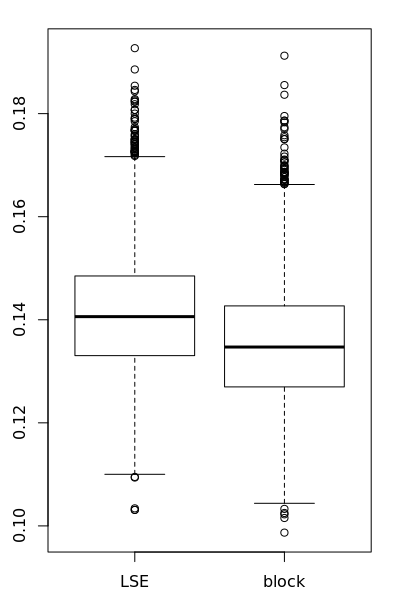}}
  \hspace{0.01\textwidth}
  \subfigure[Experiment V]{
    \label{boxplot-3d:b}
    \includegraphics[width=0.31\textwidth]{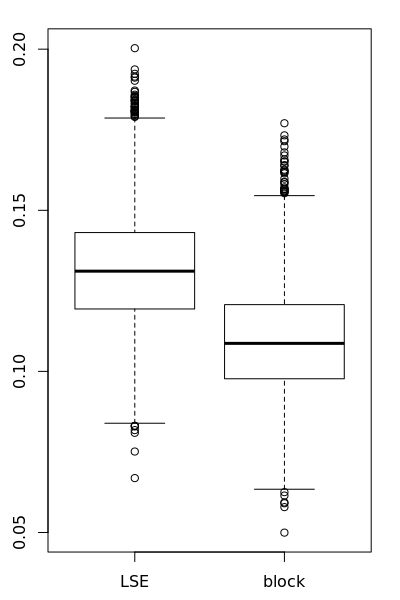}}
      \hspace{0.01\textwidth}
    \subfigure[Experiment VI]{
    \label{boxplot-3d:c}
    \includegraphics[width=0.31\textwidth]{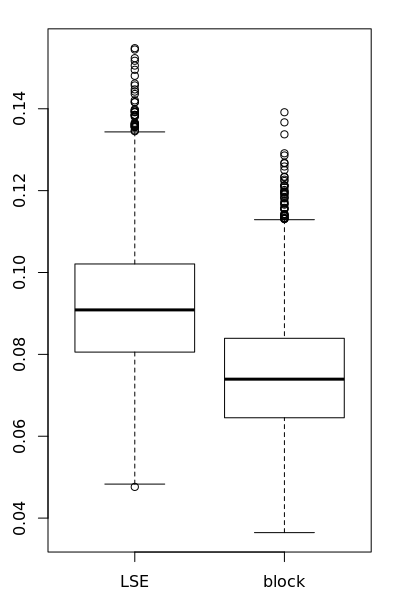}}
\caption{Boxplots for the losses of LSE and block estimator in $d=3$.}
\label{boxplot-3d}
\end{figure}

In $d=3$, we consider isotonic regression with $n_1 \times n_2 \times n_3$ lattice designs where $n_1 = n_2 = n_3 =10$, so that the number of design points in total is also $n = 1000$. We choose the true mean functions in a similar manner to $d=2$. In Experiment IV, we consider $f(\bx) = c(x_1 + x_2 + x_3)^{2/3}$. In Experiment V, we randomly assign $1, \ldots, 10$ to $2 \times 2 \times 5$ small blocks of size $5 \times 5 \times 2$ conditionally on the isotonic constraint. Lastly, the true mean function is $f(\bx) = f_1(x_1) = c\log(x_1)$ in Experiment VI. See Figure \ref{boxplot-3d} for boxplots of mean squared losses of both estimators in Experiment IV, V and VI.

Two basic statistics, mean and standard deviation of the losses of 
the LSE and the block max-min estimator and the loss difference 
of the two estimators are listed in Table \ref{table}, 
along with the two-sided $p$-value for the difference. 
In Experiment I and IV which are less favorable to the block estimator, 
the block estimator still yields slightly smaller risk, although the risk difference is insignificant 
(with $p$-values $0.6190$ and $0.1600$ respectively)
In all other four experiments the block max-min estimator significantly outperforms the LSE with $p$-values 0.0062 or smaller, supporting our theoretical analysis. 
It is worthwhile to mention that, although the risk values are incomparable due to different dimension $d$, we observe more significant difference in the mean squared losses between the LSE and the block max-min estimator in $d=3$ than in $d=2$, in view of the $p$-values and box plots.   
This observation coincide with Theorem \ref{coro-fix-piecewise} and its comparison to the existing risk bounds for the LSE.

\begin{table}
\begin{center}
\begin{tabular}{ c  c  c  c  c  c  c  c  c  c}
  \hline
$(d=2)$ & \multicolumn{3}{c}{Experiment I} & \multicolumn{3}{c}{Experiment II} & \multicolumn{3}{c}{Experiment III} \\
\cline{2-10}
 & LSE & block & diff & LSE & block & diff &  LSE & block & diff\\
 \cline{2-10}
mean &  0.0822 & 0.0807 & 0.0016 & 0.1029 & 0.0918 & 0.0111 & 0.0713 & 0.0603 & 0.0110\\
s.d.& 0.0096 & 0.0095 & 0.0031 & 0.0156 & 0.0149 & 0.0041 & 0.0115 & 0.0109 & 0.0033\\
  p-value & & & 0.6190 & & & 0.0062 & & & 0.0007\\
  \hline
$(d=3)$  & \multicolumn{3}{c}{Experiment IV} & \multicolumn{3}{c}{Experiment V} & \multicolumn{3}{c}{Experiment VI} \\
  \cline{2-10}
   & LSE & block & diff & LSE & block & diff & LSE & block & diff \\
   \cline{2-10}
  mean &  0.1412 & 0.1353 & 0.0059 & 0.1316 & 0.1096 & 0.0220 & 0.0917 & 0.0746 & 0.0170\\
  s.d. & 0.0119 & 0.0117 & 0.0042 & 0.0178 & 0.0169 & 0.0059 & 0.0160 & 0.0147 & 0.0045\\
  p-value & & & 0.1600 & & & 0.0002 & & & 0.0002\\
    \hline
\end{tabular}
\end{center}
\caption{The mean and standard deviation (s.d.) of the mean squared losses 
for the LSE and the block max-min estimator (block), and the mean, s.d. and two-sided p-value for the loss differences (diff = loss of LSE - loss of block estimator).}
\label{table}
\end{table}

\begin{figure}
  \centering
  \subfigure[true $f$ (unknown)]{
    \label{simu-4:a}
    \includegraphics[width=0.48\textwidth]{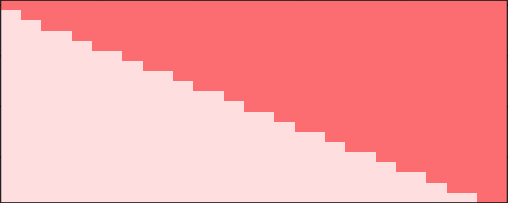}}
  \hspace{0.01\textwidth}
  \subfigure[$\by$ (observed)]{
    \label{simu-4:b}
    \includegraphics[width=0.48\textwidth]{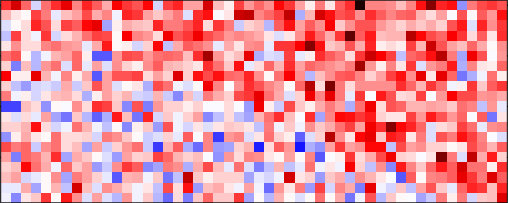}}
    \subfigure[the LSE]{
    \label{simu-4:c}
    \includegraphics[width=0.48\textwidth]{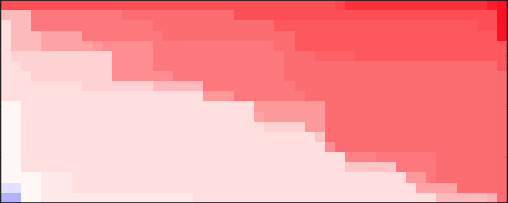}}
     \hspace{0.01\textwidth}
     \subfigure[the block max-min estimator]{
    \label{simu-4:d}
    \includegraphics[width=0.48\textwidth]{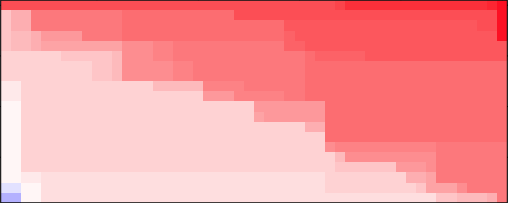}}

  \subfigure{
    \includegraphics[width=0.5\textwidth]{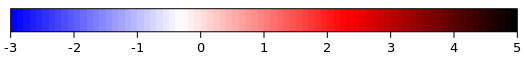}}

\caption{Heatmaps for the true two-piece function $\bff$, an observed $\by$, and its LSE and max-min estimate.}
\label{simu-4}
\end{figure}

We end this section with an example in which the LSE actually yields slightly smaller mean squared risk than the block max-min estimator. In Experiment VII, we consider the two-piece function $f(x_1, x_2) = I\{ x_1/n_1 + x_2/n_2 \ge 1\}$ 
on an $n_1 \times n_2$ lattice. 
Same as in Experiment I, II and III, we take $(n_1,n_2)=(50,20)$ and add standard gaussian noises to $f(x_1,x_2)$. See the heat maps in Figure \ref{simu-4}.

We shall recall $\fhat^{(lse)}_n(\bx) = \ybar_{U\cap L}$ for some upper set $U$ and lower set $L$. Suppose $x_1/n_1 + x_2/n_2 \ge 1$ so that $f(\bx)=1$, then the best level set $U\cap L$ for this design point is the upper red triangle in Figure \ref{simu-4:a}. In contrast, as $\fhat^{(block)}_n(\bx) = \ybar_{[\bu, \bv]}$ for some $\bu$ and $\bv$, the best possible block contains at most half design points of the upper triangle (when $\bu = (n_1/2,n_2/2)$ and $\bv = (n_1, n_2)$). Therefore, the variability of the block estimator at each design point may be larger than the LSE, resulting in a greater risk. Indeed, when we compare them on $5000$ replications of $\by$ 
as in Experiments I-VI, the mean squared losses for the LSE has mean $0.0420$ and standard deviation $0.0090$, while for the block max-min estimator the mean is $0.0440$ and the standard deviation is $0.0087$. 
However, the difference is not significant as the mean and standard deviation for the loss difference are $-0.0020$ and $0.0040$, and the two-sided p-value is $0.6163$. 

It would be difficult to characterize settings or general examples in which the LSE outperforms the block estimator. 
When we set $f(\bx) = 0.5 I\{x_1/n_1 + x_2/n_2 \ge 1\}$, 
the average normalized $\ell_2$ loss for the LSE is $0.0298$, slightly greater than $0.0280$ for the block max-min estimator, but the difference is insignificant as the two-sided $p$-value is $0.5568$. 


\newpage
\centerline{\large SUPPLEMENTARY MATERIAL}
\medskip
\centerline{\bf \large Supplement to ``Isotonic Regression in Multi-Dimensional Spaces and Graphs''}
\noindent
This supplement contains proofs of all the theoretical results stated in the main body of the paper.

\newpage
\bibliographystyle{apalike}
\bibliography{multiple-iso-aos}

\newpage
\centerline{\large SUPPLEMENTARY MATERIAL}
\medskip
\centerline{\bf \large Supplement to ``Isotonic Regression in Multi-Dimensional Spaces and Graphs''}
\bigskip

\noindent
The supplement contains proofs of all the theoretical results stated in the main body of the paper.

\bigskip
{\bf A1. Proofs of the results in Subsection 3.1}

\medskip
{\bf A1.1. Proof of Theorem \ref{th-1}.}

By the definition of $\bv_{\bx}$ in (\ref{k-m}),
\bes
\fhat^{(block)}_n(\bx)
\le \max_{\bu \preceq \bx}  \sum_{\bx_i\in [\bu,\bv_{\bx}]} \frac{y_i}{n_{\bu,\bv_{\bx}}}
\le f(\bx) + r_{q,+}^{1/q}(m_{\bx})
+ \max_{\bu \preceq \bx}  \sum_{\bx_i\in [\bu,\bv_{\bx}]} \frac{\veps_i}{n_{\bu,\bv_{\bx}}},
\ees
where $\bx = \bx_1, \ldots, \bx_n$.
Thus, by the definition of $r_{q,+}(m)$ in (\ref{r-plus}),
\bes
\E \Big\{\fhat^{(block)}_n(\bx) - f(\bx)\Big\}_+^q
&\le&  \E\bigg(r_{q,+}^{1/q}(m_{\bx})
+\max_{\bu \preceq \bx}  \sum_{\bx_i\in [\bu,\bv_{\bx}]} \frac{\veps_i}{n_{\bu,\bv_{\bx}}}\bigg)_+^q
\cr &\le&  2^{q-1}r_{q,+}(m_{\bx}) + 2^{q-1}\E\bigg(
\max_{\bu \preceq \bx}  \sum_{\bx_i\in [\bu,\bv_{\bx}]} \frac{\veps_i}{n_{\bu,\bv_{\bx}}}\bigg)_+^q
\cr &\le & 2^qr_{q,+}(m_{\bx}).
\ees
Similarly, we can have the second inequality in (\ref {th-1-1}). It follows that with the $\ell_+(m)$ in (\ref{ell-pm}),
\bes
T_{q,+}(V_0) \le \sum_{m=1}^{\infty}2^q r_{q,+}(m)\Big\{\ell_+(m)-\ell_+(m-1)\Big\} + \sum_{m=1}^{\infty}2^q r_{q,-}(m)\Big\{\ell_-(m)-\ell_-(m-1)\Big\}.
\ees
Hence (\ref {th-1-2}) follows as $r_{q,\pm}(m)$ is non-increasing. $\hfill\square$

\bigskip
{\bf A2. Proofs of the results in Subsection 3.2}

\medskip
{\bf A2.1. Proof of Lemma \ref{lm-lower-bd}.} 
For $\bm = (m_1,\ldots,m_d)^T\in \N_+^d$, define $m^*=\prod_{j=1}^d m_j$ and 
\bes
\scrK_{\bm}= \Big\{\bk\in \N_+^d:\, k_j\le n_j/m_j\,\forall\,j\le d, 1 \le k_1+\cdots+k_d - k^* +1 \le \sqrt{m^*}\Delta_n^*\Big\},
\ees
For a certain integer $k^*\ge d$ to be determined later 
and $N_{\bm} = |\scrK_{\bm}|$, we shall first prove that
\bel{pf-lm-lower-bd-1}&&
\inf_{\hbf}\, \sup\Big\{\E \big\|\hbf - \bff_n \big\|_q^q:
\bff_n\in\scrF_n, \Delta(\bff_n/\sigma)\le\Delta_n^*\Big\}
\ge c_q \sigma^q
\max_{\bm \in \scrM}\Big( (m^*)^{1-q/2} N_{\bm}\Big). 
\eel
Let $K_j = \lfloor n_j/m_j\rfloor$ and $\bn'=(K_1m_1,\ldots,K_dm_d)^T$.
The lattice $[{\bf 1},\bn']$, contained in $V=[{\bf 1},\bn]$,
is a lattice of $K_1\times\cdots\times K_d$ blocks of size
$m_1\times\cdots\times m_d$,
indexed by $\bk=(k_1,\ldots,k_d)^T$,
$1\le k_j\le K_j$. Suppose $f(\bx)$ is known to be piecewise constant and
non-decreasing
in this partition of blocks.
Let $g(\bk)$ be the value of $f$ on block $\bk$
and $\scrG_n$ be the class of $g(\cdot)$ satisfying
\bes
g(\bk) = \sigma \min\Big\{\Delta_n^*, (m^*)^{-1/2}
\big[\theta(\bk) + (k_1+\cdots+k_d - k^*)_+\big]
 \Big\},\ \theta(\bk) \in \{0,1\}.
\ees
As $k_1+\cdots+k_d - k^*$ is strictly increasing in $k_j$ for each $j$ with increment 1,
$g(\bk)$ is non-decreasing in $k_j$ for each $j$ 
and $0\le g(\bk)\le \sigma \Delta^*_n$. 
Set $f(\bx)= \sigma \Delta^*_n$ for $\bx\in [{\bf 1},\bn]\setminus [{\bf 1},\bn']$. 
As $f({\bf 1})=0$, $f(\bx)$ is non-decreasing in the entire $V$ and 
$\{f(\bn)-f({\bf 1})\}/\sigma = \Delta^*_n$. 
Note that $g(\bk) = \sigma (m^*)^{-1/2}\big[\theta(\bk) + (k_1+\cdots+k_d - k^*)_+\big]$ on $\scrK_{\bm}$
for all $g\in\scrG_n$.
Let $\ybar_{\bk}$ be the sample mean in the block indexed by $\bk$.
As $\ybar_{\bk}\sim N(g(\bk),\sigma^2/m^*)$ are sufficient for the estimation of $\bg$,
\bes
\inf_{\hbf} \sup_{\bff_n\in\scrF_n} \E \big\|\hbf - \bff_n \big\|_q^q
&\ge& \inf_{\hbg} \sup_{\bg \in\scrG_n} m^* \,\E \big\|\hbg - \bg \big\|_q^q
\cr &\ge& \inf_{\hbtheta} \sup_{\btheta \in \{0,1\}^{N_{\bm}}}
m^*(\sigma/\sqrt{m^*})^q \E \big\|\hbtheta - \btheta \big\|_q^q
\cr &\ge& m^*(\sigma/\sqrt{m^*})^q \E_{\pi} \big\|\hbtheta - \btheta \big\|_q^q 
\cr &=& m^*(\sigma/\sqrt{m^*})^q N_{\bm} \E_\pi \big|\hmu - \mu\big|_q^q,
\ees
where the infimum is taken over
$\htheta(\bk)$ based on $\ybar_{\bk}$, 
$\E_\pi$ is the joint expectation under which $\btheta$ 
has iid Bernoulli$(1/2)$ (prior) distribution, 
and $\hmu$ 
is the Bayes rule based on a single observation $X$ with $X|\mu \sim N(\mu,1)$ 
and $\mu\sim$Bernoulli$(1/2)$. 
This gives (\ref{pf-lm-lower-bd-1}).

Consider fixed $m_1,\ldots,m_d$ with $\sqrt{m^*}\Delta_n^* \ge 1$.
Assume without loss of generality $K_1\ge \cdots\ge K_d$.
For $K_1<2d$, we take $k^*=d$ so that $N_{\bm} \ge 1$.
For $K_1\ge 2d$, we take $k^*=\lfloor K_1/2\rfloor\ge d$, so that for all
$k^* \le k \le \min(k^*-1+\sqrt{m^*}\Delta_n^*,K_1+d-1)$
and $k_2+\cdots+k_d < k^*$, a solution $\bk$ exists satisfying
$k_1+\cdots+k_d = k$, $0\le k-k^* < k_1 \le k - (d-1) \le K_1$ and
\bes
k_1+\cdots+k_d-k^*+1= k-k^*+1 \le \sqrt{m^*}\Delta_n^*.
\ees
Such $\bk$ belongs to $\Lambda$ iff $k_j\le K_j$ for all $2\le j\le d$.
As $2(k^*+1)\ge K_1$ and $k^*\ge d$,
$\#\{(k_2,\ldots,k_d): 1\le k_j\le K_j\ \forall\, j\ge 2, k_2+\cdots+k_d < k^*\}\ge c_d \prod_{j=2}^dK_j$.
As the number of allowed $k$ is
$\min(\lfloor\sqrt{m^*}\Delta_n^*\rfloor,K_1+d-k^*)
\ge \min(\lfloor\sqrt{m^*}\Delta_n^*\rfloor, K_1/2+d)$, we find that
\bes
N_{\bm} \ge  \min(\sqrt{m^*}\Delta_n^*, K_1/2+d) c_d K_1^{-1}\prod_{j=1}^dK_j
\ge \frac{c_d\, n}{2^{d+1}m^*}\min\bigg( \frac{\sqrt{m^*}\Delta_n^*}{K_1}, 1\bigg).
\ees
This gives (\ref{lm-lower-bd-1}). As (\ref{lm-lower-bd-1}) is decreasing in $m^*$ given
$K_1=\max_{1\le j\le d}\lfloor n_j/m_j\rfloor$, its optimal configuration is attained
when either $K_j=K_1$ or $m_j=1$ for each $j$.
$\hfill\square$

\medskip
{\bf A2.2. Proof of Proposition \ref{prop-lower-bd}.}
For the optimal configuration of $\bm$ in (\ref{lm-lower-bd-1}),
there exist integers $s \in [1,d]$ and $K_1 = \lfloor n_s/m_s\rfloor \in [n_{s+1},n_s]$
such that the lower bound in (\ref{lm-lower-bd-1}) is maximized
with $\lfloor n_j/m_j\rfloor = K_1$ for $j\le s$ and $m_j = 1$ for $s<j\le d$. Thus,
$2^{-s}n_s^*/m^* < K_1^s \le n_s^*/m^*$ and
\bel{pf-prop-lower-bd-1}
&& \max_{\bm \in \scrM} \bigg\{ \frac{1}{(m^*)^{q/2}}
\min\bigg( \frac{\sqrt{m^*}\Delta_n^*}{\max_j \lfloor n_j/m_j\rfloor}, 1\bigg) \bigg\}
\\ \nonumber &\ge& c_d\max_{1\le s \le d}\,\max_{n_{s+1}\le (n_s^*/m^*)^{1/s} \le n_s}
\bigg\{\frac{1}{(m^*)^{q/2}}
\min\bigg(1, \frac{\sqrt{m^*}\Delta_n^*}{(n_s^*/m^*)^{1/s}}\bigg):
\sqrt{m^*}\Delta_n^*\ge 1\bigg\}
\\ \nonumber &=& c_d\max_{1\le t\le n}\Big\{ t^{-q/2}H(t):\ h_0(t) \ge 1\Big\}
\\ \nonumber &=& c_d \max_{1\le t \le n}
\Big\{\min\big(h_1(t), h_2(t)\big):  h_0(t)\ge 1\Big\},
\eel
where $h_0(t) = \Delta_n^*t^{1/2}$ and
$H(t) = \min\big\{1, h_0(t)/(n^*_s/t)^{1/s}\big\}, t\in [t_s,t_{s+1}]$ with
$t_s = n^*_s/n_s^s$ are as stated, $h_1(t) = t^{-q/2}$, and
\bes
h_2(t) = \frac{h_1(t)h_0(t)}{(n_s^*/t)^{1/s}}
= \frac{\Delta_n^* t^{1/s - (q-1)/2}}{(n^*_s)^{1/s}},\ t\in [t_s,t_{s+1}],\ s = 1,\ldots, d.
\ees
We note that $t_s \uparrow s$, $t_1=1$ and $t_{d+1}=n$.
As $(n_s^*/t_s)^{1/s} = n_s$ and $(n_s^*/t_{s+1})^{1/s} = n_{s+1}$, $H(t)$ and $h_2(t)$ are
continuos in $t$.
As $s_q = \lceil 2/(q-1)\rceil \wedge (d+1)$, $1/s-(q-1)/2\le 0$ iff $s\ge s_q$
for $1\le s\le d$.
It follows that $h_2(t)$ is increasing in $t$ for $t\le t_{s_q}$ and non-increasing in $t$ for $t\ge t_{s_q}$.

Thus, the optimal solution is given by
\bel{pf-prop-lower-bd-2}
&& \max_{h_0(t)\ge 1}\min\{h_1(t),h_2(t)\}
\cr &=& \begin{cases} h_1(1), & n_1 \le h_0(1),
\cr  h_1(t_*), & h_0(t_*) = (n^*_{s}/t_*)^{1/s} \hbox{ with }t_s\le t_* \le t_{s+1}\le t_{s_q}, 
\cr  h_2(t_{s_q}), & 1\le h_0(t_{s_q}) \le (n^*_{s_q}/t_{s_q})^{1/s_q} = n_{s_q},
\cr  h_2(t^*), & h_0(t^*) = 1,\ t_{s_q} \le t^*\le n.
\end{cases}
\eel

We note that $h_2(t)\ge h_1(t)$ for all $1\le t\le n$ in the first case above,
$h_1(t_*)=h_2(t_*)$ and $h_0(t_*)\ge 1$ in the second case,
$h_2(t_{s_q}) \le h_1(t_{s_q})$ in the third case, and  
the maximizer $t = t^*$ is determined by the constraint $h_0(t)\ge 1$ 
in the fourth case.  
Moreover, for $1\le t_* \le t_{s_q}$ we have 
\bes
t_* = \{(n^*_s)^{1/s}/\Delta^*_n\}^{2s/(2+s)}\in [t_s,t_{s+1}]
\ \hbox{ and }\ n_{s+1}/t_{s+1}^{1/2}
\le \Delta^*_n \le n_{s}/t_{s}^{1/2},
\quad s < s_q,
\ees
and for $t_{s_q} \le t^*\le n$ we have
\bes
t^* = (\Delta_n^*)^{-2}
\ \hbox{ and }\ t_{s+1}^{-1/2}
\le \Delta^*_n \le t_{s}^{-1/2},
\quad s_q\le s\le d.
\ees
Thus, the right-hand side of (\ref{pf-prop-lower-bd-2}) is
\bes
\begin{cases} 1, & n_1 \le \Delta^*_n, \qquad\qquad\qquad\qquad\ \ \, (s=0), 
\cr  \displaystyle
\big(\Delta^*_n/(n_s^*)^{1/s}\big)^{qs/(2+s)}
&\displaystyle n_{s+1}/t_{s+1}^{1/2}
\le \Delta^*_n \le n_{s}/t_{s}^{1/2},\quad
(1\le s < s_q),
\cr \displaystyle
\Delta_n^*/\big(n_s t_{s}^{(q-1)/2}\big),
&\displaystyle t_{s}^{-1/2}
\le \Delta^*_n \le n_{s}/t_{s}^{1/2}, \qquad\quad\, (s=s_q\le d),
\cr \displaystyle
(\Delta_n^*)^{q-2/s}/(n^*_s)^{1/s},
& \displaystyle t_{s+1}^{-1/2}
\le \Delta^*_n \le t_{s}^{-1/2},\qquad\qquad\, (s_q\le s\le d).
\end{cases}
\ees

This gives (\ref{prop-lower-bd-1}) for $\Delta^*_n\ge n^{-1/2}$
through (\ref{lm-lower-bd-1}), (\ref{pf-prop-lower-bd-1}) and (\ref{pf-prop-lower-bd-2}).
As $n^*_d=n$ and $n_{d+1}=1$, the rate is $n^{-q/2}$ for
$\Delta^*_n=n^{-1/2}$. As the minimax rate for $\Delta^*_n=0$ is the same
$n^{-q/2}$ due to the unknown average of $\bff_n$, (\ref{prop-lower-bd-1})
also holds for $0\le \Delta^*_n \le n^{-1/2}$.
Finally, for (\ref{prop-lower-bd-2}), we note that $s_q\ge d$ iff $q\le 1+2/d$.
$\hfill\square$

\bigskip
{\bf A3. Proofs of the results in Subsection 3.3}

\medskip
{\bf A3.1. Proof of Lemma \ref{lm-doob}}.
As $\scrF_t^{(j)}$ does not depend on $(s_1,\ldots,s_{j-1}, t_{j+1},\ldots,t_d)$,
\bes
\E\Big[ \max_{s_1\le t_1,\ldots,s_{j-1}\le t_{j-1}}\big|f_{s_1,\ldots,s_{j-1}, t, t_{j+1},\ldots,t_d}\big|
\Big| \scrF_s^{(j)}\Big]
\ge \max_{s_1\le t_1,\ldots,s_{j-1}\le t_{j-1}}\big|f_{s_1,\ldots,s_{j-1}, s, t_{j+1},\ldots,t_d}\big|
\ees
for all $t>s$. Thus, repeated application of the Doob inequality gives
\bes
\E \max_{s_1\le t_1,\ldots,s_{j}\le t_{j}}\big|f_{s_1,\ldots,s_{j}, t_{j+1},\ldots,t_d}\big|^q
\le \bigg(\frac{q}{q-1}\bigg)^q\E \max_{s_1\le t_1,\ldots,s_{j-1}\le t_{j-1}}
\big|f_{s_1,\ldots,s_{j-1}, t_j, t_{j+1},\ldots,t_d}\big|^q.
\ees
The conclusion for the general $f_{\bt}$ follows. 

For independent $\veps_i$, define $f_{\bt} = \sum_{\bx_i\le \bt} \veps_i$ and
$\scrF_t^{(j)}=\sigma\{\veps_i: x_{i,j}\le t\}$ where $x_{i,j}$
is the $j$-th component of $\bx_i$.
As $\E \veps_i=0$, $\{f_{s_1,\ldots,s_{j-1}, t, t_{j+1},\ldots,t_d}, t\in\R\}$ is a sub-martingale
with respect to the filtration $\{\scrF_t^{(j)}, t\in\R\}$.
We apply H\"older's inequality to avoid the singularity at $q=1+$.
$\hfill\square$

\medskip
{\bf A3.2. Proof of Proposition \ref{prop-single-point}.} 
As \eqref{new-prop-2-3} is a consequence of \eqref{new-prop-2-1} and \eqref{new-prop-2-2}, 
we only need to prove \eqref{new-prop-2-1} by symmetry. 
Moreover, by Theorem \ref{th-1} it suffices to prove that for any $\bv_0\in [\bx_i,\bn]$ 
\bes
r_{q,+}(m_{\bx_i}) \lesssim_{q,d} \E\Big(\ybar_{[\bx_i,\bv_0]}-f(\bx_i)\Big)_+^q  
\ees
Let $\bv_1$ be a design point in $[\bx_i,\bv_0]$ satisfying 
$\min(n_{\bx_i,\bv_1},n_{\bv_1,\bv_0})/ n_{\bx_i,\bv_0}\ge 1/2^d $. 
Such $\bv_1$ always exists as the minimum is attained when $ \bv_0=\bx_i+  {\bf 1} $. 
As $\ybar_{[\bx_i,\bv_0]} - f(\bx_i)$ is Gaussian, 
\bes
\E\Big(\ybar_{[\bx_i,\bv_0]}-f(\bx_i)\Big)_+^q 
&\ge& \frac{1}{2}\bigg\{\Big(\fbar_{[\bx_i,\bv_0]}-f(\bx_i)\Big)^q 
+ \frac{\sigma^q \E(N(0,1)_+^q)}{n_{\bx_i,\bv_0}^{q/2}}\bigg\}
\cr &\ge& 2^{-1} 2^{-qd} \min\big\{1,\E(N(0,1)_+^q)\big\}
\Big\{\big(f(\bv_1)-f(\bx_i)\big)^q + \sigma^q n_{\bx_i,\bv_1}^{-q/2}\Big\}. 
\ees
If (\ref{r-plus}) holds with $r_{q,+}(m) = C_{q,d}\sigma^q m^{-q/2}\ge \sigma^q m^{-q/2}$, 
by the definition of $m_{\bx_i}$ 
\bes
r_{q,+}(m_{\bx_i})
&\le & \min_{\bx\le \bv\le \bv_1}
\Big\{ C_{q,d}\sigma^q n_{\bx,\bv}^{-q/2}: C_{q,d}^{1/q}\sigma n_{\bx,\bv}^{-1/2}
\ge f(\bv)-f(\bx_i)\Big\}
\cr &\le&\begin{cases} 
C_{q,d}\sigma^q n_{\bx,\bv_1}^{-q/2}, & f(\bv_1)-f(\bx_i) \le C_{q,d}^{1/q}\sigma n_{\bx_i,\bv_1}^{-1/2}, \cr
f(\bv_1)-f(\bx_i), & f(\bv_1)-f(\bx_i) > C_{q,d}^{1/q}\sigma n_{\bx_i,\bv_1}^{-1/2}
\end{cases} 
\cr &\le& \Big(2C_{q,d} 2^{qd}\big/\min\big\{1,\E(N(0,1)_+^q)\big\}\Big)
\E\Big(\ybar_{[\bx_i,\bv_0]}-f(\bx_i)\Big)_+^q. 
\ees
Thus, the proof is complete if we verify (\ref{r-plus}) with the 
$r_{q,+}(m)$ 
in \eqref{r-bd}. 
The rest of the proof is devoted to this task. 
In fact, we prove 
\bel{pf2-1-2}
\E \Big( \max_{ \bu \le \bx} \sum_{\bx_i \in [\bu,\bv]} \frac{\veps_i}{n_{\bu,\bv}} \Big)^q_{+}
\le C_{q,d} \sigma^qm^{-q/2} = r_{q,+}(m)
\eel
under the moment condition $\E \veps_i=0$ and $\E|\veps_i|^{q\vee 2}\le\sigma^{q\vee 2}$, 
without assuming the normality of $\veps_i$. 

To control the maximization over $\bu \in [{\bf 1},\bx]$ in (\ref{r-plus}),
we cover the set by a collection $U_{\bx,\bv}$ of blocks
$\big\{[\ubu_{\bi}, \obu_{\bi}]\big\}$ indexed by vector $\bi\in \N_+^d$ as follows. Define
\bes
U'_{\bx,\bv} &=& \Big\{[\ubu'_{\bi},\obu'_{\bi}]: \frac{v_j -(\uu_{\bi})_j}{v_j - x_j} = 2^{i_j}, \frac{(\ou_{\bi})_j -(\uu_{\bi})_j}{v_j - x_j} = 2^{i_j-1}, \bi=(i_1,\ldots, i_d), i_j \in \N^+ \Big\},
\ees
which covers $\{\bu: \bu\le \bx\}$. The covering collection for $[{\bf 1},\bx]$ is defined as
\bel{Uxv}
U_{\bx,\bv} &=& \Big\{[\ubu_{\bi},\obu_{\bi}] \neq \emptyset: [\ubu_{\bi},\obu_{\bi}]
=[\ubu'_{\bi},\obu'_{\bi}] \cap [{\bf 1},\bx], [\ubu'_{\bi},\obu'_{\bi}] \in U'_{\bx,\bv}\Big\}.
\eel
Here the edge of the blocks are allowed to overlap to simplify the discussion.
Observe
\bes
&& \E \Big( \max_{ \bu \le \bx} \sum_{\bx_i \in [\bu,\bv]} \frac{\veps_i}{n_{\bu,\bv}} \Big)^q_{+}
\cr &\le& \bigg(\E \max_{ \bu \le \bx}
\bigg| \sum_{\bx_i \in [\bu,\bv]} \frac{\veps_i}{n_{\bu,\bv}} \bigg|^{q\vee 2}\bigg)^{q/(q\vee 2)}
\cr &\le&
\bigg(\sum_{[\ubu_{\bi},\obu_{\bi}]\in U_{\bx,\bv}} \E \max_{\bu \in [\ubu_{\bi},\obu_{\bi}] }
\bigg| \sum_{\bx_i \in [\bu,\bv]} \frac{\veps_i}{n_{\obu_{\bi},\bv}}\bigg|^{q\vee 2}\bigg)^{q/(q\vee 2)}
\ees
due to $n_{\bu,\bv}\ge n_{\obu_{\bi},\bv}$ for $\bu \in [\ubu_{\bi},\obu_{\bi}]$.
It follows from the definition of $U_{\bx,\bv}$ that $n_{\obu_{\bi},\bv} = 2^{|\bi| - d}m$
and $n_{\ubu_{\bi},\bv} \le 2^{|\bi| }m$,
where $m=n_{\bx,\bv}$ and $|\bi| = \sum_{j=1}^d i_j$. Therefore, by Lemma \ref{lm-doob},
\bes
&& \sum_{[\ubu_{\bi},\obu_{\bi}]\in U_{\bx,\bv}}
\E \max_{\bu \in [\ubu_{\bi},\obu_{\bi}] }\Big( \frac{1}{2^{|\bi|-d}m}
\Big|\sum_{\bx_i \in [\bu,\bv]} {\veps_i} \Big|\Big)^{q\vee 2}
\\
\nonumber &\le& 2^{(q\vee 2)d} \sum_{[\ubu_{\bi},\obu_{\bi}]\in U_{\bx,\bv}}
\frac{1}{(2^{|\bi|}m)^{q\vee 2}}
\Big(\frac{{q\vee 2}}{{q\vee 2}-1}\Big)^{(q\vee 2)d}
\E  \Big|\sum_{\bx_i \in [\ubu_{\bi},\bv]} {\veps_i} \Big|^{q\vee 2}
\cr
&\le& C^d_q \sigma^{q\vee 2}
m^{-(q\vee 2)/2} \sum_{[\ubu_{\bi},\obu_{\bi}]\in U_{\bx,\bv}}
\frac{1}{2^{(q\vee 2)|\bi|/2}}
\E  \Big|\frac{1}{\sqrt{n_{\ubu_i, \bv}}}\sum_{\bx_i \in [\ubu_{\bi},\bv]}
\frac{\veps_i}{\sigma} \Big|^{q \vee 2},
\ees
where $C_q = \{2(q\vee 2)/(q\vee 2-1)\}^{q\vee 2}$. The Rosenthal inequality gives
\bes
\E \Big|\frac{1}{\sqrt{n_{\ubu_{\bi},\bv}}} \sum_{\bx_i \in [\ubu_{\bi},\bv]} \frac{\veps_i}{\sigma} \Big|^{q \vee 2}
\le C'_{q \vee 2} \max \Big\{ \sum_{\bx_i \in [\ubu_{\bi}, \bv]}
\E \Big| \frac{\veps_i}{\sqrt{n_{\ubu_{\bi}, \bv}}\sigma}\Big|^{q \vee 2}, 1 \Big\}
\le C'_{q \vee 2},
\ees
where constant $C_{q \vee 2}'$ is continuous in $q \ge 1$.
It follows that
\bes
&& \E \Big( \max_{ \bu \le \bx} \sum_{\bx_i \in [\bu,\bv]} \frac{\veps_i}{n_{\bu,\bv}} \Big)^q_{+}
\\
\nonumber &\le& \big(C^d_q C_{q \vee 2}'\big)^{q/(q\vee 2)} \sigma^qm^{-q/2}
\sum_{[\ubu_{\bi},\obu_{\bi}]\in U_{\bx,\bv}} \frac{1}{2^{q|\bi|/2}}
\cr
&\le&\big(C^d_q C_{q \vee 2}'\big)^{q/(q\vee 2)} \sigma^qm^{-q/2}
\sum_{i_1=1 }^{\infty} \cdots \sum_{i_d=1 }^{\infty}  \frac{1}{2^{q(i_1+\cdots + i_d) /2}}
\cr
&\le& C_{q,d} \sigma^qm^{-q/2},
\ees
where $C_{q,d} = \big(C^d_q C_{q \vee 2}'\big)^{q/(q\vee 2)}(2^{q/2}-1)^d$
remains continuous in $q \ge 1$. 
This gives \eqref{pf2-1-2}.  
$\hfill\square$ 

\begin{remark}
If $\veps_i$'s are i.i.d., Lemma \ref{lm-doob} gives
specifically $C_{q,d} = \{(q\vee 2)/((q\vee 2)-1)\}^{qd}$ in \eqref{pf2-1-2} as the process below is a multi-indexed reverse martingale,  
\bes
W(\bu) = \sum_{\bx_i \in [\bu,\bv]} \frac{\veps_i}{n_{\bu,\bv}}, \quad \bu \le \bx. 
\ees
\end{remark}

\medskip
{\bf A3.3. Proof of Theorem \ref{th-lattice-general}.}
Let $V_0 = [\ba, \bb]$.
It follows from Theorem \ref{th-1} that
\bel{Rqab}&& 
T_q(V_0) = \sum_{\bx_i \in [\ba,\bb]} \E\big|\fhat_n(\bx_i)-f(\bx_i)\big|^q
\le 2^q\sum_{\bx_i \in [\ba,\bb]} r_{q,+}(m_{\bx_i}) + 2^q\sum_{\bx_i \in [\ba,\bb]}r_{q,-}(m_{\bx_i,-}).
\eel
As the analysis of $r_{q,-}(m_{\bx_i,-})$ is symmetric to that of $r_{q,+}(m_{\bx_i})$,
it suffices to consider the first term on the right-hand side above. To this end, 
we divide $V_0$ into two parts: $V_{0,+}$ and $V_0 \setminus V_{0,+}$, where 
$V_{0,+}$ is composed of design points $\bx_i$ in $V_0$ for which 
the positive-side $L_q$ risk measure $r_{q,+}(m_{\bx_i})$ 
can be controlled by the variation of $f$ within 
$V_0$. More precisely,
\bes
V_{0,+} = \Big\{\bx_i \in[\ba,\bb]: f(\bb) > f(\bx_i) + r_{q,+}^{1/q}(n_{\bx_i,\bb})\Big\}. 
\ees
According to the definition of $m_{\bx_i}$ in (\ref {k-m}), 
$m_{\bx_i}\ge n_{\bx_i,\bb}$ for design points $\bx_i$ in $V_0 \setminus V_{0,+}$, 
so that by the monotonicity of the variability measure
$r_{q,+}(m_{\bx_i}) \le r_{q,+}(n_{\bx_i,\bb})$, representing the edge effect 
when $f$ is flat within $V_0$ but the data for the design points beyond $V_0$ are not used. 
In our analysis, 
\bel{Tq+1}
T_{q,+}(V_{0,+}) = 2^q \sum_{\bx_i \in V_{0,+}} r_{q,+}(m_{\bx_i})
\eel
is controlled by the range-to-noise ratio $\Delta_n^*$ and the dimensions $\ntil_j$ of the block $V_0=[\ba,\bb]$, and 
\bel{Tq+2}
T_{q,+}(V_0 \setminus V_{0,+}) = 2^q \sum_{\bx_i \notin V_{0,+}} r_{q,+}(m_{\bx_i}) 
\le 2^q \sum_{\bx_i \in [\ba,\bb]} r_{q,+}(n_{\bx_i,\bb})
\eel
is relatively easy to bound 
as $r_{q,+}(m)$ is explicitly given in \eqref{pf2-1-2}. 

%
Let $\ell^*_{0,+}(m)$ be any upper bound for $\ell_{0,+}(m)=\#\{\bx\in V_{0,+}: m_{\bx} \le m\}$.
We have
\bes
T_{q,+}(V_{0,+})
&=& \sum_{m=1}^{n_{\ba, \bb}} 2^qC_{q,d}\sigma^q m^{-q/2}
\Big\{\ell_{0,+}(m)-\ell_{0,+}(m-1)\Big\}
\cr &\le & 2^qC_{q,d}\sigma^q \sum_{m=1}^{\ntil_d^*} m^{-q/2}
\Big\{\ell^*_{0,+}(m)-\ell^*_{0,+}(m-1)\Big\}.
\ees

In the following four steps, we respectively bound $\ell_{0,+}(m)$ 
when $d$ is the effective dimension, bound $T_{q,+}(V_{0,+})$,
bound $T_{q,+}(V_{0} \setminus V_{0,+})$, and draw the final conclusion on the sum $T_{q,+}(V_0) = T_{q,+}(V_{0,+} ) + T_{q,+}(V_0 \setminus V_{0,+}) $.

{\bf Step 1}:   
Bound $\ell_{0,+}(m)$ for $m\ge t_d = \ntil/\ntil_d^d$.  We 
first partition $V_0 = [\ba,\bb]$ 
into a lattice of $\ntil_d^d$ hyper-rectangular blocks, which we call unit blocks and write as  
\bes
B_{\bk} = \Big\{\bx_i: a_j + (k_j-1)\ntil_j/\ntil_d \le x_{i,j} < a_j +  k_j\ntil_j/\ntil_d, 1\le j\le d\Big\}, 
\ees
indexed by $\bk=(k_1, \ldots, k_d)^T \in \{1,\ldots,\ntil_d\}^d$, 
where $x_{i,j}$ is the $j$-th element of $\bx_i$. 
We note that while $\ba$ and $\bb$ are integer-valued vectors, 
$a_j +  k_j\ntil_j/\ntil_d$ are not necessarily integers unless 
$\ntil_j$ are all multipliers of $\ntil_d$, and in this special case 
the dimensions of $B_{\bk}$ are exactly proportional to that of $[\ba,\bb]$.  
In general $t_d$ is the average number of data points in $B_{\bk}$. 
In the simplest case where $\ntil_1=\cdots=\ntil_d$, $t_d=1$ and 
$B_{\bk}$ contains just a single design point 
$\ba + \bk - {\bf1}$. 

We then partition $[\ba, \bb]$ into diagonal sequences of unit blocks and denote by $L_{\bk}$ the sequence starting from the unit block $B_{\bk}$ on the lower half boundary of $[\ba, \bb]$.
Formally, we define $L_{\bk}$ as
\bes
L_{\bk} = \bigcup_{i=0}^{\min_j (\ntil_d - k_j)} B_{\bk +i{\bf 1} }, \hbox{ where }
\bk \in \{\bk: \exists\, j \le d\hbox{ s.t. } k_j=1\}.
\ees
Note that $\#\{\bk: \exists\, j \le d\hbox{ s.t. } k_j=1\} \le d\ntil_d^{d-1}$, 
so that there are at most $d \ntil_d^{d-1}$ such $L_{\bk}$'s. 
In the case of $d=2$ and $\ntil_1=\ntil_2$, each $L_{\bk}$ is the intersection 
of $V_0 = [\ba,\bb]$ and a 45 degree line in $\R^2$. 

Finally we partition $V_{0,+}$ according to the value of the function $f$ into
\bes
D_j = \Big\{\bx \in V_{0,+}: f(\ba) + (j-1) r_{q,+}^{1/q}\big((3k)^d t_d\big) \le  f(\bx)
< f(\ba) + j r_{q,+}^{1/q}\big((3k)^d t_d\big)\Big\},
\ees
$j=1, \ldots, J$, where $r_{q,+}(m)$ is as in \eqref{pf2-1-2}, 
$J=\big\lceil \{f(\bb)- f(\ba)\}/r_{q,+}^{1/q}\big((3k)^d t_d \big)\big\rceil$
and $k$ is a positive integer satisfying $(k-1)^dt_d < m \le k^dt_d$.
Figure~\ref{fig:proof} depicts a segment of $L_{\bk}$ passing through $D_j$ 
with a diagonal line of unit blocks in red color, in the case of $d=2$.

On $D_j \cap L_{\bk}$ consider design points 
$\bx \in D_j \cap B_{\bk + c{\bf 1}}$ and
$\bv \in D_j \cap B_{\bk + (c + 1+k ){\bf 1}}$  
for some integers $c\ge 0$ and $k\ge 1$. 
The lower and upper bounds for the individual coordinates of the design points in the two unit blocks provide 
\bes
k^dt_d < n_{\bx, \bv}
\le \big(k +2\big)^d t_d \le (3k)^d t_d. 
\ees
As $\bx$ and $\bv$ are both in $D_j$, 
$f(\bv)-f(\bx) \le r_{q,+}^{1/q}\big((3k)^d t_d\big)
\le r_{q,+}^{1/q}\big(n_{\bx, \bv}\big)$,
so that $m_{\bx} \ge n_{\bx, \bv} > k^dt_d$
by the definition of $m_{\bx}$ in (\ref{k-m}).
Thus, for $\bx\in D_j\cap L_{\bk}$, $m_{\bx}\le k^dt_d$ implies that
$\bx$ is within $k $ blocks away from the upper contour of $D_j$,
\bes
\#\big\{\bx \in V_{0,+} \cap L_{\bk} \cap D_j: m_{\bx} \le k^dt_d\big\}
\le (k +1)t_d.
\ees

The above bound holds for all $D_j$, but actually we
can replace the upper bound by 0 for $j=J$.
Let $\bx \in V_{0,+} \cap D_J$. We have
$f(\bb) \le f(\bx) + r_{q,+}^{1/q}\big((3k)^d t_d\big)$
by the definition of $D_J$ and $f(\bb) > f(\bx) + r_{q,+}^{1/q}\big(n_{\bx,\bb}\big)$
by the definition of $V_{0,+}$, so that $(3k)^d t_d < n_{\bx, \bb}$.
Consequently, there must exist two adjacent design points $\bv_1$ and $\bv_2=\bv_1+\bfe_j$
in $[\bx,\bb]$ for some canonic unit vector $\bfe_j$ such that
$n_{\bx,\bv_1} \le (3k)^d t_d
< n_{\bx,\bv_2} \le 2n_{\bx,\bv_1}$. 
It follows that
$f(\bv_1) \le f(\bb) \le f(\bx) + r_{q,+}^{1/q}\big((3k)^d t_d\big) \le f(\bx) + r_{q,+}^{1/q}(n_{\bx, \bv_1})$, so that $m_{\bx} \ge n_{\bx, \bv_1} \ge
n_{\bx,\bv_2}/2 > k^d t_d$.
Thus, 
\bes
\#\{\bx\in V_{0,+} \cap  D_J: \ m_{\bx}\le k^d t_d\}=0.
\ees

As $r_{q,+}(m) = C_{q,d}\sigma^q m^{-q/2}$ and $k = \lceil (m/t_d)^{1/d} \rceil$, 
we have
\bes
\ell_{0,+}(m) &\le & \# \{\bx \in V_{0,+}: m_{\bx} \le k^d t_d\}
\cr &=& \sum_{\bk} \sum_{j=1}^{J-1}
\#\big\{\bx \in V_{0,+} \cap L_{\bk} \cap D_j: m_{\bx}\le k^d t_d\}
\cr &\le& d \ntil_d^{d-1} (k +1 )t_d(J-1)
\cr &\le& d \ntil_d^{d-1} 2k t_d \frac{f(\bb)- f(\ba)}{r_{q,+}^{1/q}\big((3k)^d t_d \big)}
\cr &\le& d \ntil_d^{d-1} 4m^{1/d}t_d^{1-1/d} C_{q,d}^{-1/q}
\sqrt{(3k)^d t_d}\{f(\bb)- f(\ba)\}/\sigma
\cr &\le& d \ntil_d^{d-1}4 m^{1/d}t_d^{1-1/d}C_{q,d}^{-1/q}
3^{d/2} \sqrt{2}\sqrt{m}\Delta^*_n
\cr &=& C_{q,d}' \ntil_d^{d-1}t_d^{1-1/d} m^{1/d+1/2}\Delta^*_n
\cr &=& C_{q,d}' \Delta_n^* m^{1/d+1/2}\ntil^{1-1/d}
\ees
due to $\ntil_d^dt_d = \ntil$. As the above inequality holds for all integers $k\ge 1$, 
\bes
\ell_{0,+}(m) \le C_{q,d}'\Delta_n^* m^{1/d+1/2}\ntil^{1-1/d}\quad \forall\, m \ge t_d. 
\ees

{\bf Step 2}: Bound $T_{q,+}(V_{0,+})$.  
If $[\ba, \bb]$ is a hyper-cube, i.e., $\ntil_1 = \cdots = \ntil_d$, then $t_d=1$; otherwise, 
we still need to bound $\ell_{0,+}(m)$ for $1\le m < t_d = \ntil/\ntil^d$. 
Recall that $\ntil^*_s=\prod_{j=1}^s \ntil_j$ and $t_s=\ntil^*_s/\ntil_s^s$. 
As $1\le t_1\le \cdots\le t_d$, we just consider $t_s\le m\le t_{s+1}$ for some $1\le s<d$. 
We partition $V_0=[\ba,\bb]$ into $\ntil_{s+1}\cdots \ntil_d = \ntil/\ntil_s^*$ lattice slices 
of dimension $\ntil_1\times \cdots\times \ntil_s$ as follows,
\bes
V_{0,k_{s+1},\ldots,k_d} = \big\{\bx_i\in [\ba,\bb]: x_{i,j}=a_j+(k_j-1), s<j\le d\big\},\ 
k_j = 1,\ldots, \ntil_j, s<j\le d. 
\ees
We apply Step 1 
to each $V_{0,k_{s+1},\ldots,k_d}$ 
so that for all $t_s\le m\le t_{s+1}$ 
\bes
\#\Big\{\bx_i\in V_{0,+}\cap V_{0,k_{s+1},\ldots,k_d}: m_{\bx_i}\le m\Big\} 
\le C_{q,s}' \Delta_n^* m^{1/s+1/2}|V_{0,k_{s+1},\ldots,k_d}|^{1-1/s}
\ees
As $V_{0,k_{s+1},\ldots,k_d}$ is of $s$-dimensional and contains $\ntil_s^*$ data points, 
Step 1 
yields 
\bes
\ell_{0,+}(m) \le \min\Big\{\ntil,
(\ntil/\ntil_s^*) \times C_{q,s}' \Delta_n^* m^{1/s+1/2} (\ntil_s^*)^{1-1/s}\Big\}.
\ees
As $\Htil(t) = \min\big\{1, \Delta^*_n t^{1/2}(t/\ntil^*_s)^{1/s}\big\}$ for
$t_s \le t\le t_{s+1}$, it follows that
\bes
\ell_{0,+} (m) \le \ell_{0,+}^*(m) = C_{q,d}'\, \ntil \Htil(m),\ t_s\le m\le t_{s+1},\ 1\le s \le d.
\ees
Hence, as $\ntil = n_{\ba,\bb}$,
\bel{Tq+V0+}
T_{q,+}(V_{0,+})
&\le & C_{q,d}2^q\sigma^q \sum_{m=1}^{n_{\ba, \bb}} m^{-q/2}
\Big\{\ell^*_{0,+}(m)-\ell^*_{0,+}(m-1)\Big\}
\\
\nonumber  &\le& C_{q,d}2^qC_{q,d}'\sigma^q n_{\ba, \bb}
\bigg(\Htil(1) + \sum_{s=1}^d \int_{t_s}^{t_{s+1}}t^{-q/2}\Htil(dt)\bigg).
\eel

\medskip
{\bf Step 3:} We bound $T_{q,+}(V_0 \setminus V_{0,+})$ by
\bes
T_{q,+}(V_0 \setminus V_{0,+} )
&\le& 2^q\sum_{\bx_i\in [\ba,\bb]} r_{q,+}(n_{\bx_i,\bb})
\cr
&\le& 2^q C_{q,d} \sigma^{q} \sum_{\bx_i\in [\ba,\bb]}\prod_{j=1}^d \frac{1}{(b_j-x_{i,j}+1)^{q/2}}
\cr
&\le& 2^q C_{q,d}\prod_{j=1}^d\sum_{m=1}^{b_j-a_j+1} m^{-q/2}.
\ees

\medskip
{\bf Step 4:}
In view of (\ref {Rqab}), (\ref{Tq+1}) and (\ref{Tq+2}), the main
conclusion (\ref {th2-main}) directly follows from Steps 2 and 3 by summing \eqref{Tq+V0+} and the above inequality. 
Note that we may take in (\ref {th2-main})  
$C^*_{q,d} = \max\{ 2 C_{q,d} 2^q C_{q,d}', 2^{1+q}C_{q,d}\}$, 
which remains continuous in $q\in [1,\infty)$ and decreasing in $d$.
Note that $t_1=1$, $t_s=\ntil^*_s/\ntil_s^s\uparrow$ in $s$, $t_{d+1} = \ntil$, 
and $\Htil(t) = \min\big\{1, \Delta^*_n t^{1/2}(t/\ntil^*_s)^{1/s}\big\}$ for $t\in [t_s,t_{s+1}]$
is defined in the same way as in Proposition \ref{prop-lower-bd}.

To obtain more explicit bounds, we write
\bel{H-integral}
\Pi = \Htil(1) + \int_1^{n_{\ba,\bb}}t^{-q/2}\Htil(dt)
= \Htil(1) + \int_1^{t_*}t^{-q/2}\Htil(dt),
\eel
where $t_*=\min\{ t \ge 1: \Htil(t)=1 \hbox{ or } t =n_{\ba,\bb}\}$. 
Note that $n_{\ba,\bb}=\ntil$ and $t_*\in [1,\ntil]$ exists as $H(t)$ is strictly increasing 
and continuous in $t\in [1,\ntil]$.  

For $\Delta^*_n\ge \ntil_1=\ntil_1^*$, we have $t_*=\Pi=1$, which gives the first case of (\ref{th2-1}).

For $\Delta^*_n \le \ntil_1$, 
$t_* =((\ntil_{s}^*)^{1/s}/\Delta^*_n)^{2s/(2+s)}\wedge \ntil$ 
for some $s$ satisfying $t_*\in [t_{s},t_{s+1}]$. 
For $t_*<t_{s_q}$, this matches the $t_*$ in (\ref{pf-prop-lower-bd-2}) 
where the lattice is of size $n_1 \times \cdots \times n_d$.  
Inside the interval $[t_{s},t_{s+1}]$, 
$t^{-q/2}\Htil(dt) = (1/2 + 1/s)\big( \Delta^*_n/(\ntil^*_s)^{1/s}\big)t^{1/s-1/2-q/2}dt$. 
If $1/s-1/2-q/2 \neq -1$ for all $s=1,\ldots, d$, 
$\int_1^{t_*}t^{-q/2}\Htil(dt) = \max_{1\le t\le t^*}\big\{t^{-q/2}\Htil(t)\big\}$ as we analyzed 
in the proof of Proposition~\ref{prop-lower-bd}; otherwise, the integration may have an extra logarithmic 
factor when the maximum is attained with $s = s_q= \lceil 2/(q-1) \rceil \wedge (d+1)$. 
In the simplest case where $\ntil_j$ are all equal, $1=t_1=\cdots=t_d$, $s=d$ is the effective dimension of the lattice 
$[\ba,\bb]$ and 
\bes
\Pi = \big( \Delta^*_n/\ntil^{1/d}\big)\bigg(1+(1/2 + 1/d)\int_1^{t_*} t^{1/d-1/2-q/2}dt\bigg)
\ees
with $t_*= (\ntil^{1/d}/\Delta^*_n)^{2d/(2+d)}\wedge \ntil$.  Here are the details for the general case. 

For $1 \le t_*\le t_{s_q}$, or equivalently $t^*\in [t_s,t_{s+1}]$ for some $1\le s <s_q$, 
\bes
\Pi \lesssim_{q,d} \Htil(t^*)/(t^*)^{q/2} = \max_{1\le t \le n}t^{-q/2}\Htil(t)
\ees
as in the second case of (\ref{pf-prop-lower-bd-2}),
or the second case of (\ref{th2-1}), 
or the second case of (\ref{prop-lower-bd-1}).

Similarly, for $t_{s_q} \le t_*\le \ntil$ and $2/(q-1)\not\in \{1,\ldots,d\}$,
\bes
\Pi \le  \int_0^{t_{s_q}} t^{-q/2}\Htil(dt) + \int_{t_{s_q}}^{t_*} t^{-q/2}\Htil(dt)
\lesssim_{q,d} t_{s_q}^{-q/2}\Htil(t_{s_q}) = \Delta^*_n t_{s_q}^{-q/2+1/2}/\ntil_{s_q}
\ees
as in the third case of (\ref{pf-prop-lower-bd-2}),
or the third case of (\ref{th2-1}) with $\Lambda_{s_q}=1$, 
or the third case of (\ref{prop-lower-bd-1}). 

Finally, for $t_{s_q} \le t_*\le \ntil$ and $2/(q-1) = s_q \in \{1,\ldots,d\}$, 
the integration of $t^{-1}$ in the critical interval $[t_{s_q},t_{s_q+1}\wedge t_*]$ 
may result in an extra logarithmic term. With $s=s_q$ 
\bes
\int_{t_{s}}^{t_*\wedge t_{s+1}} t^{-q/2}\Htil(dt)
\le \frac{(1/2+1/s)\Delta^*_n}{(1/s)(\ntil_s^*)^{1/s}}\log_+\bigg(\bigg(\frac{t_*\wedge t_{s+1}}{t_{s}}\bigg)^{1/s}\bigg)
= \frac{(s+2)\Delta^*_n\Lambda_s}{2(\ntil_s^*)^{(q-1)/2}}
\ees
in view of (\ref{lambda}), as $(t_{s+1}/t_s)^{1/s}=\ntil_s/\ntil_{s+1}$ and 
$(t_*/t_s)^{1/s}=((\ntil_{s}^*)^{1/s}/\Delta^*_n)^{2/(2+s)}\ntil_s/(\ntil^*_s)^{1/s}$.
For $t_*>t_{s_q+1}$, $\int_{t_{s_q+1}}^{t_*} t^{-q/2}\Htil(dt) 
\lesssim_{q}\Htil(t_{s_q+1})/t_{s_q+1}^{q/2}=\Htil(t_{s_q})/t_{s_q}^{q/2}$ 
as $\Htil(t)/t^{q/2}=\Delta_n^*/(\ntil^*_{s_q})^{1/s_q}$ is 
a constant in $[t_{s_q},t_{s_q+1}]$ for $q/2=1/s+1/2$. Therefore, 
\bes
  \Pi \lesssim_{q} \frac{\Delta^*_n\Lambda_s}{(\ntil_s^*)^{(q-1)/2}} 
= \frac{\Delta^*_n\Lambda_s}{t_s^{(q-1)/2}\ntil_s},\quad s=s_q=2/(q-1)\le d,
\ees
which gives the third case of (\ref{th2-1}) when the integration of $t^{-1}$ is involved.  

$\hfill\square$

\medskip
{\bf A3.4. Proof of Theorem \ref{th-fix-worst}.}
It follows from (\ref{prop-lower-bd-1}) of Proposition \ref{prop-lower-bd} that
\bes
&&  \sigma^q\, \max
\Big\{(t\wedge n)^{-q/2}H(t):\ t \wedge h_0(t) \ge 1\Big\}
\cr &\lesssim_{q,d}&
\inf_{\hbf}\, \sup\Big\{ R_q(\hbf, \bff_n):
\bff_n \in\scrF_n, \Delta(\bff_n/\sigma)\le\Delta_n^*\Big\}.
\ees
Thus, the main claim (\ref {th-fix-worst-0}) holds when
\bes
&& \sup\Big\{ R_q(\hbf^{(block)}_n, \bff_n):
\bff_n \in\scrF_n, \Delta(\bff_n/\sigma)\le\Delta_n^*\Big\}
\cr &\lesssim_{q,d}& \,\sigma^q\bigg( H(1) + \int_1^{n}\frac{H(dt)}{t^{q/2}}\bigg)
+ \frac{\sigma^{q}}{n}\prod_{j=1}^d \int_0^{n_j}\frac{dt}{(t\vee 1)^{q/2}}
\cr &\lesssim_{q,d}& \Lambda_{s_q}
\sigma^q\, \max
\Big\{(t\wedge n)^{-q/2}H(t):\ t \wedge h_0(t) \ge 1\Big\}
+ \Big(\frac{\sigma^q}{n}\prod_{j=1}^d \log_+(n_j)\Big)^{I\{q=2\}}.
\ees
The first inequality above is (\ref{th2-main}) in Theorem \ref{th-lattice-general} with $[\ba, \bb] = [{\bf 0}, {\bf 1}]$.
The second follows from a comparison between the
upper bound (\ref{th2-1}) in Theorem \ref{th-lattice-general}
with $[\ba, \bb] = [{\bf 1}, \bn]$ and the lower bound (\ref{prop-lower-bd-1})
in the respective scenarios, covering $\Delta_n^*\ge t_{s_q}^{-1/2} = \Big(\prod_{j=1}^{s_q}(n_j/n_{s_q})\Big)^{-1/2}$.

The rate in (\ref {th-fix-worst-1}) follows directly from (\ref{th2-1}) with $[\ba, \bb] = [{\bf 1}, \bn]$. Note when $2/(q-1)=s_q \le d-1$, $n_{s_q+1}/n_{s_q} \asymp 1$ but $n_{d}/n_{d+1} \asymp n^{1/d}$.  $\hfill\square$

\bigskip
{\bf A4. Proofs of the results in Subsection 3.4}

\medskip
{\bf A4.1. Proof of Theorem \ref{coro-fix-piecewise}.}
As $\Delta_{\ba_k,\bb_k}=0$ for all $k$, it follows from Theorem \ref{th-fix-worst} and (\ref {th2-2}) that
\bes
T_{q}([\ba_k,\bb_k])
\lesssim_{q,d} \sigma^q \bigg[n_{\ba_k, \bb_k}^{1-q/2} + \Big(\prod_{j=1}^d \log_+(b_{k,j} -a_{k,j}+1)\Big)^{I\{q=2\}}\bigg],
\ees
where $a_{k,j}$ and $b_{k,j}$ are the $j_{th}$ element of $\ba_k$ and $\bb_k$ respectively. 
The first conclusion follows. 
When $1\le q < 2$, we have
\bes
T_q(V) \lesssim_{q, d} \sigma^q  \sum_{k=1}^K n_{\ba_k, \bb_k}^{1-q/2} \lesssim_{q, d} \sigma^q  K (n/K)^{1-q/2}
\ees
and when $q=2$,
\bes
T_2(V) &\lesssim_d& {\sigma^2} \sum_{k=1}^K \log_+^{d_K} \big(\max_j (b_{k,j} -a_{k,j})+1 \big)
\cr
&\lesssim_d&  {\sigma^2} \sum_{k=1}^K \log_+^{d_K}( n_{\ba_k,\bb_k}\big)
\cr
&\lesssim_d& \sigma^2 {K} \log_+^{d_K}(\frac{1}{K} \sum_{k=1}^K n_{\ba_k,\bb_k}\big),
\ees
where the last inequality follows as $\log_+^{d_K}(x)$ is a concave function when $x$ is greater than a certain constant $C_{d}$. The second conclusion follows as $q>2$ simply gives rate $K/n$.  $\hfill\square$

\bigskip
{\bf A5. Proofs of the results in Subsection 3.5}

\medskip
{\bf A5.1. Proof of Theorem \ref{th-fix-selection}.}
As $f(\bx)=f_S(\bx_S)$, we can always take $\bv$ with the largest $\bv_{S^c}$,
so that
\bes
r_{q,+}(m_{\bx}) = r_{q,+}(m_{S,\bx_S})(C_{q,d} / C_{q,s})n_{\bx_{S^c},\bn_{S^c}}^{-q/2}
\ees
where $r_{q,+}(m_{S,\bx_S})$ is the risk bound at $\bx_S$ in model $S$,
$n_{\bx_{S^c},\bn_{S^c}}$ is the size of $[\bx_{S^c},\bn_{S^c}]$ in model $S^c$, and $C_{q,d}$ is from the definition of $r_{q,+}(m)$ as in (\ref {r-bd}). We note $C_{q,s} \le C_{q,d}$ for all $s \le d$.
Thus, in the sheet of $\bx$ with fixed $\bx_{S^c}$, the risk bound is identical
to that of model $S$ with $\sigma^q$ reduced by a factor $n_{\bx_{S^c},\bn_{S^c}}^{-q/2}$.
Let $\sigma_{\bx_{S^c}}^q =(C_{q,d}/C_{q,s}) \sigma^q /n_{\bx_{S^c},{\bn}_{S^c}}^{q/2}$.

\bes
&& T_q(V)
\cr &\lesssim_{q,d}& \sum_{\bx_{S^c}}n^{s/d}
\sigma_{\bx_{S^c}}^q
\min\bigg\{1, \Big(\Delta(\bff_n/\sigma_{\bx_{S^c}})n^{-1/d}\Big)^{\min\{1, qs/(s+2)\}} (\log n)^{I\{qs=s+2\}}
\cr
&& \qquad\qquad +\,\big(n^{s/d}\big)^{-\min\{1,  q/2\}} (\log n)^{sI\{q=2\}} \bigg\}
\cr &\lesssim_{q,d}& \sum_{\bj\in[1,n^{1/d}]^{d-s}}n^{s/d}\sigma^q(j_1\cdots j_{d-s})^{-q/2}
\cr &&\qquad \times
\min\bigg\{1, \Big(\Delta(\bff_n/\sigma)
(j_1\cdots j_{d-s})^{1/2}n^{-1/d}\Big)^{\min\{1, qs/(s+2)\}}(\log n)^{I\{qs=s+2\}}\cr
&& \qquad\qquad +\,\big(n^{s/d}\big)^{-\min\{1,  q/2\}} (\log n)^{sI\{q=2\}}\bigg\},
\ees
where the first inequality follows from (\ref {th-fix-worst-1}) in Theorem \ref{th-fix-worst}.
Hence we obtain (\ref {th-fix-selection-1}) and (\ref {th-fix-selection-2}). $\hfill\square$

\bigskip
{\bf A6. Proofs of the results in Subsection 3.6}

\medskip
{\bf A6.1. Proof of Proposition \ref{prop-point-risk}}. 
We first prove that \eqref{r-plus-new} holds with the variability bound in \eqref{r-plus-random}. 
As the $V=[0,1]^d$ for the random design, we modify the proof of \eqref{pf2-1-2} in the proof of Proposition~\ref{prop-single-point} as follows. 

Consider fixed ${\bf 0}\le\ba\le\bb\le {\bf 1}$ and $\bx\in [\ba,\bb]$. 
The modification involves vectors $\bv$ and $\bw$ in $\R^d$ 
with $\bx\le\bw\le \bv\le\bb$, 
and a probability measure $\P^*$ which endows the same conditional distribution of 
$\{\veps_i, i\le n\}$ given $\{\bx_i, i\le n\}$ as $\P$ does 
and iid $\bx_i\in [0,1]^d$ distributed according to the conditional 
distribution under $\P$ given $n_{\bx,\bw}=0$. 
This covers three cases, $\bv=\bv_{\bx}$ and $\bw=\bx$ for $\P^*=\P$, 
$\bv=\bb$ and $\bw=\bv_{\bx}$ for the conditional probability given $n_{\bx,\bv_{\bx}}=0$, 
and $\bv={\bf 1}$ and $\bw=\bb$ for the conditional probability given $n_{\bx,\bb}=0$. 

Consider the case of 
$\P\big\{\bx_i\in [\bx,\bv]\setminus [\bx,\bw]\big\}>0$. 
Let $\mu^*$ be the measure in $\R^d$ given by 
\bes
\mu^*([\bu,\bv]) = \mu^L([\bu,\bv]\setminus [\bx,\bw])/\mu^L([{\bf 0},{\bf 1}]\setminus [\bx,\bw]). 
\ees
For $t\ge 0$, let $h(t)$ be the function satisfying 
\bes
h(t)\mu^L([\bx,\bv]) - \mu^L([\bx,\bw]) = e^{t}\mu^L\big([\bx,\bv]\setminus [\bx,\bw]\big). 
\ees
As $\P\big\{\bx_i\in [\bx,\bv]\setminus [\bx,\bw]\big\} >0$,  
$\mu^L([\bx,\bv]\setminus [\bx,\bw]>0$ by \eqref{P-condition}. 
We have $h(0)=1$ and 
\bes
(d/dt)\log h(t) = \frac{\mu^L\big([\bx,\bv]\setminus [\bx,\bw]\big)}
{\mu^L\big([\bx,\bv]\setminus [\bx,\bw]\big) + e^{-t}\mu^L([\bx,\bw)}. 
\ees
Thus, $\log h(t)$ is convex in $t$. 
To modify the construction of $U_{\bx, \bv}$ in (\ref{Uxv}), 
we define $\ubu_{\bi}$ by
\bes
\big(\ubu_{\bi}\big)_{(j)} = v_{(j)} - h_{(j)}(\bi)\big(v_{(j)}-x_{(j)}\big). 
\ees
where $(j)$ is determined by $i_{(1)}\le \ldots \le i_{(d)}$ and $h_{(j)}(\bi)$ 
are given by 
\bel{h_{(j)}}
 h_{(j)}(\bi) = \bigg(\frac{h(i_{(1)}+\cdots+i_{(j)} + (d-j)i_{(j)})}
 {h_{(1)}(\bi) \times \cdots \times h_{(j-1)}(\bi)}\bigg)^{1/(d+1-j)}
\eel
with $h_{(0)}(\bi)=1$. More explicitly, 
with $g_{(j)}(\bi)=\log h(i_{(1)}+\cdots+i_{(j)} + (d-j)i_{(j)})$, 
\bel{g_{(j)}}
\log h_{(j)}(\bi) = \frac{g_{(j)}(\bi)}{d+1-j}-\sum_{\ell=1}^{j-1}\frac{g_{(\ell)}(\bi)}{(d+1-\ell)(d-\ell)}
\eel
Furthermore, we define $\obu_{\bi} = \ubu_{\bi-{\bf 1}}$, and 
\bes
U_{\bx,\bv} = \Big\{[\ubu_{\bi},\obu_{\bi}]: \obu_{\bi}\in [{\bf 0},\bx], i_j > 0\ \forall j\Big\}. 
\ees
Next, we verify the following properties of $U_{\bx,\bv}$, 
\bel{new-cover}
\frac{\mu^*([\ubu_{\bi},\bv])}{\mu^*([\bx,\bv])} = e^{|\bi|},\quad  
\cup\Big\{[\ubu_{\bi},\obu_{\bi}]: [\ubu_{\bi},\obu_{\bi}]\in U_{\bx,\bv}\Big\}\supseteq [{\bf 0},\bx]. 
\eel
As $\prod_{j=1}^d h_{(j)}(\bi) = h(|\bi|)$, the first part above follows from 
\bes
\frac{\mu^*([\ubu_{\bi},\bv])}{\mu^*([\bx,\bv])}
=\frac{\mu^L([\ubu_{\bi},\bv]) - \mu^L([\bx,\bw])}{\mu^L\big([\bx,\bv]\setminus [\bx,\bw]\big)} 
=\frac{h(|\bi|)\mu^L([\bx,\bv]) - \mu^L([\bx,\bw])}{\mu^L\big([\bx,\bv]\setminus [\bx,\bw]\big)} 
=e^{|\bi|}. 
\ees
For the second part of \eqref{new-cover}, 
it suffices to consider $1 \le i_1\le \dots \le i_{j-1}<i_j = \cdots = i_d$ 
by symmetry and prove that the block $[\ubu_{\bi},\obu_{\bi}]$ overlaps with blocks $[\ubu_{\bi'},\obu_{\bi'}]$ with $\bi\neq \bi'\le\bi$ on sides $j,\ldots,d$. 
This holds iff $\ubu_{\bi'}\le \ubu_{\bi-{\bf 1}} = \obu_{\bi}$  
with $\bi'=(i_1,\ldots,i_{j-1},i_j-1,\ldots,i_d-1)^T$, iff $h_{(j)}(\bi-{\bf 1}) \le h_{(j)}(\bi')$, 
and by \eqref{g_{(j)}} iff 
\bes
\frac{g_{(j)}(\bi-{\bf 1})}{d+1-j}-\sum_{\ell=1}^{j-1}\frac{g_{(\ell)}(\bi-{\bf 1})}{(d+1-\ell)(d-\ell)}
\le \frac{g_{(j)}(\bi')}{d+1-j}-\sum_{\ell=1}^{j-1}\frac{g_{(\ell)}(\bi')}{(d+1-\ell)(d-\ell)}, 
\ees
which can be written as 
\bes
\int_{d(i_1-1)}^{i_1+\ldots+i_{j-1}+(d+1-j)i_{j-1}}   
g_0(t) d\log h(t)
\le \int_{i_1+\ldots+i_j+(d-j)i_j-d}^{i_1+\ldots+i_{j-1}+(d-j+1)(i_j-1)} 
\frac{d\log h(t)}{d+1-j}
\ees
with $g_0(t)=\sum_{\ell=1}^{j-1} 
I\{i_1+\ldots+i_\ell+(d-\ell)i_\ell-d\le t\le i_1+\ldots+i_\ell+(d-\ell)i_\ell\}/(d+1-\ell)(d-\ell)$. 
This inequality holds as $g_0(t)$ is decreasing in $t$  
in the domain of the integration, $\log h(t)$ is convex in $t$  
and the two sides are equal when $d\log h(t)$ is replaced by the Lebesgue measure $dt$. 
For example, for $j=d=2$ and $i_1<i_2$, 
$\ubu_{i_1,i_2-1} \le \ubu_{i_1-1,i_2-1}$ follows from 
\bes
2\int_{i_1+i_2-2}^{i_1+i_2-1} d\log h(t) \ge \int_{2i_1-2}^{2i_1} d \log h(t)
\ees
by direct computation from \eqref{h_{(j)}}. 
This completes the proof of \eqref{new-cover}.  

Similar to the calculation below (\ref{Uxv}), it holds that, when $\P\big\{\bx_i\in [\bx,\bv]\setminus [\bx,\bw]\big\}>0$, or equivalently $\mu_{\bx, \bv} > \mu_{\bx, \bw}$, 
\bes
&&  \E^* 
\bigg(\max_{ \bu \le \bx} \sum_{\bx_i \in [\bu,\bv]} \frac{\veps_i}{n_{\bu,\bv}\vee 1} \bigg)^q_{+}
\\
\nonumber &\le& \sum_{[\ubu_{\bi},\obu_{\bi}]\in U_{\bx,\bv} }   \E^*  \bigg(\max_{\bu \in [\ubu_{\bi},\obu_{\bi}] }\Big( \frac{1}{n_{\obu_{\bi}, \bv} \vee 1} \Big|\sum_{\bx_i \in [\bu,\bv]} {\veps_i} \Big|\Big)^q \bigg)
\cr 
&\le& \sum_{[\ubu_{\bi},\obu_{\bi}]\in U_{\bx,\bv} }  (C'_q)^d 
\Bigg[ \E^*  \bigg( \frac{1}{n_{\obu_{\bi}, \bv}\vee 1} 
\Big|\sum_{\bx_i \in [\ubu_{\bi},\bv]} {\veps_i} \Big| \bigg)^{q \vee 2}\Bigg]^{q/(q \vee 2)}
\cr 
&\le& \sum_{[\ubu_{\bi},\obu_{\bi}]\in U_{\bx,\bv} }  (C'_q)^d 
\Bigg[ \E^* \big[(n_{\obu_{\bi}, \bv}\vee 1)^{-2(q\vee 2)}\big]
 \E^*  \bigg|\sum_{\bx_i \in [\ubu_{\bi},\bv]} {\veps_i} \bigg|^{2(q \vee 2)}\Bigg]^{(q/2)/(q\vee 2)}
\cr 
&\le& \sum_{[\ubu_{\bi},\obu_{\bi}]\in U_{\bx,\bv} } (C'_q)^d C_q'' \sigma^q 
\Bigg[ \E^* \big[(n_{\obu_{\bi}, \bv}\vee 1)^{-2(q\vee 2)}\big]
 \E^* \big[n_{\ubu_{\bi}, \bv}^{(q\vee 2)}\big]\Bigg]^{(q/2)/(q \vee 2)}, 
\ees
where the first inequality follows from the second part of \eqref{new-cover}, the second from Lemma \ref{lm-doob}, the third from Cauchy-Schwarz inequality, 
and the fourth from Rosenthal's inequality. 
As $n_{\bu,\bv}$ is a binomial random variable, 
\bes
&& \E^*[(1\vee n_{\bu,\bv})^{\pm q}] \lesssim_q [1\vee \big(\E^*[n_{\bu,\bv}]\big)^{\pm q}],\quad  \forall q>0, 
\cr && \E^*[n_{\bu,\bv}^q] \lesssim_q \E^*[n_{\bu,\bv}],\ \hbox{ when $q > 0$ and $\E^*[n_{\bu,\bv}] \le 1$.}
\ees
Moreover, by \eqref{P-condition}  and the first part of \eqref{new-cover}, 
\bes
\frac{\E^*\big[n_{\ubu_{\bi}, \bv}\big]}{\E^*\big[n_{\obu_{\bi}, \bv}\big]}
\lesssim_{\rho_1,\rho_2} \frac{\mu^*([\ubu_{\bi}, \bv])}{\mu^*([\obu_{\bi}, \bv])} = e^d,\quad 
\E\big[n_{\obu_{\bi}, \bv}\big] 
\asymp_{\rho_1,\rho_2} \mu^*([\obu_{\bi}, \bv]) = e^{|\bi|-d}\mu^*([\bx, \bv]). 
\ees
Thus, when $\mu_{\bx, \bv} > \mu_{\bx, \bw}$,
\bel{pf3-1-1}
&& \E^* 
\bigg( \max_{ \bu \le \bx} \sum_{\bx_i \in [\bu,\bv]} \frac{\veps_i}{n_{\bu,\bv}\vee 1} \bigg)^q_{+}
\le C_{q,d, \rho_1,\rho_2} \sigma^q \big( 1\vee \E^*[n_{\bx, \bv}] \big) ^{-q/2}. 
\eel
By arguments similar to what we have in the proof of Theorem \ref{th-lattice-general}, 
we can maintain $C_{q,d,\rho_1,\rho_2}$ continuous in $q$. This gives \eqref{r-plus-random} with $\P^* = \P$ and $\bw = \bx$ when $\mu_{\bx, \bv} > \mu_{\bx, \bw} = 0$. It also holds for $\mu_{\bx, \bv}=0$ as the bound $r_{q, +}(n\mu_{\bx, \bv})$ can be right-continuous due to
\bes
&& \E^* \bigg( \max_{ \bu \le \bx} \sum_{\bx_i \in [\bu,\bv]} \frac{\veps_i}{n_{\bu,\bv}\vee 1} \bigg)^q_{+}
\cr
&\le& \E^* \bigg( \max_{ \bu \le (1- \epsilon) \bx} \sum_{\bx_i \in [\bu,\bv]} \frac{\veps_i}{n_{\bu,\bv}\vee 1} \bigg)^q_{+}  + \E^* \bigg( \max_{ \bu \in [{\bf 0}, \bx] \setminus [{\bf 0}, (1-\epsilon)\bx]} \sum_{\bx_i \in [\bu,\bv]} \frac{\veps_i}{n_{\bu,\bv}\vee 1} \bigg)^q_{+},
\ees 
where the second term goes to zero when $\epsilon \to 0+$.

We then prove \eqref{point-risk} by considering three cases that can happen when estimating $f(\bx)$: (1) the block $[\bx, \bv_{\bx}]$ is non-empty; (2) the block $[\bx, \bv_{\bx}]$ is empty but not $[\bx, \bb]$; and (3) the block $[\bx, \bb]$ is empty.

Consider the first case where $n_{\bx, \bv_{\bx}} >0$.
By the definition of $m_{\bx}$, there exists $\bv_{\bx}\le\bb$ such that in the event $n_{\bx, \bv_{\bx}} >0$,
\bes
\fhat_n^{(block)}(\bx) \le \max_{\bu \le \bx} \sum_{\bx_i\in [\bu,\bv_{\bx}]} \frac{y_i}{n_{\bu,\bv_{\bx}}}
\le f(\bx) + r_{q,+}^{1/q}(m_{\bx}) + \max_{\bu \le \bx} \sum_{\bx_i\in [\bu,\bv_{\bx}]} \frac{\veps_i}{n_{\bu,\bv_{\bx}}}
\ees
with the $r_{q,+}(m)$ in \eqref{r-plus-new}. Thus, by \eqref{pf3-1-1} 
with $\P^*=\P$, $\bv=\bv_{\bx}$, $\bw = \bx$ and $\mu_{\bx, \bv_{\bx}} > 0$ (otherwise $\P\{n_{\bx, \bv_{\bx}} >0 \}=0$),  
\bel{point-risk-1}
&&\E\big(\fhat_n^{(block)}(\bx) - f(\bx)\big)_+^q I\{n_{\bx, \bv_{\bx}}>0\} \le 2^q r_{q,+}(m_{\bx}). 
\eel

Consider the second case where $n_{\bx, \bv_{\bx}} =0$ but $n_{\bx, \bb} > 0$. It follows that
\bes
\fhat_n^{(block)}(\bx) \le  f(\bb) + \max_{\bu\le \bx} \sum_{\bx_i\in [\bu,\bb]} \frac{\veps_i}{n_{\bu,\bb}},
\ees
so that by \eqref{pf3-1-1} with $\P^*$ being the conditional probability under $n_{\bx,\bv_{\bx}}=0$, $\bw=\bv_{\bx}$, $\bv=\bb$ and $\mu_{\bx, \bb} > \mu_{\bx, \bv_{\bx}}$ (otherwise $\P\{n_{\bx, \bv_{\bx}}=0,\, n_{\bx, \bb}>0 \} = 0$),  
\bel{point-risk-2}
&&\E\big(\fhat_n^{(block)}(\bx) - f(\bx)\big)_+^q I\{n_{\bx, \bv_{\bx}}=0,\, n_{\bx, \bb}>0\}
\\ \nonumber 
&\le& 2^{q-1} \bigg(\big(f(\bb) - f(\bx)\big)^q  + 
\E \bigg[\bigg(\max_{\ba \le \bu \le \bx  } \sum_{\bx_i \in [\bu, \bb]} 
\frac{\veps_i}{n_{\bu, \bb} \vee 1}\bigg)_+^q\bigg|n_{\bx, \bv_{\bx}} =0\bigg] \bigg) 
\P\{n_{\bx, \bv_{\bx}} =0\}
\cr
&\le& 2^{q-1} C_{q, d, \rho_1, \rho_2} \sigma^q \big(\Delta_{\ba, \bb}^q + 1\big) e^{-m_{\bx}}.
\eel

Finally, we consider the third case where $n_{\bx, \bb}=0$. By the definition of $\fhat_n^{(block)}$,
$\fhat_n^{(block)}(\bx) \le  f({\bf 1}) + \max_{\bu\le \bx} \sum_{\bx_i\in [\bu,{\bf 1}]} \veps_i/n_{\bu,{\bf 1}}$.
Similar to the first two cases,
\bel{point-risk-3}
\E\big(\fhat^{(block)}_n(\bx) - f(\bx)\big)_+^q I\{n_{\bx, \bb}=0\}
\le 2^{q-1} C_{q, d, \rho_1, \rho_2} \sigma^q \big(\Delta_{{\bf 0}, {\bf 1}}^q  + 1 \big)e^{-n\mu_{\bx, \bb}}
\eel
for $\mu_{\bx,{\bf 1}}> \mu_{\bx, \bb}$. It remains true for $\mu_{\bx,{\bf 1}} = \mu_{\bx, \bb}$ by a similar right-continuity argument to the one below \eqref{pf3-1-1}; we omit the details.
Therefore (\ref {point-risk}) follows from (\ref {point-risk-1}), (\ref {point-risk-2}) and  (\ref {point-risk-3}).
$\hfill\square$

\medskip
{\bf A6.2. Proof of Theorem \ref{th-random}.} In this proof, we may re-define some notation to fit in with the random design scenario. Such new notation supersedes definitions elsewhere, but is applicable only in this proof. 
Throughout the proof $r_{q,+}(m)$ is defined as in \eqref{r-plus-random}. 

As the risk over block $[\ba, \bb]$ is
\bes
R^*_{q}([\ba,\bb]) = \int_{[\ba,\bb]} \E \big|\fhat^{(block)}_n(\bx) - f(\bx)\big|^q d\bx,
\ees
it suffices by symmetry to only bound $R^*_{q,+}([\ba,\bb])$, where 
\bes
R^*_{q,+}(A) = \int_{A} \E \big(\fhat_n^{(block)}(\bx) - f(\bx)\big)_+^q d\bx,\quad A\subseteq [{\bf 0},{\bf 1}]. 
\ees
This will be done through Proposition \ref{prop-point-risk}. 
A direct consequence of Proposition \ref{prop-point-risk} is
\bes
R^*_{q}([\ba,\bb]) \lesssim_{q,d,\rho_1, \rho_2} \sigma^q \big(\Delta^q_{{\bf 0}, {\bf 1}}+1)\mu_{\ba, \bb},
\ees
which serves as the trivial upper bound in (\ref{th-random-rate}).

Parallel to the proof of Theorem \ref{th-lattice-general}, we partition $V_0 = [\ba, \bb]$ into $V_{0,+}$ and $V_0 \setminus V_{0,+}$, where $V_{0,+} = \{\bx \in [\ba, \bb]: f(\bb) > f(\bx) + r_{q,+}^{1/q}(n\mu_{\bx, \bb})\}$. 
By \eqref{m_x-random}, $m_{\bx} \le n\mu_{\bx, \bb}$ for $\bx \in V_{0,+}$ (equality may hold only if $f(\bx)$ is not continuous at $\bx=\bb$), 
and $m_{\bx} = n\mu_{\bx,\bb}$ for $\bx \notin V_{0,+}$.
In what follows, we first bound $R^*_{q,+}(V_{0,+})$ and then $R^*_{q,+}(V_0\setminus V_{0,+})$. The conclusions follow from summing the two bounds up. 

To derive bound for $R^*_{q,+}(V_{0,+})$, which is an integral over $V_{0,+}$, we first integrate over lines parallel to $\bb-\ba$ and starting from points in the lower-half boundary of $[\ba, \bb]$, and integrate them over the lower-half boundary. Formally, let $\pa_{Lower} = \{ \bx \in [\ba, \bb]: x_j = a_j \hbox{ for some } j\}$ denote the lower-half boundary of $[\ba, \bb]$. We define for each $\bc \in \pa_{Lower}$, the anti-diagonal line segment starting from $\bc$ as $L_{\bc} = \big\{\bc + k(\bb-\ba): k \in [0,k_{\bc}]\big\}$, where $k_{\bc} = \sup \big\{k: \bc+k(\bb-\ba) \in V_{0,+}\}$ is the length of the segment.
It follows from the above definitions that $k_{\bc} \le 1$, $L_{\bc} \subset V_{+}$ and
\bes
V_{0,+} = \bigcup_{\bc \in \pa_{Lower}} L_{\bc}.
\ees

For simplicity, let $m_{\bc, k} = m_{\bc+k(\bb-\ba)}$ and $g(k) = f(\bc + k(\bb-\ba))$.
Observe
\bes
&& R^*_{q,+}(L_{\bc})
\cr
&\le& C_{q,d,\rho_1, \rho_2}'  \sigma^q \int_{L_{\bc}} \Big((m_{\bx} \vee 1)^{-q/2} + \Delta_{\ba, \bb}^q e^{-m_{\bx}}\Big) d\bx + C_q' \sigma^q \Delta_{{\bf 0}, {\bf 1}}^q \int_{L_{\bc}} e^{-n\mu_{\bx, \bb}} d\bx
\cr
&\le& C_{q,d,\rho_1, \rho_2}' \sigma^q \int_0^{k_{\bc}} \Big((m_{\bc, k} \vee 1)^{-q/2} + \Delta_{\ba, \bb}^qe^{-m_{\bc,k}}\Big) dk + C_q' \sigma^q \Delta_{{\bf 0}, {\bf 1}}^q \int_{L_{\bc}} e^{-n\mu_{\bx, \bb}} d\bx
\cr
&\le& C_{q,d,\rho_1, \rho_2}'\sigma^q \int_0^{n\mu_{\ba, \bb}} \Big( (m \vee 1)^{-q/2} + \Delta_{\ba, \bb}^q e^{-m}\Big) d\ell_{\bc, +}(m) + C_q' \sigma^q \Delta_{{\bf 0}, {\bf 1}}^q \int_{L_{\bc}} e^{-n\mu_{\bx, \bb}} d\bx,
\ees
where $\ell_{\bc, +}(m) = \int_{0}^{k_{\bc}} I\{m_{\bc, k} < m \}dk$.

It then suffices to bound $\ell_{\bc, +}(m)$. To this end, we shall divide divide $V_{0,+}$ into
\bes
D_j = \Big\{\bx \in V_{0,+}: f(\ba) + (j-1) r_{q,+}^{1/q}(m) \le  f(\bx)
< f(\ba) + j r_{q,+}^{1/q}(m)\Big\},
\ees
$j=1, \ldots, J$, where
$J=\big\lceil \{f(\bb)- f(\ba)\}/r_{q,+}^{1/q}(m)\big\rceil$.
Consider $D_j \cap L_{\bc}$ and let $\bv$ be the right end point of this segment, i.e., $\bx \le \bv$ for all $\bx \in D_j \cap L_{\bc}$. If we can find $\bx_{\bv} \in D_j \cap L_{\bc}$ such that
$n\mu_{\bx_{\bv}, \bv} =m$, then any point $\bx \le \bx_{\bv}$ in $D_j \cap L_{\bc}$ has $m_{\bx} \ge m$. Let $\bv = \bx_{\bv} + t(\bb - \ba)$. It follows that
\bes
\int_{\bc + k{\bf 1} \in D_j \cap L_{\bc}} I\{m_{\bc, k} < m \}dk \le t \le \Big[(\rho_2/\rho_1)\frac{m}{n\mu_{\ba, \bb}}\Big]^{1/d}.
\ees
The above bound is trivial if there is no such $\bx_{\bv}$.
For $\bx \in D_J \cap L_{\bc} $, we have $f(\bb) \le f(\bx) + r_{q,+}^{1/q}(m)$ by the definition of $D_J$ and $f(\bb) > f(\bx) + r_{q,+}^{1/q}(n\mu_{\bx, \bb})$ by the definition of $V_{0,+}$,
which implies $m < n\mu_{\bx, \bb}$.
However, $m_{\bx} \ge n\mu_{\bx, \bb}$ due to $\bx \in V_{0,+}$ so that $m_{\bx} > m$ and
\bes
\int_{\bc + k{\bf 1} \in D_J \cap L_{\bc}} I\{m_{\bc, k} < m \}dk  =0.
\ees
Overall, we have
\bes
\ell_{\bc, +}(m) &=& \sum_{j=1}^J \int_{\bc + k{\bf 1} \in D_j \cap L_{\bc}} I\{m_{\bc, k} < m \}dk
\cr
&\le& \Big[(\rho_2/\rho_1)\frac{m}{n\mu_{\ba, \bb}}\Big]^{1/d} \frac{f(\bb)- f(\ba)}{r_{q,+}^{1/q}(m)}
\cr
&\le& (\rho_2/\rho_1)^{1/d}C_{q,d,\rho_1,\rho_2}^{-1/q}\Delta_n^* (n\mu_{\ba, \bb})^{-1/d} m^{1/2+1/d}
\cr
&\le& (\rho_2/\rho_1)^{1/d}C_{q,d,\rho_1,\rho_2}^{-1/q} H^*(m),
\ees
where $H^*(m) = \min\Big\{1, \Delta_{\ba, \bb}(n\mu_{\ba, \bb})^{-1/d}m^{1/2+1/d}\Big\}$.

Consequently,
\bel{Rq+V0+}
 R^*_{q,+}(V_{0,+})
&=& \int_{\bc \in \pa_{Lower}} R_{q,+}(L_{\bc})  d\bc
\\
\nonumber  &\le& C_{q,d,\rho_1, \rho_2}''\sigma^q \int_0^{n\mu_{\ba, \bb}} \Big( (m \vee 1)^{-q/2} + \Delta_{\ba, \bb}^q e^{-m}\Big) H^*(dm)
\cr
&& \quad +\, C_q' \sigma^q \Delta_{{\bf 0}, {\bf 1}}^q \int_{V_{0,+}} e^{-n\mu_{\bx, \bb}} d\bx\Big\}.
\eel

We then bound $R^*_{q,+}([\ba,\bb]\setminus V_{0,+})$. As $f(\bb) - f(\bx) \le r_{q,+}^{1/q}(m_{\bx})$ and $m_{\bx} = n\mu_{\bx, \bb}$ for $\bx \in [\ba, \bb]\setminus V_{0,+}$, it follows from Proposition \ref{prop-point-risk} that
\bes
\E \big(\fhat_n(\bx) - f(\bx) \big)_+^q  \le  2^{q} r_{q,+}(n\mu_{\bx, \bb}) + 2^q \sigma^q(\Delta_{{\bf 0}, {\bf 1}}^q +1) e^{-n\mu_{\bx, \bb}}.
\ees

Therefore
\bel{Rq+res}
&& R^*_{q,+}([\ba,\bb] \setminus V_{0,+})
\\
\nonumber  &\le& C_{q,d, \rho_1,\rho_2}''' \sigma^q \int_{\bx \in [\ba, \bb] \setminus V_{0,+}} \Big( \big( (n \mu_{\bx, \bb} )\vee 1\big)^{-q/2} + \Delta_{{\bf 0}, {\bf 1}}^q  e^{-n\mu_{\bx, \bb}} \Big) d\bx
 \cr
 &\le& C_{q,d, \rho_1,\rho_2}''' \sigma^q \int_{\bx \in [\ba, \bb] } \Big( \big( (n \mu_{\bx, \bb} )\vee 1\big)^{-q/2} + \Delta_{{\bf 0}, {\bf 1}}^q  e^{-n\mu_{\bx, \bb}} \Big) d\bx.
\eel
The main conclusion (\ref {th-random-main}) directly follows from (\ref {Rq+V0+}) and (\ref {Rq+res}), with appropriately chosen $C_{q,d,\rho_1,\rho_2}^*$ so that it remains continuous in $q \ge 1$ and non-decreasing in $d$.

We then specifically derive its rate in (\ref {th-random-rate}).
As $H^*(m) = 1$ implies $m = \big(n\mu_{\ba, \bb}/\Delta_{\ba, \bb}^d \big)^{2/(d+2)}$, we calculate the first integral in (\ref{th-random-main}) from $m=0$ to $m=\min\{ n\mu_{\ba, \bb}, \big(n\mu_{\ba, \bb}/\Delta_{\ba, \bb}^d \big)^{2/(d+2)}\} $.
The first term in (\ref {th-random-rate}) hence follows.
The last term follows from a straightforward calculation of (\ref {Rq+res}) using
\bes
\int_{\bx \in [{\bf 0}, \bb - \ba]}  \Big(n \prod_j x_j \vee 1\Big)^{-q/2} d\bx
\lesssim_{q,d,\rho_1, \rho_2} \begin{cases}
\big( \log_+(n\mu_{\ba, \bb})\big)^{d-I\{q>2\}} \big/ n& q \ge 2,
\cr
(n\mu_{\ba, \bb})^{-q/2+1}\big/n & 1\le q <2.
\end{cases}
\ees
This completes the proof. 
$\hfill\square$

\medskip
{\bf A6.3. Proof of Theorem \ref{coro-random-worst}.} We omit the proof as it's a direct result of Theorem \ref{th-random}.

\medskip
{\bf A6.4. Proof of Theorem \ref{coro-random-piecewise}.} We omit the proof as it's similar to the proof of Theorem \ref{coro-fix-piecewise}.

\bigskip
{\bf A7. Proofs of the results in Subsection 3.7}

\medskip
{\bf A7.1. Proof of Theorem \ref{th-model-mis}.} 
We provide the proof only for the $\{\fhat^{(block)}_n, \fbar_n^*\}$ pair given by 
\eqref{block-mid} and \eqref{block-mis} as the proofs 
for the $\{\fhat^{(max-min)}_n, \fbar_n^{(max-min)}\}$ 
and $\{\fhat^{(min-max)}_n, \fbar_n^{(min-max)}\}$ pairs are nearly identical and slightly simpler.  
By the definitions of $m_{\bx}$ and $\bv_{\bx}$ in (\ref{k-m-mis}),
\bes
\fhat^{(block)}_n(\bx)
&\le& \frac{1}{2}\Big\{\max_{\bu\preceq \bx}\ \min_{\bv_{\bx}\preceq\bv,}
\Ybar_{[\bu,\bv]} + \min_{\bv_{\bx}\preceq\bv}\ \max_{\bu\preceq \bx}
\Ybar_{[\bu,\bv]}\Big\}
\cr &\le& \fbar^*_n(\bx) + r_{q,+}^{1/q}(m_{\bx})
+ \max_{\bv\succeq\bv_{\bx}}\bigg(\max_{\bu \preceq \bx}\sum_{\bx_i\in [\bu,\bv]} \frac{\veps_i}
{n_{\bu,\bv}}\bigg)_+,
\ees
for data points $\bx = \bx_i$.
Thus, by the definition of $r_{q,+}(m)$ in (\ref{r-plus-mis}),
\bes
\E \Big\{\fhat^{(block)}_n(\bx) - \fbar_n^*(\bx)\Big\}_+^q
\le 2^qr_{q,+}(m_{\bx})
\ees
as in the proof of Theorem \ref{th-1}. 
Similarly, we can have the inequality on the negative side. 
This gives the $\fbar_n^*$ version of \eqref{th-1-1} and (\ref {th-1-2}). 

It remains to prove (\ref{r-plus-mis}) holds with $r_{q,\pm}(m) = C_{q,d}\sigma^q m^{-q/2}$, 
as the counterparts to the rest of the proof in the proofs of Theorems \ref{th-lattice-general}, \ref{th-fix-worst}, 
\ref{coro-fix-piecewise} and \ref{th-fix-selection} are all based 
on \eqref{th-1-1} and (\ref {th-1-2}) with 
$r_{q,\pm}(m)$ of this form. 
To this end, we notice that for fixed $\bx$ and $\bv_{\bx}$ in the lattice design, 
the partial sum $\sum_{\bu\le \bx_i \le \bv_{\bx} \le \bv}\veps_i$ indexed by $\bu$ and $\bv$ 
is a martingale in each index $u_j$ or $v_k$ while holding other $2d-1$ indices fixed. 
Thus, similar to 
the proof of \eqref{pf2-1-2} in the proof of Proposition~\ref{prop-single-point} we can group 
$\bu$ and $\bv$ in blocks of sizes $2^{|\bi|}m$ and $2^{|\bj|}m$ and bound 
$\E \max_{\bv\succeq\bv_{\bx}}\big(\max_{\bu \preceq \bx}\sum_{\bx_i\in [\bu,\bv]} 
\veps_i/n_{\bu,\bv}\big)_+^q$ by 
\bes
\sum_{\bi\ge {\bf 0}, \bj\ge {\bf 0}}\frac{C_{q,d}\sigma^q}{((2^{|\bi|}\vee 2^{|\bj|})m)^{q/2}}
\lesssim_{q,d} \sum_{\bj\ge {\bf 0}}\frac{\sigma^q |\bj|^d}{(2^{|\bj|}m)^{q/2}}
\lesssim_{q,d} \frac{\sigma^q}{m^{q/2}}, 
\ees
with $m = n_{\bx,\bv_{\bx}}$. This completes the proof. 
$\hfill\square$

\end{document}